\documentclass[12pt]{article}
\usepackage{fullpage}
\usepackage{enumerate}
\usepackage{graphicx}

\usepackage{epic}
\usepackage{eepic}
\usepackage[toc,page]{appendix} 
\usepackage{amsmath,amssymb,mathrsfs}
\usepackage[francais]{babel}
\usepackage[T1]{fontenc}
\usepackage[latin1]{inputenc}
\usepackage{soul}
\usepackage{pst-eucl}
\usepackage[all, ps]{xy}
\newtheorem{lema}{Lemme}[section]
\newtheorem{lem}{Lemme}[subsection]

\newtheorem{theoart3}{Théorème}
\newtheorem{theo}{Théorème}[subsection]
\newtheorem{prop}{Proposition}[subsection]

\author{Raphaël Beuzart-Plessis}
\title{Endoscopie et conjecture locale raffinée de Gan-Gross-Prasad pour les groupes unitaires}
\begin{document}

\maketitle

\bigskip

\textbf{Introduction}

\vspace{4mm}

Soient $F$ un corps local non archimédien de caractéristique nulle et $E$ une extension quadratique de $F$. Considérons deux espaces hermitiens $(V,h)$ et $(V',h')$ de dimensions respectives $d$ et $d'$ et de groupes unitaires respectifs $G$ et $G'$. On suppose $d$ pair et $d'$ impair et que le plus petit des deux espaces hermitiens s'injecte dans le second. Alors, une fois une telle injection choisie, le plus petit des groupes $G$ et $G'$ est naturellement un sous-groupe du plus gros. Soient $\sigma$ et $\sigma'$ des représentations lisses irréductibles de $G(F)$ et $G'(F)$ respectivement. Gan, Gross et Prasad ont défini dans [GGP] une multiplicité $m(\sigma,\sigma')$. La définition de cette multiplicité est rappelée en 4.1. Indiquons seulement ici sa signification dans deux cas extrêmes :

\vspace{2mm}

\begin{itemize}
\renewcommand{\labelitemi}{$\bullet$}

\item Si $d'=d-1$ (resp. $d=d'-1$), alors $G'$ est un sous-groupe de $G$ (resp. $G$ est un sous-groupe de $G'$) et on a

\vspace{2mm}

\begin{center}
$m(\sigma,\sigma')=dim_{\mathbb{C}}\; Hom_{G'}(\sigma\otimes \sigma',1)$
\end{center}

\begin{center}
(resp. $m(\sigma,\sigma')=dim_{\mathbb{C}}\; Hom_{G}(\sigma\otimes \sigma',1)$)
\end{center}

\vspace{2mm}

\noindent où $Hom_{G'}(\sigma\otimes \sigma',1)$ (resp. $Hom_{G}(\sigma\otimes \sigma',1)$) désigne l'espace des formes linéaires $G'(F)$-invariantes (resp. $G(F)$-invariantes) sur $\sigma\otimes \sigma'$.

\item Si $d=0$, alors $m(\sigma,\sigma')$ est la dimension de l'espace des fonctionnelles de Whittaker sur $\sigma'$ (pour un certain sous-groupe unipotent maximal de $G'$ et un caractère générique sur les $F$-points de ce sous-groupe).
\end{itemize}

\vspace{2mm}

\noindent En général, la définition de $m(\sigma,\sigma')$ est un "mixte" entre ces deux cas extrêmes (cf 4.1). D'après [GGP] et [AGRS], on a toujours $m(\sigma,\sigma')\leqslant 1$. \\

Nous nous intéressons ici seulement au cas où les représentations $\sigma$ et $\sigma'$ sont tempérées. Rappelons la correpondance, toujours conjecturale, de Langlands pour les représentations irréductibles tempérées d'un groupe unitaire. Enonçons ces conjectures pour $G$, celles-ci étant analogues pour $G'$ . Notons $WD_E=W_E\times SL_2(\mathbb{C})$ le groupe de Weil-Deligne de $E$ et fixons un élément $t\in W_F\backslash W_E$ que l'on laisse agir par conjugaison sur $WD_E$ (avec action triviale sur le facteur $SL_2(\mathbb{C})$). Soit $\mathbf{\Phi}_{temp}(G)$ l'ensemble des classes de conjugaison d'homomorphismes $\varphi: WD_E\to GL_d(\mathbb{C})$ vérifiant des conditions classiques de semi-simplicité et d'algébricité ainsi que les conditions supplémentaires suivantes

\begin{enumerate}[(i)]
\item $\varphi$ est tempérée, c'est-à-dire que l'image de $W_E$ par $\varphi$ est relativement compacte;
\item $\varphi$ est $(-1)^{d+1}$-conjuguée-duale, cela signifie qu'il existe une forme bilinéaire non dégénérée $B:\mathbb{C}^d\times \mathbb{C}^d\to \mathbb{C}$ telle que

\begin{itemize}
\item $B(\varphi(\tau)w,\varphi(t\tau t^{-1})w')=B(w,w')$;
\item $B(w,\varphi(t^2)w')=(-1)^{d+1}B(w,w')$
\end{itemize}

\noindent pour tous $w,w'\in\mathbb{C}^d$, $\tau\in WD_E$.
\end{enumerate}

\noindent Cette dernière condition ne dépend pas du choix de $t$. Soit $\varphi\in \mathbf{\Phi}_{temp}(G)$, on peut alors décomposer $\varphi$ en somme directe de représentations irréductibles de $WD_E$. Notons $\varphi=\bigoplus_{i\in I} l_i\varphi_i$ cette décomposition avec les $\varphi_i$ deux à deux non isomorphe. Notons $I'$ le sous-ensemble des $i\in I$ tels que $\varphi_i$ soit $(-1)^{d+1}$-conjuguée-duale. Fixons une forme $B$ vérifiant la condition (ii) précédente. Soit $S_{\varphi}$ le sous-groupe de $GL_d(\mathbb{C})$ des éléments qui commutent à $\varphi$ et qui préservent $B$. Notons $\mathcal{S}_{\varphi}$ son groupe des composantes. On a alors un isomorphisme naturel entre $\mathcal{S}_{\varphi}$ et $\big(\mathbb{Z}/2\mathbb{Z}\big)^{I'}$, en particulier $\mathcal{S}_{\varphi}$ ne dépend pas des choix de $B$ et $t$. Notons $z_\varphi$ l'image de $-I\in S_{\varphi}$ dans $\mathcal{S}_{\varphi}$. Définissons $\mathcal{E}^G(\varphi)$ comme l'ensemble des caractères $\epsilon$ de $\mathcal{S}_\varphi$ vérifiant $\epsilon(z_\varphi)=\mu(G)$, où on a posé

\[\begin{aligned}
\mu(G)=\left\{
    \begin{array}{ll}
        1 & \mbox{si } G \mbox{ est quasi-déployé} \\
        -1 & \mbox{sinon} 
    \end{array}
\right.
\end{aligned}\]

Soit $Temp(G)$ l'ensemble des classes d'isomorphisme de représentations tempérées irréductibles de $G(F)$. On conjecture que $Temp(G)$ est réunion disjointe de $L$-paquets $\Pi^G(\varphi)$ indexés par les $\varphi\in \mathbf{\Phi}_{temp}(G)$. On conjecture aussi que $\Pi^G(\varphi)$ admet une paramétrisation par $\mathcal{E}^G(\varphi)$ que l'on note $\epsilon\mapsto \sigma(\varphi,\epsilon)$. On impose à cette paramétrisation des conditions de types endoscopiques. Il y a tout d'abord l'endoscopie usuelle entre $G$ et ses groupes endoscopiques elliptiques. Ces groupes endoscopiques sont de la forme $G_+\times G_-$ où $G_+$ et $G_-$ sont des groupes unitaires quasi-déployés. On considère aussi une endoscopie tordue. Plus précisément soient $G_+$ et $G_-$ des groupes unitaires quasi-déployés d'espaces hermitiens de dimensions $d_+$ et $d_-$. Posons $n=d_++d_-$ et $M_{n}=R_{E/F} GL_{n}$ (où $R_{E/F}$ désigne la restriction des scalaires à la Weil). Notons $\theta$ l'automorphisme de $M_n$ donné par $g\mapsto {}^t(g^c)^{-1}$ et $\widetilde{M}_n$ la composante non neutre du groupe non connexe $M_n\rtimes \{1,\theta\}$. Alors $G_+\times G_-$ est un groupe endoscopique du groupe tordu $\widetilde{M}_n$. Soient $\varphi_+\in\mathbf{\Phi}_{temp}(G_+)$, $\varphi_-\in\mathbf{\Phi}_{temp}(G_-)$ et $\Pi^{G_+}(\varphi_+)$, $\Pi^{G_-}(\varphi_-)$ les $L$-paquets correspondants de représentations tempérées de $G_+(F)$ et $G_-(F)$ respectivement. Posons $\varphi=\varphi_+\oplus\varphi_-$. La correspondance locale de Langlands pour les groupes linéaires, prouvée par Harris-Taylor et Henniart, associe à $\varphi$ une représentation tempérée irréductible $\pi(\varphi)$ de $M_n(F)$. On peut étendre cette représentation en une représentation unitaire $\pi^+(\varphi)$ de $M_n(F)\rtimes\{1,\theta\}$ et on note $\tilde{\pi}(\varphi)$ sa restriction à $\widetilde{M}_n(F)$. D'après Clozel ([C]), $\tilde{\pi}(\varphi)$ possède un caractère que l'on note $\Theta_{\tilde{\pi}(\varphi)}$. On impose alors à nos $L$-paquets les deux conditions suivantes

\vspace{2mm}

\begin{itemize}
\renewcommand{\labelitemi}{$\bullet$}

\item Les caractères $\displaystyle\Theta_{\varphi_\pm}=\sum_{\sigma_\pm \in\Pi^{G_\pm}(\varphi_\pm)} \Theta_{\sigma_\pm}$ sont stables;

\item Il existe un nombre complexe $c$ de module $1$ tel que $c\Theta_{\tilde{\pi}(\varphi)}$ soit un transfert de $\Theta_{\varphi_+}\times \Theta_{\varphi_-}$. 
\end{itemize}

\vspace{2mm}

\noindent Pour que cela ait un sens, il faut bien entendu fixer précisément les facteurs de transfert. On renvoie pour cela à 3.2. Les conjectures sont les mêmes pour $G'$. Il suffit de remplacer $G$ par $G'$ et $d$ par $d'$ dans les définitions sauf en ce qui concerne l'ensemble de caractères $\mathcal{E}^{G'}(\varphi')$: on le définit comme l'ensemble des caractères $\epsilon$ de $\mathcal{S}_{\varphi'}$ vérifiant $\epsilon(z_{\varphi'})=\mu(G)$. On renvoie à 7.2 pour l'énoncé précis des conjectures, notamment en ce qui concerne l'endoscopie classique. On a introduit dans nos conjectures des constantes non précisées. Cela entraîne que la correspondance de Langlands n'est pas uniquement déterminée par nos hypothèses. Plus précisément, la composition des $L$-paquets $\Pi^G(\varphi)$ est entièrement déterminée par nos hypothèses mais pas les bijections $\epsilon\mapsto \sigma(\varphi,\epsilon)$ qui elles ne sont déterminées qu'à une "translation" près. On montre en 7.7 et 8.3 comment fixé précisément la correspondance. Dans le cas quasi-déployé, on demande que la représentation correspondant au caractère trivial soit l'unique représentation du $L$-paquet admettant un certain type de modèle de Whittaker. Dans le cas non quasi-déployé, on demande à notre paramétrisation de vérifier une instance de la conjecture de Gan-Gross-Prasad. Ces choix dépendent d'un caractère additif non trivial $\psi_E^\delta$ de $E$ qui est trivial sur $F$ (qui nous permet par exemple de choisir un type de modèle de Whittaker en rang pair). Dès lors, les $L$-paquets et les paramétrisations sont uniquement déterminés. \\

Soient $\varphi\in\mathbf{\Phi}_{temp}(G)$ et $\varphi'\in\mathbf{\Phi}_{temp}(G')$ des paramètres de Langlands tempérés pour $G$ et $G'$ respectivement. Soit $\varphi=\bigoplus_{i\in I} l_i\varphi_i$ la décomposition de $\varphi$ en composantes isotypiques. Comme indiqué précédemment, on a une identification naturelle $\mathcal{S}_\varphi=\big(\mathbb{Z}/2\mathbb{Z}\big)^{I'}$. On définit un caractère $\epsilon_{\varphi,\varphi'}^G$ de $\mathcal{S}_{\varphi}$ par

$$\displaystyle\epsilon^G_{\varphi,\varphi'}((e_i)_{i\in I'})=\prod_{i\in I'} \epsilon(1/2,\varphi_i\otimes \varphi',\psi_E^\delta)^{e_i}$$

\noindent pour tout $(e_i)_{i\in I'}\in \mathcal{S}_\varphi=\big(\mathbb{Z}/2\mathbb{Z}\big)^{I'}$ et où $\psi_E^\delta$ est le caractère additif qui nous a permis de fixer la correspondance de Langlands. Le facteur epsilon d'une représentation du groupe de Weil-Deligne est défini comme dans [GGP]. On définit de la même façon un caractère $\epsilon^{G'}_{\varphi,\varphi'}$ de $\mathcal{S}_{\varphi'}$. On vérifie que

$$\epsilon^G_{\varphi,\varphi'}(z_\varphi)=\epsilon^{G'}_{\varphi,\varphi'}(z_{\varphi'})=\epsilon(1/2,\varphi\otimes\varphi',\psi_E^\delta)$$

En admettant les conjectures ci-dessus (précisées en 7.2), notre résultat principal est le suivant (théorème 8.5.1)

\begin{theoart3}
Soient $\varphi\in\mathbf{\Phi}_{temp}(G)$ et $\varphi'\in\mathbf{\Phi}_{temp}(G')$. Alors

\begin{enumerate}[(i)]
\item Si $\mu(G)\neq\epsilon(1/2,\varphi\otimes\varphi',\psi_E^\delta)$, on a $m(\sigma,\sigma')=0$ pour tous $\sigma\in\Pi^G(\varphi)$, $\sigma'\in \Pi^{G'}(\varphi')$;
\item Si $\mu(G)=\epsilon(1/2,\varphi\otimes\varphi',\psi_E^\delta)$, on a

$$m\big(\sigma(\varphi,\epsilon_{\varphi,\varphi'}^G),\sigma(\varphi',\epsilon_{\varphi,\varphi'}^{G'})\big)=1$$

\noindent et $m(\sigma,\sigma')=0$ pour tous $\sigma\in \Pi^G(\varphi)$, $\sigma'\in \Pi^{G'}(\varphi')$ avec $(\sigma,\sigma')\neq \big(\sigma(\varphi,\epsilon_{\varphi,\varphi'}^G),\sigma(\varphi',\epsilon_{\varphi,\varphi'}^{G'})\big)$.
\end{enumerate}
\end{theoart3}

Il s'agit de la conjecture 17.3 de [GGP] pour les modèles de Bessel des groupes unitaires et restreinte aux représentations tempérées. Expliquons brièvement comment on prouve le théorème 1. Rappelons que l'on a montré en [B1] une formule intégrale calculant la multiplicité $m(\sigma,\sigma')$ en fonction des caractères de $\sigma$ et $\sigma'$. Par endoscopie classique on peut exprimer ces caractères en fonction de caractères stables vivant sur des groupes endoscopiques elliptiques de $G$ et $G'$. Par endoscopie tordue cette fois, on peut exprimer ces caractères stables en fonction de caractères sur les groupes tordus $\widetilde{M}_n$. Dans [B2], on a montré une formule intégrale pour certains facteurs epsilon de paires, similaire à celle pour la multiplicité, mais où on remplace les caractères de $\sigma$ et $\sigma'$ par des caractères sur des groupes tordus $\widetilde{M}_n$ et $\widetilde{M}_m$ avec $n$ et $m$ de parités différentes. Finalement, on obtient une expression de $m(\sigma,\sigma')$ à partir de facteurs epsilon de paires. C'est de cette expression dont découle le théorème 1. \\

 Mentionnons pour finir un autre résultat, le lemme 7.7.1, qui possède un intérêt propre. On peut l'énoncer ainsi : si $G$ est quasi-déployé, pour chaque type de modèle de Whittaker sur $G(F)$, il existe dans le $L$-paquet $\Pi^G(\varphi)$ une unique représentation admettant un tel modèle de Whittaker. Il en va de même pour $G'$ (il n'existe qu'un seul type de modèle de Whittaker pour $G'$). Signalons que cette propriété des $L$-paquets tempérés a déjà été obtenue pour les groupes classiques par T.Konno ([K]).

Décrivons le plan de l'article. Dans la première section, on a rassemblé les définitions, notations et résultats généraux sur les groupes, groupes tordus et groupes unitaires dont nous aurons besoin dans la suite. On paramètre dans la deuxième section les classes de conjugaison et conjugaison stable des groupes unitaires et des groupes tordus $\widetilde{M}_n$. Dans la troisième partie, on fixe les facteurs de transfert et on les exprime en fonction des paramétrages précédents. Les formules prouvées dans [B1] et [B2] sont réexprimées dans la quatrième section sous des formes plus appropriées pour la suite. Dans la cinquième partie, on prouve, grâce à un résultat de Rodier, une formule donnant le nombre de modèles de Whittaker d'un type fixé d'une représentation tempérée irréductible d'un groupe unitaire de rang pair. Dans la dernière section de cette partie (la section 5.6), on démontre une formule de transfert pour cette multiplicité de Whittaker. Dans la sixième partie, on prouve des formules de transfert analogues pour la multiplicité $m(\sigma,\sigma')$ et pour certains facteurs epsilon de paires. On y utilise de façon cruciale les principaux résultats de [B1] et [B2]. On précise nos hypothèses de travail sur la correspondance de Langlands dans la septième partie et on en déduit un certains nombre de conséquences. C'est dans la huitième est dernière partie que l'on démontre le théorème 1. Enfin, l'article contient une appendice où on démontre une formule intégrale pour certains facteurs epsilon de caractères. La proposition de l'appendice est utilisée dans la preuve de la proposition 7.6.1. \\

Cet article s'inspire bien évidemment grandement de la preuve par Jean-Loup Waldspurger de la conjecture de Gross-Prasad pour les groupes spéciaux orthogonaux ([W4]). Il apparaîtrait même évident à qui compare les deux articles que les structures des preuves sont tout à fait semblables. Il s'agit en fait d'une partie de la thèse de l'auteur sous la direction de Jean-Loup Waldspurger. Je voudrais exprimer ici ma profonde gratitude à l'égard de ce dernier autant pour avoir partager ses idées si fécondes que pour sa patience à répondre à mes questions. Je tiens aussi à remercier Wee Teck Gan, Atsushi Ichino ainsi que le rapporteur pour m'avoir fait remarquer qu'il y avait une erreur de signe dans le lemme 5.3.1 (entraînant une contradiction entre le choix des paramétrisations en 7.7 et les prédictions initiales de Gan, Gross et Prasad).

\section{Groupes, mesures, notations}

\subsection{Notations générales}

Soit $F$ un corps local non archimédien de caractéristique nulle et $E$ une extension quadratique de $F$. On notera $\mathcal{O}_F$, $\mathfrak{p}_F$, $\varpi_F$, $k_F$, $q_F$, $val_F$ et $|.|_F$ l'anneau des entiers de $F$, son idéal maximal, une uniformisante, le corps résiduel, le cardinal de $k_F$, la valuation normalisée par $val_F(\varpi_F)=1$ et la valeur absolue définie par $|x|_F=q_F^{-val_F(x)}$ pour tout $x\in F$. On définit de même $\mathcal{O}_E$, $\mathfrak{p}_E$, $\varpi_E$, $k_E$, $q_E$, $val_E$ et $|.|_E$. On notera $N$ et $Tr_{E/F}$ les applications norme et trace de l'extension $E/F$, $x\mapsto \overline{x}$ l'unique $F$-automorphisme non trivial de $E$ que l'on notera aussi $\tau_{E/F}$ et $sgn_{E/F}$ l'unique caractère quadratique de $F^\times$ dont le noyau est $N(E^\times)$. On fixe pour toute la suite un caractère continu non trivial $\psi$ de $F$, une clôture algébrique $\overline{F}$ de $F$ et un élément $\delta\in E$ non nul de trace nulle. On pose $\psi_E=\psi\circ Tr_{E/F}$. Pour tout groupe algébrique $G$ défini sur $E$, on notera $R_{E/F}G$ le groupe algébrique sur $F$ obtenu par restriction des scalaires à la Weil à partir de $G$. \\

\subsection{Groupes et groupes tordus}

Soit $G$ un groupe réductif défini sur $F$. On note $\mathfrak{g}$ son algèbre de Lie et on pose $\delta(G)=dim(G)-rg(G)$ où $rg(G)$ désigne le rang (absolu) de $G$. On désignera par $G_{ss}$ la sous-variété des éléments semi-simples de $G$. L'action adjointe de $G$ sur $\mathfrak{g}$ sera notée ainsi

$$G\times \mathfrak{g}\to \mathfrak{g}$$
$$(g,X)\mapsto ad(g)X=gXg^{-1}$$

Pour $x\in G$, on notera $Z_G(x)$, respectivement $G_x$, le centralisateur de $x$ dans $G$, respectivement le centralisateur connexe de $x$ dans $G$. Un élément $x\in G$ sera dit fortement régulier si $Z_G(x)$ est un sous-tore (forcément maximal) de $G$. On définit de la même façon la notion d'élément fortement régulier de $\mathfrak{g}$. On note $G_{reg}$, resp. $\mathfrak{g}_{reg}$, la sous-variété des éléments fortement réguliers de $G$, resp. $\mathfrak{g}$. On dira que deux éléments $x,y\in G_{reg}(F)$ sont stablement conjugués, s'ils sont conjugués par un élément de $G(\overline{F})$. On désigne par $G_{reg}(F)/conj$, resp. $G_{reg}(F)/stconj$, l'ensemble des classes de conjugaison, resp. classes de conjugaison stable, dans $G_{reg}(F)$. On note les projections naturelles $G_{reg}(F)\to G_{reg}(F)/conj$ et $G_{reg}(F)\to G_{reg}(F)/stconj$ par $\phi_G$, resp. $\phi^{st}_G$. On dispose d'une application naturelle $\phi^{st}_G/conj$ entre  $G_{reg}(F)/conj$ et $G_{reg}(F)/stconj$ qui rend le diagramme suivant commutatif

$$
\xymatrix{
G_{reg}(F)/conj \ar[rr]^{\phi^{st}_G/conj} & & G_{reg}(F)/stconj \\
 & G_{reg}(F) \ar[ur]^{\phi^{st}_G} \ar[ul]^{\phi_G}
}
$$

Soit $T$ un tore. On notera $A_T$ son sous-tore déployé maximal. On munit $T(F)$ d'une mesure de Haar de la façon suivante: si $T(F)$ est compact, on choisit la mesure de masse totale $1$; si $T$ est déployé, on choisit la mesure qui donne au sous-groupe compact maximal de $T(F)$ la mesure $1$; dans le cas général on choisit la mesure sur $T(F)$ compatible avec les mesures sur $A_T(F)$ et $T(F)/A_T(F)=(T/A_T)(F)$ que l'on vient de fixer. \\

On peut munir $G_{reg}(F)/conj$ d'une structure d'espace $F$-analytique et d'une mesure caractérisées par la condition suivante: pour tout $x\in G_{reg}(F)$ il existe un voisinage $\omega$ de $1$ dans $Z_G(x)$ tel que l'application

\begin{center}
$\omega\to G_{reg}(F)/conj$ \\
$t\mapsto \phi_G(xt)$
\end{center}

\noindent soit un isomorphisme qui préserve les mesures de $\omega$ sur un voisinage ouvert de $\phi_G(x)$. On munit de même $G_{reg}(F)/stconj$ d'une structure de variété $F$-analytique et d'une mesure. L'application $\phi^{st}_G/conj: G_{reg}(F)/conj\to G_{reg}(F)/stconj$ est alors un revêtement analytique qui préserve localement les mesures. Un élément $x\in G_{reg}(F)$ sera dit anisotrope si $G_x$ est un tore anisotrope. On note $G(F)_{ani}$ le sous-ensemble des éléments (fortement réguliers) anisotropes de $G(F)$. On vérifie que la propriété d'être anisotrope est invariante par conjugaison et conjugaison stable. On notera $G(F)_{ani}/conj$, resp. $G(F)_{ani}/stconj$, l'ensemble des classes de conjugaison, resp. conjugaison stable, dans $G(F)_{ani}$. Ces espaces s'identifient naturellement à des ouverts de $G_{reg}(F)/conj$ et $G_{reg}(F)/stconj$ respectivement. \\

Pour $x\in G(F)$ un élément semi-simple, on pose

$$D^G(x)=|det\big((1-ad(x))_{|\mathfrak{g}/\mathfrak{g}_x}\big)|_F$$

On désignera par $Nil(\mathfrak{g}(F))$ l'ensemble des orbites nilpotentes de $\mathfrak{g}(F)$. Fixons une forme bilinéaire non dégénérée $<.,.>$ sur $\mathfrak{g}(F)$ invariante par l'action adjointe de $G(F)$. Pour $f\in C_c^\infty(\mathfrak{g}(F))$, on pose

$$\hat{f}(X)=\displaystyle\int_{\mathfrak{g}} f(Y)\psi(<Y,X>) dY, \;\;\; \forall X\in\mathfrak{g}(F)$$

\noindent où $dY$ est la mesure de Haar sur $\mathfrak{g}(F)$ telle que $\hat{\hat{f}}(X)=f(-X)$. Pour $\mathcal{O}\in Nil(\mathfrak{g}(F))$, on munit $\mathcal{O}$ d'une mesure de la façon suivante. En tout point $X\in \mathcal{O}$, la forme bilinéaire $(Y,Z)\mapsto <X,[Y,Z]>$ sur $\mathfrak{g}(F)$ se descend en une forme symplectique sur $\mathfrak{g}(F)/\mathfrak{g}_X(F)$, c'est-à-dire sur l'espace tangent à $\mathcal{O}$ au point $X$. Ainsi $\mathcal{O}$ est muni d'une structure de variété $F$-analytique symplectique et on en déduit une mesure "autoduale" sur $\mathcal{O}$ via le caractère $\psi$. On définit l'intégrale orbitale sur $\mathcal{O}$ par

$$\displaystyle J_{\mathcal{O}}(f)=\int_{\mathcal{O}} f(X) dX$$

\noindent pour tout $f\in C_c^\infty(\mathfrak{g}(F))$. D'après Harish-Chandra, [HCvD], il existe une fonction $\hat{j}(\mathcal{O},.)$ localement intégrable sur $\mathfrak{g}(F)$ telle que

$$\displaystyle J_{\mathcal{O}}(\hat{f})=\int_{\mathfrak{g}(F)} \hat{j}(\mathcal{O},X) f(X) dX$$

\noindent pour tout $f\in C_c^\infty(\mathfrak{g}(F))$. Lorsque l'on voudra préciser le groupe ambiant, on notera $\hat{j}^G(\mathcal{O},.)$ plutôt que $\hat{j}(\mathcal{O},.)$. Les fonctions $\hat{j}(\mathcal{O},.)$ vérifient la propriété d'homogénéité suivante:

$$\hat{j}(\mathcal{O},\lambda X)=|\lambda|_F^{-dim(\mathcal{O})/2}\hat{j}(\lambda \mathcal{O},X)$$

\noindent pour tout $\mathcal{O}\in Nil(\mathfrak{g}(F))$, pour tout $X\in \mathfrak{g}_{reg}(F)$ et pour tout $\lambda\in F^\times$. Signalons que si $\lambda\in F^{\times, 2}$ alors $\lambda\mathcal{O}=\mathcal{O}$. On notera $Nil(\mathfrak{g}(F))_{reg}$ le sous-ensemble de $Nil(\mathfrak{g}(F))$ constitué des orbites régulières c'est-à-dire des orbites qui sont de dimension $\delta(G)$. Ce sous-ensemble est non vide si et seulement si $G$ est quasi-déployé sur $F$. \\

Enfin on notera $Temp(G)$ l'ensemble des classes d'équivalence de représentations (complexes) irréductibles lisses et tempérées de $G(F)$. Pour $\sigma\in Temp(G)$, on notera $E_\sigma$ un espace vectoriel complexe sur lequel se réalise la représentation $\sigma$, $\sigma^\vee$ la représentation contragrédiente de $\sigma$ et $\Theta_\sigma$ le caractère de Harish-Chandra de $\sigma$ (vu comme une fonction localement intégrable sur le groupe). \\ 

Soit $(M,\widetilde{M})$ un groupe tordu défini sur $F$. Cela signifie que $M$ est un groupe réductif connexe défini sur $F$ et que $\widetilde{M}$ est une variété définie sur $F$, vérifiant $\widetilde{M}(F)\neq\emptyset$, et munie de deux actions commutantes à gauche et à droite de $M$ qui font chacune de $\widetilde{M}$ un espace principal homogène sous $M$. On note

$$M\times \widetilde{M}\times M\to \widetilde{M}$$
$$(m,\tilde{m},m')\to m\tilde{m}m'$$

\noindent les actions à gauche et à droite de $M$ sur $\widetilde{M}$. Pour $\tilde{m}\in \widetilde{M}$, on note $\theta_{\tilde{m}}$ l'automorphisme de $M$ tel que $\tilde{m}m=\theta_{\tilde{m}}(m)\tilde{m}$ pour tout $m\in M$. On note $Z_M(\tilde{m})$ le centralisateur de $\tilde{m}$ dans $M$ pour l'action par conjugaison, c'est-à-dire le sous-groupe des points fixes de $\theta_{\tilde{m}}$, et $M_{\tilde{m}}$ la composante neutre de $Z_M(\tilde{m})$. On dira que $\tilde{m}$ est fortement régulier si $Z_M(\tilde{m})$ est abélien et $M_{\tilde{m}}$ est un tore. On note $\widetilde{M}_{reg}$ la sous-variété des éléments fortement réguliers et $\widetilde{M}_{reg}(F)/conj$ l'ensemble des classes de conjugaison d'éléments fortement réguliers de $\widetilde{M}(F)$. On dit que deux éléments de $\widetilde{M}_{reg}(F)$ qu'ils sont stablement conjugués s'il sont conjugués par un élément de $M(\overline{F})$. On note $\widetilde{M}_{reg}(F)/stconj$ l'ensemble des classes de conjugaison stable de $\widetilde{M}_{reg}(F)$. On dispose, comme dans le cas des groupes, d'un diagramme commutatif naturel

$$
\xymatrix{
\widetilde{M}_{reg}(F)/conj \ar[rr]^{\phi^{st}_{\tilde{M}}/conj} & & \widetilde{M}_{reg}(F)/stconj \\
 & \widetilde{M}_{reg}(F) \ar[ur]^{\phi^{st}_{\tilde{M}}} \ar[ul]^{\phi_{\tilde{M}}}
}
$$

On peut munir $\widetilde{M}_{reg}(F)/conj$ d'une structure d'espace $F$-analytique et d'une mesure caractérisées par la propriété suivante: pour tout $\tilde{m}\in \widetilde{M}_{reg}(F)$, il existe un voisinage ouvert $\omega$ de $1$ dans $M_{\tilde{m}}(F)$ de sorte que l'application $t\in\omega\mapsto \phi_{\tilde{M}}(t\tilde{m})$ soit un isomorphisme qui préserve les mesures de $\omega$ sur un ouvert de $\widetilde{M}_{reg}(F)/conj$. On munit de la même façon $\widetilde{M}_{reg}(F)/stconj$ d'une structure $F$-analytique et d'une mesure. L'application $\phi^{st}_{\widetilde{M}}/conj$ est alors un revêtement analytique qui préserve localement les mesures. L'action de $\theta_{\tilde{m}}$ sur $A_G$ est la même quel que soit $\tilde{m}\in \widetilde{M}$. On note $A_{\tilde{G}}$ la composante neutre du sous-groupe des points fixes de $A_G$ pour cette action. Un élément $\tilde{x}\in \widetilde{M}_{reg}(F)$ sera dit anisotrope si $M_{\tilde{x}}$ est un tore anisotrope. On notera $\widetilde{M}(F)_{ani}$ le sous-ensemble des éléments (fortement réguliers) anisotropes de $\widetilde{M}(F)$. On vérifie que la propriété d'être anisotrope est invariante par conjugaison et conjugaison stable. On note $\widetilde{M}(F)_{ani}/conj$, resp. $\widetilde{M}(F)_{ani}/stconj$, l'ensemble des classes de conjugaison, resp. conjugaison stable, dans $\widetilde{M}(F)_{ani}$. Ces espaces s'identifient naturellement à des ouverts de $\widetilde{M}_{reg}(F)/conj$ et $\widetilde{M}_{reg}(F)/stconj$ respectivement. \\

Soit $\tilde{x}\in \widetilde{M}$. On dit que $\tilde{x}$ qu'il est semi-simple s'il existe une paire $(B,T)$ formée d'un sous-groupe de Borel $B$ de $M$ et d'un tore maximal $T$ de $B$ (on ne les suppose pas définis sur $F$) qui est stable par $\theta_{\tilde{x}}$. On notera $\widetilde{M}_{ss}$ la sous-variété des éléments semi-simple. Pour $\tilde{x}\in\widetilde{M}_{ss}(F)$, on pose

$$D^{\tilde{M}}(\tilde{x})=|det\; (1-\theta_{\tilde{x}})_{|\mathfrak{m}/\mathfrak{m}_{\tilde{x}}}|_F$$

Supposons de plus $\tilde{x}$ fortement régulier. Alors $M_{\tilde{x}}$ est un tore et son centralisateur dans $M$ est un tore maximal $T$. On pose alors

$$D_0^{\tilde{M}}(\tilde{x})=|det(1-\theta_{\tilde{x}})_{|\mathfrak{m}/\mathfrak{t}}|_F$$

\subsection{Les groupes unitaires et leurs orbites nilpotentes régulières}

Un espace hermitien sera pour nous un couple $(V,h)$ formé d'un espace vectoriel de dimension finie sur $E$ et d'une forme hermitienne non dégénérée $h$ sur $V$. Nos formes hermitiennes seront par convention linéaire en la deuxième variable. A un tel couple est associé un groupe unitaire $G=U(V,h)$. C'est un groupe réductif connexe défini sur $F$. Soit $(V,h)$ un espace hermitien de dimension $d$ et de groupe unitaire $G$. L'algèbre de Lie $\mathfrak{g}(F)$ s'identifie alors naturellement au sous-espace des endomorphismes $E$-linéaires $X\in End_E(V)$ qui vérifient

$$h(Xv,v')+h(v,Xv')=0$$

\noindent pour tous $v,v'\in V$. Pour tous $X,X'\in \mathfrak{g}(F)$, on pose

$$\displaystyle <X,X'>=\frac{1}{2}Tr_{E/F}(Tr(XX'))$$

\noindent On vérifie que $<.,.>$ est une forme bilinéaire non dégénérée $G$-invariante. On peut alors appliquer les constructions du paragraphe précédent associées à un tel choix. \\

 Si $d$ est impair ou $d=0$ alors $G$ est quasi-déployé et $\mathfrak{g}(F)$ ne possède qu'une seule orbite nilpotente régulière. Si $d\geqslant 2$ est pair et que $G$ n'est pas quasi-déployé, alors $\mathfrak{g}(F)$ ne possède aucune orbite nilpotente régulière. Enfin si $d\geqslant 2$ est pair et que $G$ est quasi-déployé, il existe deux orbites nilpotentes régulières dans $\mathfrak{g}(F)$. On peut paramétrer ces orbites par $\left(Ker\; Tr_{E/F}\backslash\{0\}\right)/N(E^\times)$ de la façon suivante: on fixe une base $(e_i)_{i=\pm 1,\ldots,\pm d/2}$ de $V$ vérifiant 

$$h(e_i,e_j)=\delta_{i,-j}$$

\noindent pour $i,j=\pm 1,\ldots,\pm d/2$ et où $\delta_{i,-j}$ est le symbole de Kronecker. Soit $\eta\in Ker\; Tr_{E/F}$ non nul. On définit un élément $N_\eta\in\mathfrak{g}(F)$ par $N_\eta e_{d/2}=0$, $N_\eta e_i=e_{i+1}$ pour $i=1,\ldots,d/2-1$, $N_\eta e_i=-e_{i+1}$ pour $i=-d/2,\ldots,-2$ et $N_\eta e_{-1}=\eta e_1$. Alors $N_\eta$ est un élément nilpotent régulier et l'orbite qu'il engendre ne dépend que de $\eta$ modulo $N(E^\times)$ (et pas du choix de la base $(e_i)$). On note $\mathcal{O}_\eta$ cette orbite. On obtient ainsi une bijection entre $\big(Ker\; Tr_{E/F}\backslash\{0\}\big)/N(E^\times)$ et $Nil(\mathfrak{g}(F))_{reg}$. \\

Supposons que $G$ est quasi-déployé. Soit $B$ un sous-groupe de Borel et $T_{qd}\subset B$ un tore maximal, tous deux définis sur $F$. Notons $W^G=W(G,T_{qd})=Norm_{G(F)}(T_{qd})/T_{qd}(F)$ le groupe de Weyl de $G$. On a alors $|W^G|=w(d)$, où on a posé pour tout entier $d\geqslant 0$,

\[\begin{aligned}
\mbox{(1)}\;\;\; w(d)=\left\{
    \begin{array}{ll}
        2^{d/2}(d/2)! & \mbox{si } d \mbox{ est pair} \\
        2^{(d-1)/2}\big((d-1)/2\big)! & \mbox{si } d \mbox{ est impair} 
    \end{array}
\right.
\end{aligned}\]

Pour tout entier $r\geqslant 0$, on fixe un $E$-espace vectoriel $Z_{2r+1}$ de dimension $2r+1$ ainsi qu'une base $(z_i)_{i=0,\pm 1,\ldots,\pm r}$ de celui-ci. Pour tout $\nu\in F^\times$, on notera alors $h_{2r+1,\nu}$ la forme hermitienne sur $Z_{2r+1}$ définie par

$$\displaystyle h\left(\lambda_0z_0+\sum_{i=1,\ldots ,r} \lambda_iz_i+\lambda_{-i}z_{-i},\lambda'_0z_0+\sum_{i=1,\ldots ,r} \lambda'_iz_i+\lambda'_{-i}z_{-i}\right)=\nu\overline{\lambda}_0\lambda'_0+\sum_{i=1}^r \overline{\lambda}_i\lambda'_{-i}+\overline{\lambda}_{-i}\lambda'_i$$

\noindent pour tous $(\lambda_i)_{i=0,\pm 1,\ldots,\pm r}$, $(\lambda'_i)_{i=0,\pm 1,\ldots,\pm r}\in E^{2r+1}$.

\section{Paramètrages}

\subsection{Espaces de paramètres}

Dans ce qui suit, on entendra par extension de $F$ un sous-corps de $\overline{F}$ qui contient $F$. On définit $\underline{\Xi}$ comme l'ensemble des quadruplets $\xi=(I,(F_{\pm i})_{i\in I}, (F_i)_{i\in I}, (y_i)_{i\in I})$ où

\begin{itemize}
\item $I$ est un ensemble fini;
\item Pour tout $i\in I$, $F_{\pm i}$ est une extension finie de $F$ et $F_i=F_{\pm i}\otimes_F E$. On note $\tau_i=Id\otimes \tau_{E/F}$ l'unique $F_{\pm i}$-automorphisme non trivial de $F_i$;
\item Pour tout $i\in I$, $y_i$ est un élément de $F_i^\times$ vérifiant $y_i\tau_i(y_i)=1$.
\end{itemize}

Pour un tel $\xi$, on notera $I^*$ l'ensemble des $i\in I$ tels que $F_i$ soit un corps et on posera

$$\displaystyle d_\xi=\sum_{i\in I} [F_{\pm i}:F]=\sum_{i\in I} [F_i:E]$$

Pour $d\in\mathbb{N}$, on notera $\underline{\Xi}_d$ l'ensemble des éléments $\xi\in\underline{\Xi}$ vérifiant $d_\xi=d$, $\underline{\Xi}^*$ l'ensemble des éléments $\xi\in\underline{\Xi}$ vérifiant $I^*=I$ et $\underline{\Xi}^*_d=\underline{\Xi}_d\cap \underline{\Xi}^*$. On définit un isomorphisme entre deux éléments $\xi=(I,(F_{\pm i})_{i\in I}, (F_i)_{i\in I}, (y_i)_{i\in I})$ et $\xi'=(I',(F'_{\pm i'})_{i'\in I'}, (F'_{i'})_{i'\in I'},(y'_{i'})_{i'\in I'})$ de $\underline{\Xi}$ comme étant un triplet $(\iota,(\iota_{\pm i})_{i\in I}, (\iota_i)_{i\in I})$ où

\begin{itemize}
\item $\iota: I\to I'$ est une bijection;
\item Pour tout $i\in I$, $\iota_{\pm i}: F_{\pm i}\to F'_{\pm \iota(i)}$ est un $F$-isomorphisme et $\iota_i=\iota_{\pm i}\otimes Id_E:F_i\to F'_{\iota(i)}$ est le $E$-isomorphisme déduit.
\end{itemize}

Ce triplet devant vérifier la condition supplémentaire

$$\iota_i(y_i)=y'_{\iota(i)},\; \forall i\in I$$

Un élément $\xi\in \underline{\Xi}$ sera dit régulier si l'identité est son seul automorphisme. On note $\underline{\Xi}_{reg}$ le sous-ensemble des éléments réguliers de $\underline{\Xi}$. On définit $\Xi$, resp. $\Xi_{reg}$, resp. $\Xi_d$, resp. $\Xi^*$ comme l'ensemble des classes d'isomorphisme de $\underline{\Xi}$, resp. $\underline{\Xi}_{reg}$, resp. $\underline{\Xi}_d$, resp. $\underline{\Xi}^*$. On posera $\Xi_{reg,d}=\Xi_d\cap \Xi_{reg}$, $\Xi_d^*=\Xi_d\cap\Xi^*$, $\Xi_{reg,d}^*=\Xi_{d,reg}\cap\Xi^*$. On identifiera toujours une classe d'isomorphisme de paramètres à un élément qui la représente. \\

Pour $\xi=(I,(F_{\pm i})_{i\in I}, (F_i)_{i\in I}, (y_i)_{i\in I})\in \Xi_{reg}$, introduisons le tore

$$\displaystyle T_\xi=\prod_{i\in I} Ker \; N_{F_i/F_{\pm i}}$$

\noindent C'est un tore défini sur $F$. L'espace $\Xi_{reg}$ peut être muni d'une structure de variété analytique et d'une mesure caractérisées de la façon suivante: pour tout $\xi=(I,(F_{\pm i})_{i\in I}, (F_i)_{i\in I}, (y_i)_{i\in I})\in \Xi_{reg}$, il existe un voisinage ouvert $\omega$ de $1$ dans $T_\xi(F)$ tel que l'application $\omega\to \Xi_{reg}$

$$(t_i)_{i\in I}\mapsto (I,(F_{\pm i})_{i\in I}, (F_i)_{i\in I}, (y_it_i)_{i\in I})$$

\noindent induise un isomorphisme qui préserve les mesures de $\omega$ sur un ouvert de $\Xi_{reg}$. Soient $\xi_+=(I_+,(F_{\pm i})_{i\in I_+},(F_i)_{i\in I_+},(y_i)_{i\in I_+})$, $\xi_-=(I_-,(F_{\pm i})_{i\in I_-},(F_i)_{i\in I_-},(y_i)_{i\in I_-})\in \Xi$. On définit alors $\xi_+\sqcup \xi_-\in \Xi$ comme étant égal à $(I,(F_{\pm i})_{i\in I},(F_i)_{i\in I},(y_i)_{i\in I})$ où $I=I_+\sqcup I_-$. Si $\xi_+\sqcup\xi_-$ est régulier, il en va de même de $\xi_+$ et $\xi_-$. De plus, l'application $(\xi_+,\xi_-)\in \Xi^2\mapsto \xi_+\sqcup \xi_-\in \Xi$ préserve localement les mesures au voisinage des éléments de $\Xi_{reg}$. \\

Soit $\xi=(I,(F_{\pm i})_{i\in I}, (F_i)_{i\in I}, (y_i)_{i\in I})\in \Xi_{reg}$. On lui associe le groupe abélien fini

$$\displaystyle C(\xi)=\prod_{i\in I} F_{\pm i}^\times /N_{F_i/F_{\pm i}}(F_i^\times)$$

Ce groupe s'identifie naturellement à $\{\pm 1\}^{I^*}$. On notera $C(\xi)^1$ le sous-groupe des éléments dont le produit des coordonnées vaut $1$ et $C(\xi)^{-1}=C(\xi)\backslash C(\xi)^1$. Pour $i\in I$, notons $\Gamma(y_i)$ l'ensemble des éléments $\gamma_i\in F_i^\times$ tels que $\gamma_i\tau_i(\gamma_i)^{-1}=y_i$. On pose

$$\displaystyle \Gamma(\xi)=\prod_{i\in I} \Gamma(y_i)/N_{F_i/F_{\pm i}}(F_i^\times)$$

\noindent C'est un espace principal homogène sous $C(\xi)$.

\subsection{Paramétrages des classes de conjugaison}

Soit $(V,h)$ un espace hermitien de dimension $d$ et de groupe unitaire $G$. Soient $\xi=(I,(F_{\pm i})_{i\in I}, (F_i)_{i\in I}, (y_i)_{i\in I})\in \underline{\Xi}_{reg,d}$ et $c=(c_i)_{i\in I}\in \prod_{i\in I} F_{\pm i}^\times$. Introduisons l'espace hermitien $(V_{\xi,c},h_{\xi,c})$ ainsi défini

\vspace{3mm}

\begin{itemize}
\renewcommand{\labelitemi}{$\bullet$}

\item $V_{\xi,c}=\bigoplus_{i\in I} F_i$

\item $h_{\xi,c}(\sum_{i\in I} x_i,\sum_{i\in I} x'_i)=\sum_{i\in I} Trace_{F_i/E}(c_ix'_i\tau_i(x_i))$
\end{itemize}

\vspace{3mm}

La classe d'isomorphisme de $(V_{\xi,c},h_{\xi,c})$ ne dépend que des images de $\xi$ et $c$ dans $\Xi_{reg,d}$ et $C(\xi)$ respectivement. De plus, il existe un unique $\epsilon\in\{\pm 1\}$ de sorte que $(V_{\xi,c},h_{\xi,c})$ soit isomorphe à $(V,h)$ si et seulement si $c\in C(\xi)^\epsilon$. Pour un tel $c$, fixons un isomorphisme $(V,h)\simeq (V_{\xi,c},h_{\xi,c})$ et notons $x(\xi,c)$ l'élément de $G_{reg}(F)$ qui via cet isomorphisme agit sur $V_{\xi,c}$ comme la multiplication par $y_i$ sur $F_i\subset V_{\xi,c}$. La classe de conjugaison de $x(\xi,c)$ ne dépend que des images de $\xi$ et $c$ dans $\Xi_{reg,d}$ et $C(\xi)$ respectivement (donc notamment pas de l'isomorphisme choisi). Toutes les classes de conjugaison dans $G_{reg}(F)$ sont ainsi atteintes pour un unique $\xi\in \Xi_{reg,d}$ et un unique $c\in C(\xi)$. De plus, la classe de conjugaison stable de $x(\xi,c)$ ne dépend elle que de $\xi\in \Xi_{reg,d}$. On obtient ainsi un diagramme commutatif

$$
\xymatrix{
G_{reg}(F)/conj \ar[rd]^{p_G} \ar[rr]^{\phi^{st}_G/conj} & & G_{reg}(F)/stconj \ar[ld]^{p_G^{st}}\\
 & \Xi_{reg,d}
}
$$

\noindent où $p^{st}_G$ est un isomorphisme $F$-analytique qui préserve les mesures de $G_{reg}(F)/stconj$ sur un ouvert de $\Xi_{reg,d}$ et $p_G$ est un revêtement au dessus du même ouvert. La fibre en un point $\xi$ de ce revêtement s'identifie à une classe $C(\xi)^\epsilon\subset C(\xi)$. Une classe de conjugaison (resp. conjugaison stable) est anisotrope si et seulement si son image par $p_G$ (resp. $p^{st}_G$) appartient à $\Xi^*_{reg,d}$.\\

Soit $U$ un espace vectoriel sur $E$ de dimension $n$, $M=R_{E/F} GL(U)$ et $\widetilde{M}$ la variété algébrique définie sur $F$ des formes sesquilinéaires non dégénérées sur $U$. Tout comme les formes hermitiennes, on ne considère que des formes sesquilinéaires qui sont linéaires en la deuxième variable. On laisse agir $M$ à gauche et à droite sur $\widetilde{M}$ de la façon suivante

$$(m\tilde{m}m')(u,u')=\tilde{m}(m^{-1}u,m'u')$$

\noindent pour tout $\tilde{m}\in \widetilde{M}$ et pour tous $m,m'\in M$. Le couple $(M,\widetilde{M})$ est un groupe tordu. Soit $\xi=(I,(F_{\pm i})_{i\in I}, (F_i)_{i\in I}, (y_i)_{i\in I}) \in \underline{\Xi}_{reg,n}$ et fixons un isomorphisme $E$-linéaire

$$M\simeq \bigoplus_{i\in I} F_i$$

\noindent Pour $\gamma=(\gamma_i)\in \prod_{i\in I} \Gamma(y_i)$, introduisons l'élément $\tilde{x}(\xi,\gamma)\in \widetilde{M}_{reg}(F)$ qui via l'isomorphisme précédent est donné par la formule

$$\tilde{x}(\xi,\gamma)\left(\sum_{i\in I}u_i,\sum_{i\in I} u'_i\right)=\sum_{i\in I} Trace_{F_i/E}(\gamma_i\tau_i(u_i)u'_i)$$

\noindent Alors la classe de conjugaison de $\tilde{x}(\xi,\gamma)$ est bien définie et ne dépend que des images de $\xi$ et $\gamma$ dans $\Xi_{reg,n}$ et $\Gamma(\xi)$ respectivement. Toutes les classes de conjugaison dans $\widetilde{M}_{reg}(F)$ sont ainsi atteintes pour un unique $\xi\in \Xi_{reg,n}$ et un unique $\gamma\in \Gamma(\xi)$. De plus, la classe de conjugaison stable de $\tilde{x}(\xi,\gamma)$ ne dépend elle que de $\xi\in \Xi_{reg,n}$. On obtient ainsi un diagramme commutatif

$$
\xymatrix{
\widetilde{M}_{reg}(F)/conj \ar[rd]^{p_{\tilde{M}}} \ar[rr]^{\phi^{st}_{\tilde{M}}/conj} & & \widetilde{M}_{reg}(F)/stconj \ar[ld]^{p_{\tilde{M}}^{st}}\\
 & \Xi_{reg,n}
}
$$

\noindent où $p^{st}_{\tilde{M}}$ est un isomorphisme $F$-analytique de $\widetilde{M}_{reg}(F)/stconj$ sur $\Xi_{reg,n}$ et $p_{\tilde{M}}$ est un revêtement. La fibre en un point $\xi$ de ce revêtement s'identifie à $\Gamma(\xi)$. Par contre, l'application $p^{st}_{\tilde{M}}$ ne préserve pas les mesures et son jacobien vaut $|2|_F^n$. En effet, soient $\xi=(I,(F_{\pm i})_{i\in I},(F_i)_{i_\in I}, (y_i)_{i\in I})\in \Xi_{reg,n}$ et $\gamma\in \Gamma(\xi)$. Posons $\tilde{m}=\tilde{x}(\xi,\gamma)$ en ayant effectués les choix nécessaires. On a alors un isomorphisme $T_\xi\simeq M_{\tilde{m}}$ qui à $t=(t_i)_{i\in I}\in T_\xi$ associe l'élément de $M_{\tilde{m}}$ qui agit sur $F_i$ comme la multiplication par $t_i$. Modulo cet isomorphisme, pour tout $t=(t_i)_{i\in I}\in M_{\tilde{m}}(F)$, la classe de conjugaison stable de $\tilde{m}t$ est paramétrée par $(I,(F_{\pm i})_{i\in I},(F_i)_{i_\in I}, (y_it_i^2)_{i\in I})$. L'application $M_{\tilde{m}}\to \widetilde{M}_{reg}(F)/stconj$ , $t\mapsto \phi_{\tilde{M}}^{st}(\tilde{m}t)$ préserve localement les mesures au voisinage de l'origine, tandis que l'application

$$T_\xi(F)\to \Xi_{reg,n}$$
$$(t_i)\mapsto (I,(F_{\pm i})_{i\in I},(F_i)_{i_\in I}, (y_it_i^2)_{i\in I})$$

\noindent est de jacobien $|2|_F^n$ à l'origine. Ceci entraîne l'affirmation sur le jacobien de $p_{\tilde{M}}^{st}$. \\

Remarquons enfin qu'une classe de conjugaison (resp. conjugaison stable) est anisotrope si et seulement si son image par $p_{\tilde{M}}$ (resp. $p^{st}_{\tilde{M}}$) appartient à $\Xi^*_{reg,n}$.

\subsection{Définition de facteurs sur les paramètres}

On vient de voir comment paramétrer les classes de conjugaison et de conjugaison stable (suffisamment régulières) des groupes unitaires et de certains groupes tordus. Dans ce paragraphe, assez aride, on définit plusieurs fonctions sur les espaces de paramètres. Il y a tout d'abord certaines fonctions déterminants qui apparaîtront maintes fois dans la suite de l'article (la fonction $\Delta$ et les fonctions $D^d$). On définit aussi deux fonctions $\Delta_{\mu_+,\mu_-,\nu}$ et $\Delta_{\mu_+,\mu_-}$ qui s'avéreront plus tard être des facteurs de transfert (cf les paragraphes 3.1 et 3.2).

\vspace{2mm}

Soient $\xi=(I,(F_{\pm i})_{i\in I}, (F_i)_{i\in I}, (y_i)_{i\in I})\in \Xi_{reg}$. Pour tout $i\in I$, on note $\Phi_i$ l'ensemble des homomorphismes non nuls de $E$-algèbres de $F_i$ dans $\overline{F}$. Définissons un polynôme $P_\xi\in E[T]$ par

$$P_\xi(T)=\prod_{i\in I}\prod_{\phi\in \Phi_i} (T-\phi(y_i))$$

\noindent Posons

$$\Delta(\xi)=|P_\xi(1)|_E$$

\noindent Soit $(V_\xi,h_\xi)$ un espace hermitien de dimension $d_\xi$, de groupe unitaire $G_\xi$ et soit $t\in G_{\xi,reg}(F)$ dont la classe de conjugaison stable est paramétrée par $\xi$. Posons

$$\Delta(t)=|det(1-t)_{|V_\xi}|_E$$

\noindent On a alors

$$\Delta(t)=\Delta(\xi)$$

\noindent Soit $d\geqslant d_\xi$ un entier et $(Z,h_Z)$ un espace hermitien de dimension $d-d_\xi$. Soit $G$ le groupe unitaire de $(Z,h_Z)\oplus (V_\xi,h_\xi)$. Alors $t$ est un élément semi-simple de $G(F)$ et on pose

$$D^d(\xi)=D^G(t)$$

\noindent Cette fonction ne dépend que de $d$ et $\xi$ (donc pas des choix de $(V_\xi,h_\xi)$, $t$ et $(Z,h_Z)$). On pose $D(\xi)=D^{d_\xi}(\xi)$. On a alors la relation

$$D^d(\xi)=\Delta(\xi)^{d-d_\xi}D(\xi)$$

Soient $\xi_+=(I_+,(F_{\pm i})_{i\in I_+}, (F_i)_{i\in I_+}, (y_i)_{i\in I_+})\in \Xi_{reg}$, $\xi_-=(I_-,(F_{\pm i})_{i\in I_-}, (F_i)_{i\in I_-}, (y_i)_{i\in I_-})\in \Xi_{reg}$ et $\mu_+$, $\mu_-$ deux caractères continus de $E^{\times}$. Posons $\xi=\xi_+\sqcup \xi_-$, $I=I_+\sqcup I_-$, $d_-=d_{\xi_-}$, $d_+=d_{\xi_+}$ et $d=d_\xi$. \\

\noindent Pour $c=(c_i)_{i\in I}\in C(\xi)$ et $\nu\in F^\times$, on pose

$$\mbox{(1)}\;\;\; \Delta_{\mu_+,\mu_-,\nu}(\xi_+,\xi_-,c)=\mu_-(P_{\xi_-}(-1))\mu_+(P_{\xi_+}(-1))\prod_{i\in I_-} sgn_{F_i/F_{\pm i}}(\nu C_i)$$

\noindent où

$$C_i=\left\{
    \begin{array}{ll}
      -\delta^{-d-1} c_iP_\xi '(y_i)P_\xi (-1)^{-1} y_i^{1-d/2} & \mbox{ si } d \mbox{ est pair}\\
      -\delta^{-d-1} c_iP_\xi '(y_i)P_\xi (-1)^{-1} y_i^{(1-d)/2}(1+y_i) & \mbox{ si } d \mbox{ est impair}
    \end{array}
\right.$$

\noindent On verra en 3.1 que les fonctions $\Delta_{\mu_+,\mu_-,\nu}$ correspondent à certains facteurs de transfert pour les groupes unitaires.

\vspace{2mm}

\noindent Pour $\gamma=(\gamma_i)_{i\in I}\in\Gamma(\xi)$, on pose

$$\mbox{(2)}\;\;\; \Delta_{\mu_+,\mu_-}(\xi_+,\xi_-,\gamma)=\mu_-(P_{\xi_-}(-1))\mu_+(P_{\xi_+}(-1))\prod_{i\in I_-} sgn_{F_i/F_{\pm i}}(C_i)$$

\noindent où cette fois

$$C_i=\left\{
    \begin{array}{ll}
      -\delta^{-d-1} \gamma_i^{-1} P_\xi '(y_i)P_\xi (-1)^{-1} y_i^{1-d/2} (1+y_i) & \mbox{ si } d \mbox{ est pair}\\
      -\delta^{-d-1} \gamma_i^{-1} P_\xi '(y_i)P_\xi (-1)^{-1} y_i^{(3-d)/2} & \mbox{ si } d \mbox{ est impair}
    \end{array}
\right.$$

\noindent Comme précédemment, les fonctions $\Delta_{\mu_+,\mu_-}$ sont en fait des facteurs de transfert, cette fois pour une endoscopie tordue (cf 3.2). 

\vspace{2mm}

Remarquons que dans les formules (1) et (2), on peut remplacer le produit sur $I_-$ par un produit sur $I_-^*$ car pour $i\in I_-\backslash I_-^*$, on a $sgn_{F_i/F_{\pm i}}=1$. On aura besoin du lemme suivant dont la démonstration, un calcul direct, est laissée au lecteur.

\begin{lem}
Soient $a\in E$ et $b\in Ker\; Tr_{E/F}$ tous deux non nuls. Pour $\lambda\in F^\times$ assez proche de $0$, on pose:

\begin{itemize}
\renewcommand{\labelitemi}{$\bullet$}

\item $\zeta_a(\lambda)=\left(E,E\times E, (e^{\lambda a},e^{-\lambda \overline{a}})\right)\in \Xi_{2,reg}$ (où on a identifié $E\otimes_F E$ et $E\times E$ via l'isomorphisme $x\otimes y\mapsto (xy,\overline{x}y)$);

\item $\zeta_b(\lambda)=\left(F,E,e^{\lambda b}\right)\in \Xi_{1,reg}$.
\end{itemize}

\noindent Soient $\xi_+,\xi_-\in \Xi_{reg}$ et $\mu_+$, $\mu_-$ deux caractères continus de $E^\times$. On pose $d_+=d_{\xi_+}$, $d_-=d_{\xi_-}$, $\underline{d}=d_++d_-$ et $\xi=\xi_+\sqcup \xi_-$. Remarquons que, pour tout $\lambda\in F^\times$ assez proche de $0$, on a

$$C(\xi\sqcup \zeta_a(\lambda))=C(\xi),\;\;\; C(\xi\sqcup \zeta_b(\lambda))=C(\xi)\times F^\times/N(E^\times)$$

$$\Gamma(\xi\sqcup \zeta_a(\lambda))=\Gamma(\xi), \;\;\; \Gamma(\xi\sqcup \zeta_b(\lambda))=\Gamma(\xi)\times \Gamma(\zeta_b(\lambda))$$

\noindent On a les égalités suivantes:

\begin{enumerate}[(A)]

\item Si $\mu_{+|N(E^\times)}=1$ et $\mu_{-|N(E^\times)}=1$, alors

\begin{enumerate}[(i)]

\item $$\Delta_{\mu_+,\mu_-,\nu}(\xi_+\sqcup\zeta_a(\lambda),\xi_-,c)=\mu_+(e^{\lambda a}) \Delta_{\mu_+,\mu_-,\nu}(\xi_+,\xi_-,c)$$

\item $$\Delta_{\mu_+,\mu_-,\nu}(\xi_+,\xi_-\sqcup \zeta_a(\lambda),c)=\mu_-(e^{\lambda a}) \Delta_{\mu_+,\mu_-,\nu}(\xi_+,\xi_-,c)$$

\item $$\Delta_{\mu_+,\mu_-}(\xi_+\sqcup\zeta_a(\lambda),\xi_-,\gamma)=\mu_+(e^{\lambda a}) \Delta_{\mu_+,\mu_-}(\xi_+,\xi_-,\gamma)$$

\item $$\Delta_{\mu_+,\mu_-}(\xi_+,\xi_-\sqcup \zeta_a(\lambda),\gamma)=\mu_-(e^{\lambda a}) \Delta_{\mu_+,\mu_-}(\xi_+,\xi_-,\gamma)$$

\end{enumerate}

\noindent pour tout $\nu\in F^\times$, pour tout $c\in C(\xi)$, pour tout $\gamma\in \Gamma(\xi)$ et pour tout $\lambda\in F^\times$ tel que $e^{\lambda a}$ soit défini.

\item Si $\mu_{+|F^\times}=sgn_{E/F}^{d_-}$, alors

\begin{enumerate}[(i)]

\item $$\displaystyle \lim\limits_{\lambda\to 0,\lambda\in F^\times} \Delta_{\mu_+,\mu_-,\nu}(\xi_+\sqcup \zeta_b(\lambda),\xi_-,(c,c_b))=\Delta_{\mu_+,\mu_-,\nu}(\xi_+ ,\xi_-,c) sgn_{E/F}\left(\delta^{-d_-} \frac{P_{\xi_-}(1)}{P_{\xi_-}(-1)}\right)$$

\item $$\displaystyle \lim\limits_{\lambda\to 0,\lambda\in F^\times} \Delta_{\mu_+,\mu_-}(\xi_+\sqcup \zeta_b(\lambda),\xi_-,(\gamma,c_be^{\lambda b/2}))=\Delta_{\mu_+,\mu_-}(\xi_+ ,\xi_-,\gamma) sgn_{E/F}\left(\delta^{-d_-} \frac{P_{\xi_-}(1)}{P_{\xi_-}(-1)}\right)$$

\end{enumerate}

\noindent pour tout $\nu\in F^\times$, pour tout $c\in C(\xi)$, pour tout $\gamma\in \Gamma(\xi)$ et pour tout $c_b\in F^\times/N(E^\times)$.

\item Si $\mu_{-|F^\times}=sgn_{E/F}^{d_+}$, alors

\[\begin{aligned}
\lim\limits_{\lambda\to 0,\lambda\in F^\times} \Delta_{\mu_+,\mu_-,\nu}(\xi_+,\xi_-\sqcup \zeta_b(\lambda), (c,c_b))=\Delta_{\mu_+,\mu_-,\nu} & (\xi_+,\xi_-,c)sgn_{E/F}\left(\delta^{-d_+} \frac{P_{\xi_+}(1)}{P_{\xi_+}(-1)}\right) \\
 & \cdot sgn_{E/F}((-1)^{\underline{d}+1}c_b\nu)
\end{aligned}\]

\noindent pour tout $\nu\in F^\times$, pour tout $c\in C(\xi)$ et pour tout $c_b\in F^\times/N(E^\times)$.

\item Si $\mu_{-|F^\times}=sgn_{E/F}^{d_++1}$, alors

\[\begin{aligned}
\lim\limits_{\lambda\to 0,\lambda\in F^\times} \Delta_{\mu_+,\mu_-}(\xi_+,\xi_-\sqcup \zeta_b(\lambda),(\gamma,c_be^{\lambda b/2}))= \Delta_{\mu_+,\mu_-} & (\xi_+,\xi_-,\gamma)sgn_{E/F}\left(\delta^{-d_+} \frac{P_{\xi_+}(1)}{P_{\xi_+}(-1)}\right) \\
 & \cdot sgn_{E/F}((-1)^{\underline{d}}c_b)
\end{aligned}\]

\end{enumerate}

\end{lem}

\section{Endoscopie}

Soit $G$ un groupe réductif connexe défini sur $F$ et soit $(H,s,{}^L\xi)$ une donnée endoscopique de $G$, cf [LS] 1.2. Cela définit une application, appelée la correspondance endoscopique, d'un sous-ensemble de $H_{reg}(F)/stconj$, dont on qualifie les éléments de $G$-fortement réguliers, vers $G_{reg}(F)/stconj$. On dira que $x\in G_{reg}(F)/conj$ et $y\in H_{reg}(F)/stconj$ se correspondent si $\phi^{st}_G/conj(x)$ est l'image de $y$ par cette application. Pour définir des facteurs de transfert, il convient de fixer une forme intérieure quasi-déployée $\underline{G}$ de $G$, un torseur intérieur $\psi_G: G\to\underline{G}$ et un épinglage défini sur $F$ de $\underline{G}$ à $\underline{G}(F)$-conjugaison près. Cette dernière donnée revient à celle d'une orbite nilpotente régulière de $\underline{\mathfrak{g}}(F)$ (c.f lemme 5.1.A [LS]). Avec ces données, on peut définir des facteurs de transfert relatifs $\Delta_{H,G}(y,x;\overline{y},\overline{x})$ où $y\in H_{reg}(F)/stconj$ et $x\in G_{reg}(F)/conj$ resp. $\overline{y}\in H_{reg}(F)/stconj$ et $\overline{x}\in G_{reg}(F)/conj$ se correspondent. D'après une remarque de Kottwitz, si on fixe un cocycle $u:Gal(\overline{F}/F)\to \underline{G}$ vérifiant $\psi_G\circ \sigma(g)=u(\sigma)\sigma\circ \psi_G(g)u(\sigma)^{-1}$ pour tous $g\in G$ et $\sigma\in Gal(\overline{F}/F)$ (il n'en existe pas toujours), on peut alors définir des facteurs de transfert absolus $\Delta_{H,G}(y,x)$. Ces facteurs de transfert dépendent du cocycle $u$ et non pas seulement de sa classe de cohomologie. On suppose que l'on a fixé un tel cocycle $u$. \\

Soient $\Theta$ une distribution localement intégrable invariante par conjugaison sur $G(F)$ et $\Theta^H$ une distribution localement intégrable invariante par conjugaison stable sur $H(F)$ (on dit alors que $\Theta^H$ est stable). On dira que $\Theta$ est un transfert de $\Theta^H$ si on a l'égalité

$$\displaystyle \Theta(x)D^G(x)^{1/2}=\sum_y \Theta^H(y)D^H(y)^{1/2} \Delta_{H,G}(y,x)$$

\noindent pour tout $x\in G_{reg}(F)/conj$, où la somme porte sur les $y\in H_{reg}(F)/stconj$ qui sont en correspondance avec $x$. \\

 Soit $(M,\widetilde{M})$ un groupe tordu. On suppose $M$ déployé et on fixe un élément $\tilde{\theta}\in \widetilde{M}(F)$. On suppose qu'il existe un épinglage défini sur $F$ de $M$ qui est invariant par $\theta_{\tilde{\theta}}$. On fixe alors une orbite nilpotente régulière de $\mathfrak{m}_{\tilde{\theta}}(F)$. Soit $(H,s,{}^L\xi)$ une donnée endoscopique de $(M,\widetilde{M})$, cf [KS] 2.1. On dispose alors d'une application, appelée correspondance endoscopique, d'un sous-ensemble de $H_{reg}(F)/stconj$, dont les éléments sont qualifiés de $\widetilde{M}$-fortement réguliers, vers $\widetilde{M}_{reg}(F)/stconj$. Deux éléments $y\in H_{reg}(F)/stconj$ et $\tilde{x}\in \widetilde{M}_{reg}(F)/conj$ sont dits en correspondance, si $\phi^{st}_{\widetilde{M}}/conj(\tilde{x})$ est l'image de $y$ par l'application précédente. On dispose d'un facteur de transfert $\Delta_{H,\tilde{M}}(y,\tilde{x})$ défini pour tout $y\in H_{reg}(F)/stconj$ et $\tilde{x}\in \widetilde{M}_{reg}(F)/conj$ et qui est non nul si et seulement si $\tilde{x}$ et $y$ se correspondent. Soient $\tilde{\Theta}$ une distribution localement intégrable invariante par conjugaison sur $\widetilde{M}(F)$ et $\Theta^H$ une distribution localement intégrable invariante par conjugaison stable sur $H(F)$. On dira que $\tilde{\Theta}$ est un transfert de $\Theta^H$, si on a l'égalité
 
$$\displaystyle \tilde{\Theta}(\tilde{x})D^{\tilde{M}}_0(\tilde{x})^{1/2}=\sum_y \Theta^H(y) D^H(y)^{1/2} \Delta_{H,\tilde{M}}(y,\tilde{x})$$

\noindent pour tout $\tilde{x}\in \widetilde{M}_{reg}(F)/conj$, où la somme porte sur les $y\in H_{reg}(F)/stconj$ qui sont en correspondance avec $\tilde{x}$.

\subsection{Endoscopie classique pour les groupes unitaires}

Soient $(V,h)$, $(V_+,h_+)$ et $(V_-,h_-)$ trois espaces hermitiens de dimensions respectives $d$, $d_+$, $d_-$ et de groupes unitaires respectifs $G$, $G_+$ et $G_-$. On fait les hypothèses suivantes

\vspace{3mm}

\begin{itemize}
\renewcommand{\labelitemi}{$\bullet$}

\item $d=d_++d_-$;

\item $G_+$ et $G_-$ sont quasi-déployés.
\end{itemize}

\vspace{3mm}

Alors $G_+\times G_-$ est un groupe endoscopique de $G$. Il y a un certain nombre de choix à faire pour fixer la donnée endoscopique (cf [W1] 1.8). Certains de ces choix sont inoffensifs (correspondant aux paramètres $z_+$ et $z_-$ de loc.cit.). Mais d'autres, dont dépendent les facteurs de transfert, doivent être précisés. Ce sont les caractères $\mu^+$ et $\mu^-$ du paragraphe 1.8 de [W1]. On suppose donc choisi deux caractères continus $\mu_+,\mu_-:E^\times\to \mathbb{C}^\times$ dont les restrictions à $F^\times$ coïncident respectivement avec $sgn_{E/F}^{d_-}$ et $sgn_{E/F}^{d_+}$. On fixe alors la donnée endoscopique comme en [W1] 1.8. \\

On doit encore fixer quelques données supplémentaires pour que les facteurs de transfert soient bien définis. A savoir une forme intérieure quasi-déployée $\underline{G}$ de $G$, un torseur intérieur $\psi_{G}: G\to \underline{G}$, un 1-cocycle $u: Gal(\overline{F}/F)\to\underline{G}$ et une orbite nilpotente régulière dans $\underline{\mathfrak{g}}(F)$. 

\begin{center}
\textbf{Fixons dorénavant (et jusqu'à la fin de l'article), un élément $\nu_0\in F^\times$}
\end{center}

Introduisons ici la condition \textbf{(QD)} suivante \\

 \textbf{(QD)}: les groupes unitaires de $(V,h)$ et $(V,h)\oplus(Z_1,h_{1,\nu_0})$ sont quasi-déployés. \\

\noindent Si $d$ est pair cette condition est équivalente à ce que $G$ soit quasi-déployé. On posera 

\[\begin{aligned}
\mu(V,h)=\left\{
    \begin{array}{ll}
        1 & \mbox{si } (V,h) \mbox{ vérifie \textbf{(QD)}} \\
        -1 & \mbox{sinon}
    \end{array}
\right.
\end{aligned}\]

\noindent On fixe les données supplémentaires $\left(\underline{G},\psi_G,u,\mathcal{O}\right)$ de la façon suivante

\vspace{3mm}

\begin{itemize}
\renewcommand{\labelitemi}{$\bullet$}

\item Si $\mu(V,h)=1$, on pose $\underline{G}=G$, $\psi_G=Id_G$ et $u=1$;

\item Si $\mu(V,h)=-1$, il existe un unique espace hermitien (à isomorphisme près) $(\underline{V},\underline{h})$ de dimension $d$ vérifiant \textbf{(QD)}. On pose $\underline{G}=U(\underline{V},\underline{h})$. Les formes sequilinéaires $h$ et $\underline{h}$ s'étendent naturellement en des formes $E\otimes_F \overline{F}$-sesquilinéaires sur $V\otimes_F \overline{F}$ et $\underline{V}\otimes_F\overline{F}$ respectivement. On note encore $h$ et $\underline{h}$ ces prolongements. Fixons un isomorphisme $E\otimes_F \overline{F}$-linéaire $\beta: V\otimes_F \overline{F}\to \underline{V}\otimes_F \underline{F}$ qui transforme $h$ en $\underline{h}$. On pose $\psi_{G}(g)=\beta g\beta^{-1}$ et $u(\sigma)=\beta\sigma(\beta)^{-1}$;

\item Si $d$ est impair ou $d=0$, il n'y a qu'une seule orbite nilpotente régulière dans $\mathfrak{g}(F)$. On la note $\mathcal{O}$;

\item Si $d\geqslant 2$ est pair, on pose $\mathcal{O}=\mathcal{O}_{\nu_0\delta}$ (où $\delta\in Ker(Tr_{E/F})-\{0\}$ est l'élément fixé en 1.1 et on se sert du paramétrage des orbites nilpotentes régulières de $\mathfrak{g}(F)$ donné en 1.3).
\end{itemize}

\vspace{3mm}

Ces données permettent de définir une correspondance entre classes de conjugaison stable semi-simples régulières de $G(F)$ et les classes de conjugaison stable semi-simples $G$-régulières  de $G_+(F)\times G_-(F)$. Via les applications $p^{st}_{G}$, $p^{st}_{G_+}\times p^{st}_{G_-}$, une classe de conjugaison semi-simple stable $(\xi_+,\xi_-)\in \Xi_{reg,d_+}\times \Xi_{reg,d_-}$ est $G$-régulière si et seulement si $\xi_+\sqcup \xi_-\in \Xi_{reg,d}$ et la classe de conjugaison stable $\xi\in Im(p^{st}_{G})$ lui correspond si et seulement si $\xi=\xi_+\sqcup \xi_-$. Soient $y_+\in G_{+,reg}(F)/stconj$, $y_-\in G_{-,reg}(F)/stconj$ et $y\in G_{reg}(F)/conj$. Posons $\xi_+=p^{st}_{G_+}(y_+)$, $\xi_-=p^{st}_{G_-}(y_-)$, $\xi=p_{G}(y)$ et supposons que $\xi=\xi_+\sqcup\xi_-$. Soit $c\in C(\xi)$ qui paramètre la classe de conjugaison de $y$. On a alors

$$\Delta_{G_+\times G_-,G}((y_+,y_-),y)=\Delta_{\mu_+,\mu_-,\nu_0}(\xi_+,\xi_-,c)$$

\noindent cf 1.10[W1], le $\eta$ de cette référence vaut ici $(-1)^{(d-2)/2}\delta\nu_0$ si $d$ est pair et $(-1)^{(d-3)/2}\nu_0$ si $d$ est impair. Remarquons que la formule 2.3(1) ne coïncide pas tout à fait avec celle de 1.10[W1]. En effet, on y a apporté la simplification

$$\mu_-\left(P_{\xi_-}(0)P_{\xi_-}(-1)^{-1}\right)\mu_+\left(P_{\xi_+}(0)P_{\xi_+}(-1)^{-1}\right)=\mu_-\left(P_{\xi_-}(-1)\right)\mu_+\left(P_{\xi_+}(-1)\right)$$

\subsection{Changement de base pour les groupes unitaires}

Soient $U$ un $E$-espace vectoriel de dimension $d$ et $(V_+,h_+)$, $(V_-,h_-)$ deux espaces hermitiens de dimension respectives $d_+$ et $d_-$. Introduisons le groupe tordu $(M,\widetilde{M})$ associé à $U$ comme en 2.2 et $G_+$, $G_-$ les groupes unitaires respectifs de $(V_+,h_+)$, $(V_-,h_-)$. On fait les hypothèses suivantes:

\vspace{3mm}

\begin{itemize}
\renewcommand{\labelitemi}{$\bullet$}

\item $G_+$ et $G_-$ sont quasi-déployés;

\item $d_++d_-=d$.
\end{itemize}

\vspace{3mm}

On peut alors considérer $G_+\times G_-$ comme un groupe endoscopique tordu de $(M,\widetilde{M})$. Pour cela il faut fixer une donnée supplémentaire (cf [W1] 1.8): deux caractères $\mu_+$ et $\mu_-$ de $E^\times$ dont les restrictions à $F^\times$ coïncident avec $sgn_{E/F}^{d_-}$ et $sgn_{E/F}^{d_++1}$ respectivement (il faudrait aussi fixer un couple de nombres complexes $(z_+,z_-)$ mais ce choix n'importe pas). La correspondance endoscopique, elle ne dépend pas du choix de ces caractères car la classe d'équivalence de la donnée endoscopique est bien définie. Elle dépend en revanche du choix d'un point base de $\widetilde{M}(F)$. Pour cela, on fixe une base $(u_j)_{j=1,\ldots,d}$ de $U$ et on prend pour point-base l'élément $\mathbf{\theta}_d\in \widetilde{M}(F)$ défini par

$$\mathbf{\theta}_d(u_j,u_k)=(-1)^{k}\delta^{d+1}\delta_{j,d+1-k}$$

\noindent On peut alors décrire la correspondance endocopique à l'aide des espaces de paramètres de la façon suivante: une classe de conjugaison stable de $G_+(F)\times G_-(F)$ paramétrée par $(\xi_+,\xi_-)\in \Xi_{reg,d_+}\times\Xi_{reg,d_-}$ est $\widetilde{M}$-régulière si $\xi_+\sqcup\xi_-\in \Xi_{reg,d}$ et de plus elle correspond à la classe de conjugaison stable $\xi\in \Xi_{d,reg}$ de $\widetilde{M}(F)$ si et seulement si $\xi=\xi_+\sqcup\xi_-$. Pour pouvoir définir des facteurs de transfert, il faut encore fixer une orbite nilpotente régulière dans $\mathfrak{m}_{\mathbf{\theta}_d}(F)$ (où, rappelons le, $M_{\mathbf{\theta}_d}$ est le centralisateur connexe de $\mathbf{\theta}_d$ dans $M$ et $\mathfrak{m}_{\mathbf{\theta}_d}$ est son algèbre de Lie). On choisit cette orbite nilpotente régulière de la façon suivante:

\begin{itemize}
\renewcommand{\labelitemi}{$\bullet$}

\item Si $d$ est impair il n'y en a qu'une seule et on choisit donc celle là;

\item si $d$ est pair on choisit l'orbite contenant l'élément $N$ défini par $Nu_1=0$ et $Nu_j=-u_{j-1}$ pour $j=2,\ldots,d$ (où $(u_j)_{j=1,\ldots,d}$ est la base de $U$ qui nous a servi à définir $\mathbf{\theta}_d$).
\end{itemize}

\vspace{4mm}

 Soient $\tilde{x}\in\widetilde{M}_{reg}(F)$ et $y=(y_+,y_-)\in G_{+,reg}(F)\times G_{-,reg}(F)$ dont les classes de conjugaison stable se correspondent. La classe de conjugaison de $\tilde{x}$ est paramétrée par $\xi \in \Xi_{reg,d}$ et $\gamma\in \Gamma(\xi)$ tandis que la classe de conjugaison stable de $y$ est paramétrée par $(\xi_+,\xi_-)\in \Xi_{d_+,reg}\times \Xi_{d_-,reg}$. On a alors

$$\Delta(\tilde{x},y)=\Delta_{\mu_+,\mu_-}(\xi_+,\xi_-,\gamma)$$

\noindent cf [W1] 1.10, le $\eta$ de cette référence vaut ici $(-1)^d\delta^{d+1}$. Remarquons que la formule 2.3(2) ne coïncide pas tout à fait avec celle de 1.10[W1]. En effet, on y a apporté la simplification

\[\begin{aligned}
\displaystyle \mu_-\left(P_{\xi_-}(0)P_{\xi_-}(-1)^{-1}\right)\mu_+\big(P_{\xi_+}(0) & P_{\xi_+}(-1)^{-1}\big) \\
 & =\mu_-\left(P_{\xi_-}(-1)\right)\mu_+\left(P_{\xi_+}(-1)\right)\prod_{i\in I_-} sgn_{F_i/F_{\pm i}}(-1)
\end{aligned}\]

\section{Deux formules revisitées}

\subsection{Rappel: une formule pour la multiplicité}

Soient $(V,h)$ et $(V',h')$ deux espaces hermitiens de dimensions respectives $d$ et $d'$ et de groupes unitaires respectifs $G$ et $G'$. On suppose que

\vspace{3mm}

\begin{itemize}
\renewcommand{\labelitemi}{$\bullet$}

\item $d$ pair et $d'$ impair
\end{itemize}

\vspace{3mm}

\noindent Posons $r=(d-d'-1)/2$ si $d>d'$ et $r=(d'-d-1)/2$ si $d'>d$. On effectue l'hypothèse supplémentaire suivante:

\vspace{4mm}

\begin{center} 
 (1) Si $d>d'$, $(V,h)$ est isomorphe à $(V',h')\oplus (Z_{2r+1},h_{2r+1,\nu_0})$ et si $d'>d$, $(V',h')$ est isomorphe à $(V,h)\oplus (Z_{2r+1},h_{2r+1,-\nu_0})$.
\end{center}

\vspace{4mm}

Fixons un tel isomorphisme, cela définit un plongement de $G$ dans $G'$ ou de $G'$ dans $G$ (suivant que $d'>d$ ou $d>d'$). Soient $\sigma\in Temp(G)$ et $\sigma'\in Temp(G')$. On définit comme en [B1] section 4 une multiplicité $m(\sigma,\sigma')$. Rappelons en la définition. On suppose pour fixer les idées que $d>d'$. Rappelons que $Z_{2r+1}$ admet une base $(z_i)_{i=0,\pm 1,\ldots,\pm r}$. Soit $P$ le sous-groupe parabolique de $G$ qui préserve le drapeau de sous-espaces totalement isotropes suivant

$$Ez_r\subset Ez_r\oplus Ez_{r-1}\subset\ldots\subset Ez_r\oplus\ldots\oplus Ez_1$$

\noindent On a alors $G'\subset P$. Notons $U$ le radical unipotent de $P$. On définit un caractère $\xi$ de $U(F)$ en posant

$$\displaystyle\xi(u)=\psi_E(\sum_{i=1}^r h(z_i,uz_{i-1}))$$

\noindent pour tout $u\in U(F)$. La conjugaison par $G'(F)$ stabilise $\xi$. Notons $Hom_{G',\xi}(\sigma,{\sigma'}^\vee)$ l'espace des applications linéaires $\ell: E_\sigma\to E_{{\sigma'}^\vee}$ (où $E_\sigma$ resp. $E_{{\sigma'}^\vee}$ désignent les espaces des représentations $\sigma$ et ${\sigma'}^\vee$ resp.) vérifiant

$$\ell\circ \sigma(g'u)=\xi(u){\sigma'}^\vee(g')\circ \ell$$

\noindent pour tous $g'\in G'(F)$, $u\in U(F)$. On pose alors $m(\sigma,\sigma')=dim\; Hom_{G',\xi}(\sigma,{\sigma'}^\vee)$. D'après [AGRS] et [GGP], cette multiplicité vaut $0$ ou $1$. Remarquons que notre convention diffère ici de celle de [B1]: la multiplicité $m(\sigma,\sigma')$ est celle notée $m(\sigma,{\sigma'}^\vee)$ dans [B1]. On a démontré en [B1](théorème 17.1.1) une formule calculant cette multiplicité. Explicitons la. Elle comprend deux ingrédients essentiels: un ensemble de tores $\mathcal{T}$ ainsi que deux fonctions $c_\sigma$ et $c_{\sigma'}$ définis à partir des caractères des représentations $\sigma$ et $\sigma'$. Commençons par rappeler la définition de $c_\sigma$ et $c_{\sigma'}$. Plaçons nous dans le cadre plus général suivant: $H$ est un groupe réductif connexe défini sur $F$ et $\pi$ est une représentation irréductible lisse de $H(F)$. A $\pi$, on peut associer une fonction

$$c_\pi: H_{ss}(F)\to\mathbb{C}$$

\noindent de la façon suivante. Soit $x\in H_{ss}(F)$. D'après Harish-Chandra ([HCDS] theorem 16.2), le caractère $\Theta_{\pi}$ de la représentation $\pi$ admet au voisinage de $x$ un développement de la forme

\begin{center}
$\displaystyle \Theta_\pi(xexp(X))=\sum_{\mathcal{O}\in Nil(\mathfrak{g}_x(F))} c_{\pi,\mathcal{O}}(x)\hat{j}(\mathcal{O},X)$, pour $X$ proche de $0$ dans $\mathfrak{g}_x(F)$
\end{center}

\noindent où les $c_{\pi,\mathcal{O}}(x)$ sont des nombres complexes. Tous les autres termes ci-dessus ont été définis dans la section 1.2. Maintenant, on pose

$$\mbox{(2)}\;\;\; \displaystyle c_\pi(x)=\frac{1}{|Nil(\mathfrak{g}_x)_{reg}|}\sum_{\mathcal{O}\in Nil(\mathfrak{g}_x)_{reg}} c_{\pi,\mathcal{O}}(x)$$

\noindent où, rappelons le, $Nil(\mathfrak{g}_x)_{reg}$ désigne l'ensemble des orbites nilpotentes régulières dans $\mathfrak{g}_x(F)$. Cette construction générale fournit en particulier deux fonctions

$$c_\sigma: G_{ss}(F)\to\mathbb{C}$$

$$c_{\sigma'}: G'_{ss}(F)\to \mathbb{C}$$

\noindent Reste à définir l'ensemble de tores $\mathcal{T}$. Définissons d'abord un ensemble d'espaces hermitiens $\mathcal{H}$. Supposons $d>d'$. Alors $\mathcal{H}$ est l'ensemble des sous-espaces hermitiens $W$ de $V'$ tels que les groupes unitaires des orthogonaux de $W$ dans $V$ et $V'$ soient quasi-déployés. Pour $W\in \mathcal{H}$, on note $U(W)$ le groupe unitaire de $W$ et $\mathcal{T}_{ell}(W)$ l'ensemble des tores maximaux elliptiques de $U(W)$. Pour $T\in \mathcal{T}_{ell}(W)$, on note $W(T)=Norm_{U(W)(F)}(T)/T(F)$ le groupe de Weyl de $T$. Posons

$$\displaystyle \underline{\mathcal{T}}=\bigsqcup_{W\in \mathcal{H}} \mathcal{T}_{ell}(W)$$

\noindent Alors, $\mathcal{T}$ sera un ensemble de représentants dans $\underline{\mathcal{T}}$ des classes de conjugaison par $G'(F)$. Les définitions sont les même dans le cas $d'>d$ en intervertissant $V$, $V'$ et $G$, $G'$.

\vspace{2mm}

\noindent Posons

\[\begin{aligned}
\displaystyle m_{geom}(\sigma,\sigma') & =m_{geom}(\sigma',\sigma)= \\
 & \sum_{T\in \mathcal{T}} |W(T)|^{-1} \lim\limits_{s\to 0^+} \int_{T(F)} c_\sigma(t)\cdot c_{\sigma'}(t) \cdot D^G(t)^{1/2} \cdot D^{G'}(t)^{1/2} \cdot \Delta(t)^{-1/2+s} dt
\end{aligned}\]

\noindent où la fonction $\Delta$ a été définie en 2.3. Cette expression a un sens et d'après le théorème 17.1.1 de [B1], on a l'égalité

$$\mbox{(3)}\;\;\; m(\sigma,\sigma')=m_{geom}(\sigma,\sigma')$$

\subsection{Reformulation de $m_{geom}(\sigma,\sigma')$}

On conserve les notations et hypothèses du paragraphe précédent. Notons $\mathcal{D}(h,h')$ l'ensemble des entiers $\underline{d}\in\mathbb{N}$ vérifiant les conditions suivantes :

\vspace{3mm}

\begin{itemize}
\renewcommand{\labelitemi}{$\bullet$}
\item $0\leqslant \underline{d}\leqslant min(d,d')$;
\item si $G$ ou $G'$ n'est pas quasi-déployé, $\underline{d}\geqslant 2$.
\end{itemize}

\vspace{3mm}

\noindent Pour tout $\underline{d}\in\mathcal{D}(h,h')$, il existe un unique espace hermitien (à isomorphisme près) $(V_{\underline{d}},h_{\underline{d}})$ vérifiant les conditions suivantes :

\vspace{3mm}

\begin{itemize}
\renewcommand{\labelitemi}{$\bullet$}
\item $dim(V_{\underline{d}})=\underline{d}$;
\item Il existe deux espaces hermitiens $(V_\natural,h_\natural)$, $(V'_\natural,h'_\natural)$ dont les groupes unitaires sont quasi-déployés, tels que $(V,h)$ soit isomorphe à $(V_{\underline{d}},h_{\underline{d}})\oplus(V_\natural,h_\natural)$ et $(V',h')$ soit isomorphe à $(V_{\underline{d}},h_{\underline{d}})\oplus (V'_\natural,h'_\natural)$.
\end{itemize}

\vspace{3mm}

Notons $G_{\underline{d}}$ le groupe unitaire de $(V_{\underline{d}},h_{\underline{d}})$. On peut choisir des isomorphismes $(V,h)\simeq (V_{\underline{d}},h_{\underline{d}})\oplus(V_\natural,h_\natural)$ et $(V',h')\simeq (V_{\underline{d}},h_{\underline{d}})\oplus (V'_\natural,h'_\natural)$ avec $(V_\natural,h_\natural)$, $(V'_\natural,h'_\natural)$ comme précédemment. Cela définit des plongements de $G_{\underline{d}}$ dans $G$ et $G'$ bien définis à conjugaison près. Quitte à changer les isomorphismes, on peut toujours supposer que ces plongements sont compatibles avec le plongement de $G$ dans $G'$ ou de $G'$ dans $G$ (suivant que $d'>d$ ou $d>d'$). Pour $x\in G_{\underline{d}}(F)_{ani}$, l'image commune du tore $G_{\underline{d},x}$ par ces deux plongements appartient à $\underline{\mathcal{T}}$. Puisque $x\in G_{\underline{d},x}(F)$, on peut définir les coefficients $c_\sigma(x)$ et $c_{\sigma'}(x)$. Introduisons les espaces 

$$\displaystyle \mathcal{C}(h,h')=\bigsqcup_{\underline{d}\in \mathcal{D}(h,h')} G_{\underline{d}}(F)_{ani}/conj$$

\noindent et

$$\displaystyle \Xi^*(h,h')=\bigsqcup_{\underline{d}\in \mathcal{D}(h,h')} \Xi^*_{reg,\underline{d}}$$

\noindent D'après 2.2, $\mathcal{C}(h,h')$ est un revêtement de $\Xi^*(h,h')$. L'égalité 4.1(3) se réécrit

$$\mbox{(1)}\;\;\; \displaystyle m(\sigma,\sigma')=\lim\limits_{s\to 0^+}\int_{\mathcal{C}(h,h')} c_\sigma(x)\cdot c_{\sigma'}(x)\cdot D^d(x)^{1/2}\cdot D^{d'}(x)^{1/2} \cdot \Delta(x)^{-1/2+s} dx$$

\subsection{Rappel: une formule pour un facteur epsilon de paire}

Soient $U$ et $U'$ deux espaces vectoriels sur $E$ de dimensions finies $d$ et $d'$ respectivement. Posons $M=R_{E/F} GL(U)$ et soit $\widetilde{M}$ la variété algébrique des formes sesquilinéaires non dégénérées sur $U$. Comme on l'a vu en 2.2, le couple $(M,\widetilde{M})$ est naturellement un groupe tordu. On introduit de même le groupe tordu $(M',\widetilde{M}')$ relatif à $U'$. On suppose que

\vspace{3mm}

\begin{itemize}
\renewcommand{\labelitemi}{$\bullet$}

\item $d$ est pair et $d'$ est impair.
\end{itemize}

\vspace{3mm}

\noindent Posons $r=(d-d'-1)/2$ si $d>d'$ et $r=(d'-d-1)/2$ si $d'>d$.

\begin{center}
\textbf{Fixons (jusqu'à la fin de l'article) un élément $\nu_1\in F^\times$}
\end{center}

\noindent On effectue la construction suivante:

\vspace{3mm}

\begin{itemize}
\renewcommand{\labelitemi}{$\bullet$}

\item Si $d>d'$, on fixe un isomorphisme $U\simeq U'\oplus Z_{2r+1}$ et on plonge $\widetilde{M}'$ dans $\widetilde{M}$ en envoyant $\tilde{m}'\in \widetilde{M}'$ sur la somme directe $\tilde{m}'\oplus h_{2r+1,\nu_1}$;

\item Si $d'>d$, on fixe un isomorphisme $U'\simeq U\oplus Z_{2r+1}$ et on plonge $\widetilde{M}$ dans $\widetilde{M}'$ en envoyant $\tilde{m}\in \widetilde{M}$ sur la somme directe $\tilde{m}\oplus h_{2r+1,-\nu_1}$;

\end{itemize}

\vspace{3mm}

 Soient $(\pi,E_\pi)$ et $(\pi',E_{\pi'})$ des représentations irréductibles tempérées de $M(F)$ et $M'(F)$. On suppose que $\pi$ et $\pi'$ se prolongent en des représentations unitaires de $\widetilde{M}(F)$ et $\widetilde{M}'(F)$. C'est-à-dire que l'on peut trouver une application $\tilde{\pi}$ de $\widetilde{M}(F)$ dans le groupe des automorphismes unitaires de $E_\pi$ vérifiant $\tilde{\pi}(m\tilde{m}m')=\pi(m)\tilde{\pi}(\tilde{m})\pi(m')$ pour tout $\tilde{m}\in\widetilde{M}(F)$ et pour tous $m,m'\in M(F)$ (et de même pour $\pi'$). On dira alors que $\pi$ et $\pi'$ sont \ul{conjuguées-duales}. Ces prolongements ne sont uniques qu'à un facteur complexe de module $1$ près. Néanmoins, il est possible de rendre le choix unique grâce à la théorie des fonctionnelles de Whittaker. Pour cela fixons une base $u_1,\ldots,u_d$ de $U$ et identifions un élément de $M$ à sa matrice dans cette base. Notons $N$ le sous-groupe des matrices unipotentes supérieures. On définit un caractère $\psi_N$ de $N(F)$ par la formule

$$\displaystyle\psi_N(n)=\psi_E(\sum_{i=1}^{d-1} n_{i,i+1})$$

\noindent où pour $n\in N(F)$, on a noté $n_{ij}$ les coefficients de la matrice associée. La théorie des fonctionnelles de Whittaker nous dit alors que l'espace $Hom_{N(F)}(\pi,\psi_N)$ est de dimension $1$. Si $\tilde{\pi}$ est un prolongement de $\pi$ à $\widetilde{M}(F)$, l'application $\ell\mapsto \ell\circ \tilde{\pi}(\theta_d)$ est un automorphisme de $Hom_{N(F)}(\pi,\psi_N)$ (on renvoie à 3.2 pour la définition de $\theta_d$). On normalise alors $\tilde{\pi}$ de sorte que cet automorphisme soit l'identité. Un simple calcul de changement de base montre que ce prolongement ne dépend pas du choix de la base $u_1,\ldots,u_d$ (en revanche il dépend de $\psi_E$). On normalise de la même façon $\tilde{\pi}'$. \\

 Définissons la quantité

$$\epsilon_{\nu_1}(\pi,\pi')=\epsilon_{\nu_1}(\pi',\pi)=\omega_\pi(\nu_1) \omega_{\pi'}(-\nu_1) \epsilon(1/2,\pi\times \pi',\psi^\delta_E)$$

\noindent où $\epsilon(s,\pi\times \pi',\psi^\delta_E)$ est le facteur epsilon défini par Jacquet, Piatetski-Shapiro et Shalika ([JPPS]) et $\psi^\delta_E$ est le caractère additif de $E$ définit par $\psi_E^\delta(x)=\psi_E(\delta x)$. Le résultat principal de [B2] calcule ce terme par une formule intégrale. Rappelons cette formule. Comme pour la formule pour la multiplicité, celle-ci comprend deux éléments essentiels: un ensemble de tores tordus $\mathcal{T}$ ainsi que deux fonctions $c_{\tilde{\pi}}$ et $c_{\tilde{\pi}'}$ définis à partir des caractères de $\tilde{\pi}$ et $\tilde{\pi}'$. On rappelle d'abord la définition de $c_{\tilde{\pi}}$ et $c_{\tilde{\pi}'}$. D'après Clozel ([C]), les représentations tordus $\tilde{\pi}$ et $\tilde{\pi}'$ admettent, comme dans le cas non tordu, des caractères $\Theta_{\tilde{\pi}}$ et $\Theta_{\tilde{\pi}'}$. Ce sont des fonctions sur les groupes tordus $\widetilde{M}(F)$ et $\widetilde{M}'(F)$ localement intégrables et localement constantes sur les lieux réguliers. Ces caractères admettent des développements locaux analogues à ceux des caractères dans le cas non tordu. Plus précisément, dans le cas de $\Theta_{\tilde{\pi}}$ par exemple, pour tout point semi-simple $\tilde{x}\in \widetilde{M}_{ss}(F)$, on a un développement au voisinage de $\tilde{x}$ de la forme

\begin{center}
$\displaystyle \Theta_{\tilde{\pi}}(\tilde{x}exp(X))=\sum_{\mathcal{O}\in Nil(\mathfrak{m}_{\tilde{x}})} c_{\tilde{\pi},\mathcal{O}}(\tilde{x}) \hat{j}(\mathcal{O},X)$, pour tout $X$ proche de $0$ dans $\mathfrak{m}_{\tilde{x}}(F)$
\end{center}

\noindent où les $c_{\tilde{\pi},\mathcal{O}}(\tilde{x})$ sont des nombres complexes et les autres termes de ce développement ont été définis en 1.2. On pose alors

$$\displaystyle c_{\tilde{\pi}}(\tilde{x})=\frac{1}{|Nil(\mathfrak{m}_{\tilde{x}})_{reg}|}\sum_{\mathcal{O}\in Nil(\mathfrak{m}_{\tilde{x}})_{reg}} c_{\tilde{\pi},\mathcal{O}}(\tilde{x})$$

\noindent Cela définit une fonction $c_{\tilde{\pi}}: \widetilde{M}_{ss}(F)\to\mathbb{C}$. On définit de la même façon une fonction $c_{\tilde{\pi}'}: \widetilde{M}'_{ss}(F)\to\mathbb{C}$. \\

 Reste maintenant à définir l'ensemble de tores tordus $\mathcal{T}$. Supposons $d>d'$. Notons $\mathcal{E}$ l'ensemble des sous-espaces $V$ de $U'$. Pour $V\in \mathcal{E}$, on note $(M(V),\widetilde{M}(V))$ le groupe tordu associé à $V$. On plonge $\widetilde{M}(V)$ dans $\widetilde{M}'$ de la façon suivante. On choisit un supplémentaire $W$ de $V$ dans $U'$ (c'est-à-dire $U'=V\oplus W$) et une forme hermitienne non dégénérée $\zeta_W$ sur $W$ de sorte que les groupes unitaires $U(\zeta_W)$ et $U(\zeta_W\oplus h_{2r+1,\nu_1})$ soient quasi-déployés. On envoie alors $\tilde{x}\in \widetilde{M}(V)$ sur la forme sesquilinéaire $\tilde{x}\oplus \zeta_W\in \widetilde{M}'$. Notons $\mathcal{T}_{ell}(V)$ l'ensemble des tores tordus maximaux $\widetilde{T}$ de $\widetilde{M}(V)$ qui sont elliptiques (i.e. tels que $A_{\widetilde{T}}=\{1\}$). Pour un tel tore tordu $\widetilde{T}$, on notera $T$ le tore associé (i.e. le tore sous lequel $\widetilde{T}$ est un espace principal homogène) et $W(\widetilde{T})=Norm_{M(V)(F)}(\widetilde{T})/T(F)$ son groupe de Weyl. Posons
 
$$\displaystyle \underline{\mathcal{T}}=\bigsqcup_{V\in \mathcal{E}} \mathcal{T}_{ell}(V)$$

\noindent Alors $\mathcal{T}$ désigne un ensemble de représentants dans $\underline{\mathcal{T}}$ des classes de conjugaison par $M'(F)$. Les définitions sont les même dans le cas $d'>d$ en intervertissant $U$ et $U'$, $\widetilde{M}$ et $\widetilde{M}'$ ainsi que $h_{2r+1,\nu_1}$ et $h_{2r+1,-\nu_1}$. \\

Pour tout $\widetilde{T}\in \mathcal{T}$, la conjugaison par $\widetilde{T}$ conserve $T$ et l'action induite sur $T$ est la même pour tout les éléments de $\widetilde{T}$. On la notera $\theta$. Soient $T^\theta$ le sous-groupe des points fixes par $\theta$ et $T_\theta$ la composante neutre de $T^\theta$. On note $\widetilde{T}(F)/\theta$ l'espace des classes de conjugaison dans $\widetilde{T}(F)$ par $T(F)$. On peut munir cet espace d'une unique mesure de sorte que pour tout $\tilde{t}\in \widetilde{T}(F)$, l'application $T_\theta(F)\to \widetilde{T}(F)/\theta$, $t\mapsto t\tilde{t}$ préserve localement les mesures

\vspace{3mm}

\noindent Posons

\[\begin{aligned}
\displaystyle  & \epsilon_{geom,\nu_1}(\tilde{\pi}, \tilde{\pi}')= \\
 & \sum_{\widetilde{T}\in \mathcal{T}} |2|_F^{dim(\widetilde{T})-(d+d')/2}|W(\widetilde{T})|^{-1} \lim\limits_{s\to 0^+} \int_{\widetilde{T}(F)/\theta} c_{\tilde{\pi}}(\tilde{t}) \cdot c_{\tilde{\pi}'}(\tilde{t})\cdot D^{\tilde{M}}(\tilde{t})^{1/2} \cdot D^{\tilde{M}'}(\tilde{t})^{1/2} \cdot \Delta(\tilde{t})^{s-1/2} d\tilde{t}
\end{aligned}\]

\noindent Malgré les apparences, le membre de droite ci-dessus dépend bien de $\nu_1$. En effet, $\nu_1$ nous a servi à fixer les plongements $\widetilde{M}(V)\hookrightarrow \widetilde{M}'\hookrightarrow \widetilde{M}$ (resp. $\widetilde{M}(V)\hookrightarrow \widetilde{M}\hookrightarrow \widetilde{M}'$) pour $d'<d$ (resp. $d<d'$). L'expression ci-dessus a un sens et d'après le théorème 6.1.1 de [B2], on a l'égalité

$$\mbox{(1)}\;\;\; \epsilon_{\nu_1}(\pi,\pi')=\epsilon_{geom,\nu_1}(\tilde{\pi},\tilde{\pi}')$$

\subsection{Reformulation de $\epsilon_{geom,\nu_1}(\tilde{\pi},\tilde{\pi}')$}

On conserve la situation du paragraphe précédent. Pour $\underline{d}$ un entier naturel vérifiant $\underline{d}\leqslant min(d,d')$, fixons:

\vspace{3mm}

\begin{itemize}
\renewcommand{\labelitemi}{$\bullet$}

\item $U_{\underline{d}}$ un $E$-espace vectoriel de dimension $\underline{d}$;

\item des identifications $U=U_{\underline{d}}\oplus W_{\underline{d}}$ et $U'=U_{\underline{d}}\oplus W'_{\underline{d}}$ de sorte que

\begin{itemize}
\renewcommand{\labelitemi}{$\bullet$}
\item Si $d<d'$, $W'_{\underline{d}}=W_{\underline{d}}\oplus Z_{2r+1}$;
\item si $d>d'$, $W_{\underline{d}}=W'_{\underline{d}}\oplus Z_{2r+1}$.
\end{itemize}

\item des formes hermitiennes non dégénérées $h_{\underline{d}}$ sur $W_{\underline{d}}$ et $h'_{\underline{d}}$ sur $W'_{\underline{d}}$ vérifiant les conditions suivantes:

\begin{itemize}
\renewcommand{\labelitemi}{$\bullet$}

\item Les groupes unitaires de $h_{\underline{d}}$ et $h'_{\underline{d}}$ sont quasi-déployés;

\item Si $d<d'$, $h'_{\underline{d}}$ est égale à $h_{\underline{d}}\oplus h_{2r+1,-\nu_1}$;

\item Si $d>d'$, $h_{\underline{d}}$ est égale à $h'_{\underline{d}}\oplus h_{2r+1,\nu_1}$.
\end{itemize}
\end{itemize}

\vspace{3mm}

\noindent (Rappelons que $h_{2r+1,\pm \nu_1}$ sont des formes hermitiennes non dégénérées sur $Z_{2r+1}$). Notons $(M_{\underline{d}},\widetilde{M}_{\underline{d}})$ le groupe tordu associé à $U_{\underline{d}}$. On plonge $\widetilde{M}_{\underline{d}}$ dans $\widetilde{M}$ en identifiant $\tilde{x}\in \widetilde{M}_{\underline{d}}$ à la somme directe de $\tilde{x}$ et de $h_{\underline{d}}$. On plonge de même $\widetilde{M}_{\underline{d}}$ dans $\widetilde{M}'$. Ces plongements sont alors compatibles avec le plongement de $\widetilde{M}$ dans $\widetilde{M}'$ ou de $\widetilde{M}'$ dans $\widetilde{M}$ (suivant que $d<d'$ ou $d>d'$). Alors si $\widetilde{T}_{\underline{d}}$ est un tore maximal elliptique de $\widetilde{M}_{\underline{d}}$ défini sur $F$, son image commune par chacun de ces plongements appartient à $\underline{\mathcal{T}}$. Ceci permet de définir $c_{\tilde{\pi}}(\tilde{x})$ et $c_{\tilde{\pi}'}(\tilde{x})$ pour tout $\tilde{x}\in \widetilde{M}_{\underline{d}}(F)_{ani}$. On définit de même les termes $D^{\tilde{M}}(\tilde{x})$ et $D^{\tilde{M}'}(\tilde{x})$. Introduisons l'espace

$$\displaystyle \mathcal{X}(d,d')=\bigsqcup_{\underline{d}\leqslant min(d,d')} \widetilde{M}_{\underline{d}}(F)_{ani}/conj$$

\noindent Cet espace est muni d'une mesure d'après 1.2. Alors, l'égalité 4.3(1) se réécrit

$$\mbox{(1)}\;\;\; \displaystyle \epsilon_{\nu_1}(\tilde{\pi},\tilde{\pi}')=\lim\limits_{s\to 0^+}\int_{\mathcal{X}(d,d')} |2|_F^{\underline{d}-(d+d')/2} c_{\tilde{\pi}}(\tilde{x})c_{\tilde{\pi}'}(\tilde{x})D^{\tilde{M}}(\tilde{x})^{1/2} D^{\tilde{M}'}(\tilde{x})^{1/2} \Delta(\tilde{x})^{-1/2+s} d\tilde{x}$$

\subsection{Rappel sur les fonctions $c_\sigma$}

Soit $(V,h)$ un espace hermitien. Donnons nous une décomposition orthogonale $V=W\oplus V_\natural$. On note $G$, $G_W$, $G_\natural$ les groupes unitaires de respectivement $V$, $W$, $V_\natural$ et on suppose que $G_\natural$ est quasi-déployé. Soit $\sigma$ une représentation irréductible de $G(F)$ et $x\in G_{W,reg}(F)$. On a défini en 4.1(2) un coefficient $c_\sigma(x)$ à partir du caractère $\Theta_\sigma$ de $\sigma$. On va rappeler une formule permettant de calculer ce terme. Fixons un sous-groupe de Borel $B_\natural$ de $G_\natural$ et un tore maximal $T_\natural$ de $B_\natural$ tous deux définis sur $F$. Soit $Y_\natural\in \mathfrak{t}_\natural(F)\cap \mathfrak{g}_{\natural,reg}(F)$ qui n'a pas $0$ comme valeur propre, on a alors

$$\mbox{(1)}\;\;\; \displaystyle c_\sigma(x)=w(d_\natural)^{-1}\lim\limits_{\lambda\in F^\times, \\ \lambda\to 0} \Theta_\sigma(xexp(\lambda Y_\natural)) D^{G_\natural}(\lambda Y_\natural)^{1/2}$$

\noindent cf 18.4 [B1].

\section{Modèles de Whittaker pour un groupe unitaire de dimension paire}

Soit $(V,h)$ un espace hermitien de dimension $d$. Notons $G$ son groupe unitaire. Dans toute cette section on fait l'hypothèse suivante

\begin{center}
$G$ est quasi-déployé et $d$ est pair
\end{center}

\noindent En 1.3, on a paramétré les orbites nilpotentes régulières de $\mathfrak{g}(F)$ par $(Ker Tr_{E/F}\backslash \{0\})/N(E^\times)$. Le choix d'un caractère continu de $F$ (on en a fixé un $\psi$) et d'une dualité $\mathfrak{g}(F)\times \mathfrak{g}(F)\to F$, $G$-invariante (on en a fixé une en 1.3) déterminent alors un paramètrage analogue des types de modèles de Whittaker pour $G$. Ce paramètrage est décrit explicitement en 5.2. On rappelle alors un résultat de Rodier exprimant l'existence ou non d'un type de modèle de Whittaker pour une représentation irréductible $\sigma$ en fonction de son caractère (ou plutôt du comportement au voisinage de l'origine de ce dernier). Le lemme 5.5.1 donne une façon explicite de calculer le coefficient qui apparaît dans la formule de Rodier. Pour arriver à cette proposition, qui ne dépend d'aucune mesure, on a besoin de changer un peu nos normalisations (notamment les mesures sur les tores). Ces changements, ainsi qu'un certain nombre de notations et résultats utiles, sont introduits dans le paragraphe suivant. Enfin, dans la section 5.6, on applique le lemme 5.5.1 pour obtenir une certaine formule de transfert qui nous sera utile par la suite.

\subsection{Intégrales orbitales, transformée de Fourier, mesures et constantes de Weil}

Rappelons que l'on a fixé en 1.3 une forme bilinéaire symétrique non dégénérée $<.,.>$ sur $\mathfrak{g}(F)$. On en a alors déduit une mesure sur $\mathfrak{g}(F)$: celle qui est autoduale pour le bicaractère $\psi(<.,.>):\mathfrak{g}(F)\times \mathfrak{g}(F)\to \mathbb{C}^\times$. Jusqu'à présent, la mesure fixée sur $G(F)$ était inessentielle. Dans cette section, on choisit sur $G(F)$ l'unique mesure de Haar telle que l'exponentielle préserve localement les mesures à l'origine. \\

Soit $X\in\mathfrak{g}_{reg}(F)$. Alors $G_X$ est un tore et la restriction de $<.,.>$ à $\mathfrak{g}_X(F)$ est non dégénérée. On munit alors $\mathfrak{g}_X(F)$ de la mesure autoduale associée. On a déjà fixé en 1.2 une mesure sur les tores. On oublie jusqu'en 5.5 cette normalisation et on relève sur $G_X(F)$ la mesure de $\mathfrak{g}_X(F)$ via l'exponentielle (comme pour $G(F)$). Avec cette normalisation des mesures, on pose

$$\displaystyle J_G(X,f)=\int_{G_X(F)\backslash G(F)} f(g^{-1}Xg) dg$$

\noindent pour tout $f\in C_c^\infty(\mathfrak{g}(F))$. C'est l'intégrale orbitale de $f$ en $X$. D'après Harish-Chandra, il existe une fonction $\hat{j}$ localement intégrable sur $\mathfrak{g}(F)\times \mathfrak{g}(F)$, localement constante sur $\mathfrak{g}_{reg}(F)\times \mathfrak{g}_{reg}(F)$, telle que on ait l'égalité   

$$\displaystyle J_G(X,\hat{f})=\int_{\mathfrak{g}(F)} f(Y) \hat{j}(X,Y) dY$$

\noindent pour tout $f\in C_c^\infty(\mathfrak{g}(F))$ et pour tout $X\in \mathfrak{g}_{reg}(F)$. Lorsque l'on voudra préciser le groupe ambiant, on notera $\hat{j}^G(.,.)$ plutôt que $\hat{j}(.,.)$. Pour $\mathcal{O}\in Nil(\mathfrak{g}(F))$, il existe une unique fonction $\Gamma_\mathcal{O}$ sur $\mathfrak{g}_{reg}(F)$, le germe de Shalika associé à $\mathcal{O}$, qui vérifie les conditions suivantes. La première est une condition d'homogénéité: pour tout $X\in\mathfrak{g}_{reg}(F)$ et pour tout $\lambda\in F^{\times}$, on a

$$\Gamma_{\mathcal{O}}(\lambda X)=|\lambda|_F^{(\delta(G)-dim(\mathcal{O}))/2}\Gamma_{\lambda \mathcal{O}}(X)$$

\noindent (rappelons que pour $\lambda\in F^{\times,2}$, on a $\lambda \mathcal{O}=\mathcal{O}$). La deuxième condition est la suivante: pour tout $f\in C_c^\infty(\mathfrak{g}(F))$ il existe un voisinage $\omega$ de $0$ tel que pour tout $X\in \omega\cap \mathfrak{g}_{reg}(F)$, on ait

$$J_G(X,f)=\displaystyle\sum_{\mathcal{O}\in Nil(\mathfrak{g}(F))} \Gamma_\mathcal{O}(X) J_\mathcal{O}(f)$$

\noindent D'après la "conjecture de Howe", si $\omega$ est un ouvert de $\mathfrak{g}(F)$ compact modulo conjugaison et contenant $0$, alors il existe un ouvert $\Omega$ de $\mathfrak{g}(F)$ compact modulo conjugaison et contenant $0$ tel que, pour tout $X\in \Omega\cap\mathfrak{g}_{reg}(F)$ et tout $Y\in \omega\cap\mathfrak{g}_{reg}(F)$, on ait l'égalité

$$\mbox{(1)}\;\;\; \displaystyle \hat{j}(X,Y)=\sum_{\mathcal{O}\in Nil(\mathfrak{g}(F))} \Gamma_\mathcal{O}(X) \hat{j}(\mathcal{O},Y)$$

\noindent On aura encore besoin d'une formule sur les fonctions $\hat{j}$. Soit $M$ un Levi de $G$. Alors la restriction de $<.,.>$ à $\mathfrak{m}(F)$ est non dégénérée. On peut donc appliquer toutes les constructions précédentes à $M$. En particulier, on dispose de fonctions $\hat{j}^M:\mathfrak{m}(F)\times \mathfrak{m}(F)\to\mathbb{C}$. Soient $X\in\mathfrak{m}(F)\cap\mathfrak{g}_{reg}(F)$ et $Y\in\mathfrak{g}_{reg}(F)$. Fixons un ensemble de représentants $(Y_i)_{i=1,\ldots,r}$ des classes de conjugaison par $M(F)$ dans $\mathfrak{m}(F)$ qui rencontrent la classe de conjugaison de $Y$ par $G(F)$. On a alors l'égalité

$$\mbox{(2)}\;\;\; \displaystyle \hat{j}^G(X,Y)D^G(Y)^{1/2}=\sum_{i=1}^r \hat{j}^M(X,Y_i)D^M(Y_i)^{1/2}$$

\noindent  On vérifie aussi que si $M$ est un tore, alors $\hat{j}^M(X,Y)=\psi(<X,Y>)$ pour tous $X,Y\in\mathfrak{m}(F)$. \\

 Soit $(\mathcal{V},q)$ un $F$-espace vectoriel de dimension finie muni d'une forme bilinéaire symétrique non dégénérée $q$. Pour $L$ un réseau de $\mathcal{V}$, on pose
 
$$\displaystyle I_\psi(L,q)=\int_L \psi(q(v,v)/2) dv$$

\noindent Alors pour $L$ assez grand, $I_\psi(L,q)/|I_\psi(L,q)|$ ne dépend pas de $L$, on note sa valeur $\gamma_\psi(q)$. C'est la constante de Weil de $q$ pour le caractère $\psi$. On vérifie facilement les propriétés suivants

\begin{itemize}
\renewcommand{\labelitemi}{$\bullet$}
\item $\gamma_\psi(q\oplus q')=\gamma_\psi(q)\gamma_\psi(q')$ pour toutes formes quadratiques $q,q'$;
\item Si $(V,q)$ est un plan hyperbolique, $\gamma_\psi(q)=1$;
\item $\gamma_\psi(q)^{-1}=\gamma_\psi(-q)$;
\item $\gamma_\psi(q)$ est une racine 8ème de l'unité.
\end{itemize}

\noindent En particulier, on note $\gamma_\psi(N_{E/F})$ la constante de Weil associée à la forme bilinéaire $(x,y)\mapsto \frac{1}{2}Tr_{E/F}(x\overline{y})$ sur $E$. On a alors

$$\gamma_\psi(\lambda N_{E/F})=sgn_{E/F}(\lambda)\gamma_\psi(N_{E/F})$$

\noindent pour tout $\lambda\in F^\times$.

\subsection{Types de modèles de Whittaker}

Fixons une base $(z_i)_{i=\pm 1,\ldots, \pm r}$ telle que $h(z_i,z_j)=\delta_{i,-j}$. Introduisons le sous-groupe parabolique $P_0$ de $G$ des éléments qui conservent le drapeau

$$Ez_r\subset Ez_r\oplus Ez_{r-1}\subset\ldots\subset Ez_r\oplus\ldots\oplus Ez_1$$

\noindent et notons $U_0$ son radical unipotent. Soit $\eta\in E^\times$ tel que $Tr_{E/F}(\eta)=0$. On définit un caractère $\psi_\eta$ de $U(F)$ par la formule suivante

$$\displaystyle\psi_\eta(u)=\psi\bigg(\sum_{i=1}^{r-1} Tr_{E/F}\big(h(z_{-i-1},uz_i)\big)+2Tr_{E/F}\big(\eta h(z_{-1},uz_{-1})\big)\bigg)$$

Définissons un élément $N^-\in\mathfrak{g}(F)$ par $N^-z_i=z_{i-1}$ pour $i=2,\ldots,r$, $N^-z_1=\eta z_{-1}$, $N^-z_{-i}=z_{-i-1}$ pour $i=1,\ldots,r-1$ et $Nz_{-r}=0$. Notons $\overline{U}_0$ le radical unipotent du sous-groupe parabolique opposé à $P_0$. On a alors $N^-\in \overline{\mathfrak{u}}_0(F)$ et $\psi_\eta(e^N)=\psi(<N,N^->)$ pour tout $N\in \mathfrak{u}_0(F)$. \\

Soit $\sigma\in Temp(G)$. On appelle fonctionnelle de Whittaker pour $\sigma$ relativement à $\eta$ toute forme linéaire

$$\ell: E_\sigma\to \mathbb{C}$$

\noindent vérifiant $\ell\circ\sigma(u)=\psi_\eta(u)\ell$ pour tout $u\in U_0(F)$. On note $Hom_{U_0}(\sigma,\psi_\eta)$ l'espace des fonctionnelles de Whittaker relativement à $\sigma$. Comme on le sait bien, cette espace est de dimension au plus $1$. Posons

$$m(\sigma,\eta)=dim\big( Hom_{U_0}(\sigma,\psi_\eta)\big)$$

\noindent Rappelons que l'on a un développement de $\Theta_\sigma$ au voisinage de l'origine de la forme

$$\displaystyle\Theta_\sigma(e^X)=\sum_{\mathcal{O}\in Nil(\mathfrak{g}(F))} c_{\sigma,\mathcal{O}}(1)\hat{j}(\mathcal{O},X)$$

\noindent pour tout $X\in \mathfrak{g}_{reg}(F)$ assez proche de $0$ et où les $c_{\sigma,\mathcal{O}}(1)$ sont des nombres complexes. D'après un résultat de Rodier ([Ro]), on a

$$\mbox{(1)} \;\; m(\sigma,\eta)=c_{\sigma,\mathcal{O}_\eta}(1)$$

\noindent On se propose dans les paragraphes qui suivent de trouver une méthode effective de calcul de ce coefficient.

\subsection{Calcul de germes de Shalika}

Soient $B$ un sous-groupe de Borel de $G$ défini sur $F$ et $T_{qd}$ un sous-tore maximal de $B$. On fixe un élément $X_{qd}\in\mathfrak{t}_{qd}(F)\cap\mathfrak{g}_{reg}(F)$. \\

Choisissons deux éléments $\nu_+,\nu_-\in F^\times$ tels que $sgn_{E/F}(\nu_\pm)=\pm 1$. Soit $Z_0$ un sous-espace de $V$ somme orthogonale de $(d-2)/2$ plans hyperboliques. Pour $\epsilon\in \{\pm\}$, on peut toujours trouver une décomposition orthogonale $V=D^\epsilon_+\oplus^\perp D^\epsilon_-\oplus^\perp Z_0$ où

\begin{itemize}
\item $D^\epsilon_+$ et $D^\epsilon_-$ sont des droites;
\item la restriction de $h$ à $D^\epsilon_+$ est équivalente à $\nu_{\epsilon}N_{E/F}$ et la restriction de $h$ à $D^\epsilon_-$ est équivalente à $-\nu_\epsilon N_{E/F}$.
\end{itemize}

Fixons de telles décompositions. On note $G_0=U(Z_0)$, c'est un groupe unitaire quasi-déployé. Soit $T_{0,qd}$ un tore maximal de $G_0$ inclus dans un sous-groupe de Borel. Fixons un élément $X_{0,qd}$ de $\mathfrak{t}_{0,qd}(F)\cap\mathfrak{g}_{0,reg}(F)$ et deux éléments distincts $a_1,a_2\in Ker(Tr_{E/F})$ qui ne sont pas valeurs propres de $X_{0,qd}$. Notons $X^{\epsilon}\in\mathfrak{g}(F)$ l'élément ainsi défini: $X^{\epsilon}$ agit par multiplication par $a_1$ sur $D^{\epsilon}_+$, par multiplication par $a_2$ sur $D^{\epsilon}_-$ et comme $X_{0,qd}$ sur $Z_0$. On a alors $X^{\epsilon}\in \mathfrak{g}_{reg}(F)$.

\begin{lem}
Pour $\eta\in Ker(Tr_{E/F})\backslash\{0\}$, on a

\begin{enumerate}[(i)]
\item $$\Gamma_{\mathcal{O}_\eta}(X_{qd})=1$$

\item $$\Gamma_{\mathcal{O}_{\eta}}(X^+)-\Gamma_{\mathcal{O}_{\eta}}(X^-)=sgn_{E/F}((a_2-a_1)\eta)$$
\end{enumerate}
\end{lem}

\noindent\ul{Preuve}: Fixons $\eta\in Ker(Tr_{E/F})\backslash\{0\}$. La première assertion de l'énoncé est le lemme 9.3.1 de [B1]. Prouvons le deuxième point. Pour $\epsilon\in\{\pm\}$, notons $T^\epsilon$ le tore maximal de $G$ tel que $X^\epsilon\in\mathfrak{t}^{\epsilon}(F)$. D'après Shelstad ([She]), $\Gamma_{\mathcal{O}_{\eta}}(X^\epsilon)$ vaut $1$ ou $0$ suivant qu'un certain invariant vaut $1$ ou non. Soit $N\in\mathcal{O}_{\eta}$. Cet invariant est noté $inv(X^{\epsilon})inv(T^{\epsilon})/inv(N)$ par Shelstad, c'est un élément de $H^1(F,T^{\epsilon})$. On a une identification naturelle $H^1(F,T^{\epsilon})=\{\pm 1\}\times\{\pm 1\}$ (correspondant à l'indentification naturelle $T^{\epsilon}=Ker(N_{E/F})\times Ker(N_{E/F})\times (R_{E/F}\mathbb{G}_m)^{(d-2)/2}$). Les invariants $inv(X^{\epsilon})inv(T^{\epsilon})$ et $inv(N)$ dépendent chacun du choix d'un épinglage à conjugaison près, ces épinglages devant être opposés l'un à l'autre. Le choix d'un épinglage à conjugaison près revient au choix d'une orbite nilpotente régulière, donc au choix d'un élément de $\big(Ker(Tr_{E/F})\backslash \{0\}\big)/N_{E/F}(E^\times)$ via la paramétrisation 1.3. Deux épinglages $\eta_0,\eta_1\in \big(Ker(Tr_{E/F})\backslash \{0\}\big)/N_{E/F}(E^\times)$ sont opposés si $\eta_0=-\eta_1$. Le lemme X.7 de [W2] calcule l'invariant $inv(X^{\epsilon})inv(T^{\epsilon})$ pour le choix d'un épinglage $\eta_0$ (le $\eta$ de cet référence vaut alors $(-1)^{(d-2)/2}2\eta_0$):

$$inv(X^{\epsilon})inv(T^{\epsilon})=\big(\epsilon sgn_{E/F}((-1)^{(d-2)/2}\eta_0 P'_{X^{\epsilon}}(a_1)),\epsilon sgn_{E/F}((-1)^{d/2}\eta_0 P'_{X^{\epsilon}}(a_2))\big)$$

\noindent où $P_{X^{\epsilon}}$ désigne le polynôme caractéristique de $X^{\epsilon}$ agissant sur $V$. Notons $\lambda_1,\ldots,\lambda_r$, $-\overline{\lambda}_1,\ldots,-\overline{\lambda}_r$, où $r=(d-2)/2$, les valeurs propres de $X_{0,qd}$ agissant sur $V_0$. On a alors

$$\displaystyle P'_{X^{\epsilon}}(a_1)=(a_1-a_2)\times \prod_{i=1}^r (a_1-\lambda_i)(a_1+\overline{\lambda}_i)$$

La norme de $(a_1-\lambda_i)$ est égale à $-(a_1-\lambda_i)(a_1+\overline{\lambda}_i)$. On en déduit que $P'_{X^{\epsilon}}(a_1)\in (-1)^{(d-2)/2}(a_1-a_2)N_{E/F}(E^\times)$. De la même façon, on a $P'_{X^{\epsilon}}(a_2)\in (-1)^{d/2}(a_1-a_2)N_{E/F}(E^\times)$. Par conséquent 

$$inv(X^{\epsilon})inv(T^{\epsilon})=(\epsilon sgn_{E/F}(\eta_0(a_1-a_2)),\epsilon sgn_{E/F}(\eta_0(a_1-a_2)))$$

\noindent Passons au calcul de $inv(N)$. On considère donc, comme on l'a dit, l'épinglage défini par $-\eta_0$. Au choix de cet épinglage est associé, suivant la construction de [She], un élément de $H^1(F,Z(G))$ que l'on voit dans $H^1(F,T^\epsilon)$, via l'application naturelle $H^1(F,Z(G))\to H^1(F,T^\epsilon)$. On a $H^1(F,Z(G))=\{\pm 1\}$ et l'application précédente est l'application diagonale. La construction de Shelstad n'est pas difficile à expliciter ici: on a $inv(N)=1$ si $\eta\in-\eta_0 N_{E/F}(E^\times)$, $-1$ sinon. En d'autre termes, on a $inv(N)=(sgn_{E/F}(-\eta_0/\eta),sgn_{E/F}(-\eta_0/\eta))$. Finalement, on obtient que

$$\Gamma_{\mathcal{O}_{\eta}}(X^{\epsilon})=\left\{
    \begin{array}{ll}
        1 & \mbox{si } sgn_{E/F}(\eta(a_2-a_1))=\epsilon \\
        0 & \mbox{sinon.}
    \end{array}
\right.
$$

\noindent On retrouve alors facilement la deuxième égalité de l'énoncé $\blacksquare$

\subsection{Calcul de fonctions $\hat{j}^G$}

On conserve les objets et notations introduits dans les sections 5.1 à 5.3.

\begin{lem}
Pour $\eta\in Ker(Tr_{E/F})\backslash\{0\}$, on a

\begin{enumerate}[(i)]
\item $$\hat{j}^G(\mathcal{O}_{\eta}, X_{qd})=\frac{w(d)}{2}D^G(X_{qd})^{-1/2}$$

\noindent (où $w(d)$ est le terme défini en 1.3(1))

\item $$\hat{j}^G(\mathcal{O}_{\eta}, X^\epsilon)=\epsilon d^{-1}w(d)\gamma_{\psi}(N_{E/F}) sgn_{E/F}(2(a_1-a_2)\eta) D^G(X^\epsilon)^{-1/2}$$

\noindent pour tout $\epsilon\in\{\pm\}$.

\end{enumerate}
\end{lem}

\noindent\ul{Preuve}: La première égalité est le lemme 18.2.1 de [B1]. Prouvons le deuxième point. Soient $X,Y\in\mathfrak{g}_{reg}(F)$. D'après 5.1(1), pour $\lambda\in F^\times$ assez proche de $0$, on a

$$\displaystyle \hat{j}^G(\lambda X,Y)=\sum_{\mathcal{O}\in Nil(\mathfrak{g})} \Gamma_{\mathcal{O}}(\lambda X) \hat{j}^G(\mathcal{O},Y)$$

\noindent D'après les propriétés d'homogénéité des germes de Shalika, on a aussi

$$\displaystyle \lim\limits_{\substack{\lambda\in F^{\times,2} \\ \lambda\to 0}} \Gamma_{\mathcal{O}}(\lambda X)=\left\{
    \begin{array}{ll}
        \Gamma_{\mathcal{O}}(X) & \mbox{si } \mathcal{O} \mbox{ est régulière} \\
        0 & \mbox{sinon.}
    \end{array}
\right.
$$

\noindent Par conséquent

$$\displaystyle \lim\limits_{\substack{\lambda\in F^{\times,2} \\ \lambda\to 0}}\hat{j}^G(\lambda X,Y)=\sum_{\eta} \Gamma_{\mathcal{O}_{\eta}}(X)\hat{j}^G(\mathcal{O}_{\eta},Y)$$

\noindent où la somme porte sur $\big(Ker(Tr_{E/F})\backslash\{0\}\big)/N(E^\times)$. Appliquons l'égalité précédente à $Y=X^\epsilon$ et $X=X^+,X^-$ et $X_{qd}$. D'après le lemme 5.3.1, on a

$$\mbox{(1)}\;\; \displaystyle \hat{j}^G(\mathcal{O}_{\eta},X^\epsilon)=\frac{1}{2}\lim\limits_{\substack{\lambda\in F^{\times,2} \\ \lambda\to 0}}\bigg(\hat{j}^G(\lambda X_{qd},X^\epsilon)+sgn_{E/F}((a_2-a_1)\eta) \big(\hat{j}^G(\lambda X^+,X^\epsilon)-\hat{j}^G(\lambda X^-,X^\epsilon)\big)\bigg)$$

\noindent pour tout $\epsilon\in\{\pm\}$. Soient $\epsilon\in\{\pm\}$ et $\lambda\in F^{\times, 2}$. La classe de conjugaison de $X^+$ ne rencontre pas $\mathfrak{t}_{qd}(F)$. Par conséquent, d'après 5.1(2), on a $\hat{j}^G(\lambda X_{qd},X^\epsilon)=0$. Il ne reste plus qu'à calculer $\hat{j}^G(\lambda X^+,X^\epsilon)-\hat{j}^G(\lambda X^-,X^\epsilon)$. \\

Supposons dans un premier temps que $d=2$. Posons $G_1=G_2=U(1)$. Alors, $H=G_1\times G_2$ est un groupe endoscopique de $G$. Pour la correspondance endoscopique associée, $X^+$ et $X^-$ correspondent aux classes de conjugaison stable de $Y_1=(a_1,a_2)$ et $Y_2=(a_2,a_1)$. On peut normaliser les facteurs de transfert de sorte que

$$\Delta_{G,H}(Y_i,X^\epsilon)=\epsilon$$

\noindent pour tout $i\in\{1,2\}$ et pour tout $\epsilon\in\{\pm\}$. Grâce à Ngô Bau Chau, la conjecture 1.2 de [W3] est maintenant un théorème. La fonction $\hat{i}^G(X,Z)$ de cette référence vaut $D^G(Z)^{1/2}\hat{j}^G(X,Z)$. Appliquée à $Y=\lambda Y_1$ et $Z=X^\epsilon$ la conjecture 1.2 de [W3] nous donne

\[\begin{aligned}
\mbox{(2)} \;\; \gamma_{\psi}(\mathfrak{g})D^G(X^\epsilon)^{1/2}\big(\Delta_{G,H}(\lambda Y_1,\lambda X^+) \hat{j}^G(\lambda X^+,X^\epsilon)+\Delta_{G,H}(\lambda Y_1,\lambda X^-) \hat{j}^G(\lambda X^-,X^\epsilon)\big) \\
=\epsilon\gamma_{\psi}(\mathfrak{h})\big(\hat{j}^H(\lambda Y_1,Y_1)+\hat{j}^H(\lambda Y_1, Y_2)\big)
\end{aligned}\]

\noindent où $\gamma_{\psi}(\mathfrak{g})$, resp. $\gamma_{\psi}(\mathfrak{h})$, désigne la constante de Weil associée à $\psi$ et à la forme quadratique $X\mapsto \frac{1}{2}Tr_{E/F}(Tr(X^2))$ sur $\mathfrak{g}$, resp. la forme quadratique $(X_1,X_2)\mapsto \frac{1}{2}Tr_{E/F}(Tr(X_1^2)+Tr(X_2^2))$ sur $\mathfrak{h}=\mathfrak{g}_1\times\mathfrak{g}_2$. On a $\Delta_{G,H}(\lambda Y_1,\lambda X^+)=\Delta_{G,H}(Y_1,X^+)=1$ et $\Delta_{G,H}(\lambda Y_1,\lambda X^-)=\Delta_{G,H}(Y_1,X^-)=-1$. Puisque $H$ est un tore, on a $\hat{j}^H(\lambda Y_1,Y_1)=\hat{j}^H(\lambda Y_1,Y_2)=1$ pour $\lambda$ assez petit. Par conséquent, (2) devient

$$\hat{j}^G(\lambda X^+,X^\epsilon)-\hat{j}^G(\lambda X^-,X^\epsilon)=2\epsilon\frac{\gamma_{\psi}(\mathfrak{h})}{\gamma_{\psi}(\mathfrak{g})}D^G(X^\epsilon)^{-1/2}$$

\noindent pour $\lambda\in F^{\times, 2}$ assez petit. Soit $\xi\in F^\times$ tel que $E=F[\xi]$. La forme quadratique sur $\mathfrak{g}$ est équivalente à la somme orthogonale d'un plan hyperbolique et de $(x,y)\mapsto 2x^2+2\xi y^2$. La forme quadratique sur $\mathfrak{h}$ est équivalente à $(x,y)\mapsto \xi x^2+\xi y^2$. La différence de ces deux formes quadratiques est équivalente à la somme d'un plan hyperbolique et de la forme $-2 N_{E/F}$. On a donc

$$\frac{\gamma_{\psi}(\mathfrak{h})}{\gamma_{\psi}(\mathfrak{g})}=\gamma_{\psi}(-2N_{E/F})= sgn_{E/F}(-2)\gamma_\psi(N_{E/F})$$

\noindent On obtient finalement

$$\hat{j}^G(\lambda X^+,X^\epsilon)-\hat{j}^G(\lambda X^-,X^\epsilon)=2\epsilon sgn_{E/F}(-2)\gamma_\psi(N_{E/F})D^G(X^\epsilon)^{-1/2}$$

\noindent pour $\lambda\in F^{\times,2}$ assez petit. \\

Revenons au cas général. Rappelons que pour définir $X^+$, on a fixé trois espaces orthogonaux $D_+$, $D_-$ et $Z_0$. Soient $G'$ le groupe unitaire de $D_+\oplus D_-$ et $T_{0,qd}$ le tore maximal de $G_0=U(Z_0)$ que l'on a fixé pour définir $X^+$ et $X^-$. Alors $M=G'\times T_{0,qd}$ est un Levi de $G$ auquel appartient $X^+$ et $X^-$. Soit ${X'}^+$ (resp. ${X'}^-$) la projection de $X^+$ (resp. $X^-$) sur $\mathfrak{g}'$. Alors la classe de conjugaison de $X^\epsilon$ rencontre $2^{d/2-1}(d/2-1)!$ classes de conjugaisons de $\mathfrak{m}(F)$ et la projection de chacune de ces classes de conjugaison sur $\mathfrak{g}'(F)$ coïncide avec la classe de conjugaison de ${X'}^\epsilon$. D'après 5.1(2), on en déduit que, pour $\lambda$ assez petit,

$$\hat{j}^G(\lambda X^+,X^\epsilon)=2^{d/2-1}(d/2-1)! \hat{j}^{G'}(\lambda {X'}^+,{X'}^\epsilon)D^{G'}({X'}^\epsilon)^{1/2}D^G(X^\epsilon)^{-1/2}$$

\noindent et

$$\hat{j}^G(\lambda X^-,X^\epsilon)=2^{d/2-1}(d/2-1)! \hat{j}^{G'}(\lambda {X'}^-,{X'}^\epsilon)D^{G'}({X'}^\epsilon)^{1/2}D^G(X^\epsilon)^{-1/2}$$

\noindent D'après le résultat obtenu dans le cas $d=2$, on a donc

\[\begin{aligned}
\hat{j}^G(\lambda X^+,X^\epsilon) & -\hat{j}^G(\lambda X^-,X^\epsilon) \\
 & =2^{d/2-1}(d/2-1)!\big(\hat{j}^{G'}(\lambda {X'}^+,{X'}^\epsilon)-\hat{j}^{G'}(\lambda {X'}^-,{X'}^\epsilon)\big)D^{G'}({X'}^\epsilon)^{1/2}D^G(X^\epsilon)^{-1/2} \\
 & =\epsilon 2^{d/2}(d/2-1)!sgn_{E/F}(-2)\gamma_{\psi}(N_{E/F})D^G(X^\epsilon)^{-1/2}
\end{aligned}\]

\noindent pour tout $\lambda\in F^{\times,2}$ assez petit. Réinjectant dans (1), on obtient

$$\hat{j}^G(\mathcal{O}_{\eta},X^\epsilon)=\epsilon 2^{d/2-1}(d/2-1)! sgn_{E/F}(2(a_1-a_2)\eta)\gamma_{\psi}(N_{E/F})D^G(X^\epsilon)^{-1/2}$$

\noindent C'est la deuxième égalité de l'énoncé $\blacksquare$

\subsection{Une formule pour la multiplicité $m(\Sigma,\eta)$}

Soit $\Sigma=\sum_k b_k\sigma_k$ une combinaison linéaire formelle de représentations irréductibles tempérées de $G(F)$. Soit $\eta\in Ker Tr_{E/F}\backslash \{0\}$. En prolongeant par linéarité l'application $\sigma\mapsto m(\sigma,\eta)$, on définit $m(\Sigma,\eta)$. On définit aussi par linéarité le caractère $\Theta_{\Sigma}$ de la représentation virtuelle $\Sigma$.

\begin{lem}
\begin{enumerate}[(i)]

\item On a l'égalité

\[\begin{aligned}
\displaystyle m(\Sigma,\eta)=(2w(d))^{-1}\lim\limits_{\substack{\lambda\in F^{\times, 2} \\ \lambda\to 0}} & \bigg(\Theta_\Sigma(e^{\lambda X_{qd}})D^G(e^{\lambda X_{qd}})^{1/2} \\
 & +d\; sgn_{E/F}(2(a_1-a_2)\eta)\gamma_{\psi}(N_{E/F})^{-1}\Theta_\Sigma(e^{\lambda X^+})D^G(e^{\lambda X^+})^{1/2}\bigg)
\end{aligned}\]

\item Si $\Theta_\Sigma$ est stable alors,

$$\displaystyle \lim\limits_{\substack{\lambda\in F^{\times,2}\\ \lambda\to 0}}\Theta_\Sigma(e^{\lambda X^+})D^G(e^{\lambda X^+})^{1/2}=0$$
\end{enumerate} 
\end{lem}

\noindent\ul{Preuve}: Le premier point est une conséquence facile de 5.2(1), du lemme 5.4.1 et des propriétés d'homogénéité des fonctions $\hat{j}^G(\mathcal{O},.)$. Supposons que $\Theta_\Sigma$ soit stable. Comme $X^+$ et $X^-$ sont stablement conjugués, pour tout $\lambda\in F^\times$, on a $\Theta_\Sigma(e^{\lambda X^+})=\Theta_\Sigma(e^{\lambda X^-})$ et $D^G(e^{\lambda X^+})=D^G(e^{\lambda X^-})$. Par conséquent dans la limite en (ii), on peut remplacer $X^+$ par $X^-$. D'après le lemme 5.4.1, la limite existe et selon qu'on la calcule avec $X^+$ ou $X^-$, on obtient deux résultats opposés. La limite est donc nulle $\blacksquare$

\subsection{Transfert et multiplicité $m(\Sigma,\eta)$}

On conserve les notations et hypothèses précédentes: $(V,h)$ est un espace hermitien de dimension $d$ paire et de groupe unitaire $G$ quasi-déployé. Soient $(V_+,h_+)$ et $(V_-,h_-)$ deux espaces hermitiens de dimensions respectives $d_+$, $d_-$ et de groupes unitaires respectifs $G_+$ et $G_-$. On fait les hypothèses suivantes

\vspace{2mm}

\begin{itemize}
\renewcommand{\labelitemi}{$\bullet$}
\item $d=d_++d_-$;
\item $G_+$ et $G_-$ sont quasi-déployés.
\end{itemize}

\vspace{2mm}

On peut alors considérer $G_+\times G_-$ comme un groupe endoscopique de $G$. Pour fixer la donnée endoscopique, on choisit un couple $(\mu_+,\mu_-)$ de caractères continus de $E^\times$ dont la restriction à $F^\times$ coïncide avec $(sgn_{E/F}^{d_-},sgn_{E/F}^{d_+})$. On normalise alors les facteurs de transfert comme en 3.1. Soient $\Sigma$, $\Sigma_+$, $\Sigma_-$ des combinaisons linéaires formelles d'éléments de $Temp(G)$, $Temp(G_+)$ et $Temp(G_-)$ respectivement. On suppose que

\begin{itemize}
\renewcommand{\labelitemi}{$\bullet$}
\item $\Theta_{\Sigma_+}$ et $\Theta_{\Sigma_-}$ sont stables;
\item $\Theta_{\Sigma}$ est un transfert de $\Theta_{\Sigma_+}\times \Theta_{\Sigma_-}$.
\end{itemize}

\begin{prop}
Soit $\eta\in Ker(Tr_{E/F})-\{0\}/N(E^\times)$. Sous ces hypothèses, on a l'égalité

$$m(\Sigma,\eta)=\left\{
    \begin{array}{ll}
        c_{\Sigma_+}(1)c_{\Sigma_-}(1) & \mbox{si } d_+ \mbox{ et } d_- \mbox{ sont pairs} \\
        sgn_{E/F}(2\nu_0\eta\delta)\gamma_{\psi}(N_{E/F})^{-1}c_{\Sigma_+}(1)c_{\Sigma_-}(1) & \mbox{si } d_+ \mbox{ et } d_- \mbox{ sont impairs}
    \end{array}
\right.
$$
\end{prop}

\noindent\ul{Preuve}: Fixons deux éléments distincts $a_1,a_2\in Ker Tr_{E/F}\backslash \{0\}$. Introduisons des éléments $X_+,X_{qd}\in\mathfrak{g}_{reg}(F)$ comme en 5.3. D'après le lemme 5.5.1, on a alors

\[\begin{aligned}
\mbox{(1)}\;\;\; \displaystyle & m(\Sigma,\eta)=(2w(d))^{-1} \\
 & \lim\limits_{\substack{\lambda\in F^{\times,2} \\ \lambda\to 0}} \big(2\Theta_\Sigma(e^{\lambda X_{qd}})D^G(e^{\lambda X_{qd}})^{1/2}+d \; sgn_{E/F}(2A\eta)\gamma_\psi(N_{E/F})^{-1} \Theta_\Sigma(e^{\lambda X_+})D^G(e^{\lambda X_+})^{1/2} \big)
\end{aligned}\]

\noindent où on a posé $A=a_1-a_2$. Il existe une famille $(\alpha_i)_{i\in I}$ d'éléments de $E$ telle que, pour tout $\lambda\in F^\times$ assez petit, la classe de conjugaison stable de $e^{\lambda X_{qd}}$ soit paramétrée par $\xi_\lambda=(I,(F_{\pm i})_{i\in I}, (F_i)_{i\in I}, (y_i(\lambda))_{i\in I})$ où

\begin{itemize}
\item $F_{\pm i}=E$, $F_i=E\times E$;
\item $y_i(\lambda)=(e^{\lambda \alpha_i},e^{-\lambda\overline{\alpha}_i})$
\end{itemize}

\noindent pour tout $i\in I$. On a $C(\xi_\lambda)=\{1\}$, c'est à-dire que la classe de conjugaison de $e^{\lambda X_{qd}}$ est entièrement déterminée par $\xi_\lambda$. Pour tout $I'\subset I$, posons $\xi_\lambda(I')=(I',(F_{\pm i})_{i\in I'},(F_i)_{i\in I'}, (y_i(\lambda))_{i\in I'})$. Puisque $\Theta_\Sigma$ est un transfert de $\Theta_{\Sigma_+}\times \Theta_{\Sigma_-}$, on a, pour tout $\lambda$ assez petit,

\[\begin{aligned}
\mbox{(2)}\;\;\; \displaystyle \Theta_\Sigma(e^{\lambda X_{qd}})D^G(e^{\lambda X_{qd}})^{1/2}=\sum_{I_1,I_2} & \Theta_{\Sigma_+}(\xi_\lambda(I_1))D^{G_+}(\xi_\lambda(I_1))^{1/2}\Theta_{\Sigma_-}(\xi_\lambda(I_2)) \\
 & D^{G_-}(\xi_\lambda(I_2))^{1/2} \Delta_{\mu_+,\mu_-}(\xi_\lambda(I_1),\xi_\lambda(I_2),1)
\end{aligned}\]

\noindent où la somme porte sur les couples $(I_1,I_2)$ vérifiant $I_1\sqcup I_2=I$, $2|I_1|=d_+$ et $2|I_2|=d_-$. Fixons un tel couple. D'après 4.5(1), on a

$$\lim\limits_{\lambda\to 0}\Theta_{\Sigma_+}(\xi_\lambda(I_1))D^{G_+}(\xi_\lambda(I_1))^{1/2}=w(d_+)c_{\Sigma_+}(1)$$

$$\lim\limits_{\lambda\to 0}\Theta_{\Sigma_-}(\xi_\lambda(I_2))D^{G_-}(\xi_\lambda(I_2))^{1/2}=w(d_-)c_{\Sigma_-}(1)$$

\noindent D'après les points (A)(i) et (A)(ii) du lemme 2.3.1 (appliqués plusieurs fois), on a, pour $\lambda$ assez petit,

$$\Delta_{\mu_+,\mu_-}(\xi_\lambda(I_1),\xi_\lambda(I_2),1)=\Delta_{\mu_+,\mu_-}(\emptyset,\emptyset,1)= 1$$

\noindent Par passage à la limite dans (2), on en déduit que

$$\mbox{(3)}\;\;\; \lim\limits_{\lambda\in F^\times, \lambda\to 0}\Theta_\Sigma(e^{\lambda X_{qd}})D^G(e^{\lambda X_{qd}})^{1/2}=\sum_{I_1,I_2}w(d_+)w(d_-)c_{\Sigma_+}(1)c_{\Sigma_-}(1)$$

\noindent la somme portant sur les mêmes couples que (2). Si $d_+$ et $d_-$ sont impairs, cette somme est vide. Si $d_+$ et $d_-$ sont pairs, cette somme comporte $\binom{|I|}{d_+/2}$ termes. On vérifie que

$$\binom{|I|}{d_+/2}w(d_+)w(d_-)=w(d)$$

\noindent Par conséquent (3) devient

$$\mbox{(4)}\;\;\; \lim\limits_{\lambda\to 0}\Theta_\Sigma(e^{\lambda X_{qd}})D^G(e^{\lambda X_{qd}})^{1/2}= \left\{
 \begin{array}{ll}
        w(d)c_{\Sigma_+}(1)c_{\Sigma_-}(1) & \mbox{si } d_+ \mbox{ et } d_- \mbox{ sont pairs} \\
        0 & \mbox{si } d_+ \mbox{ et } d_- \mbox{ sont impairs}
    \end{array}
\right.
$$

 Il existe une famille $(\beta_j)_{j\in J}$ d'éléments de $E$ avec $\{1,2\}\subset J$, $\beta_1=a_1$, $\beta_2=a_2$ telle que, pour $\lambda$ assez petit, la classe de conjugaison stable de $e^{\lambda X_+}$ soit paramétrée par $\zeta_\lambda=(J,(F_{\pm j})_{j\in J},(F_j)_{j\in J},(y_j(\lambda))_{j\in J})$, où
 
\begin{itemize}
\item $F_{\pm 1}=F_{\pm 2}=F$, $F_1=F_2=E$, $y_1(\lambda)=e^{\lambda a_1}$, $y_2(\lambda)=e^{\lambda a_2}$;
\item $F_{\pm j}=E$, $F_j=E\times E$ $y_j(\lambda)=(e^{\lambda \beta_j},e^{-\lambda\overline{\beta}_j})$ pour tout $j\in J\backslash\{1,2\}$.
\end{itemize}

\noindent On a alors $C(\zeta_\lambda)=F_{\pm 1}^\times/N(F_1^\times)\times F_{\pm 2}^\times/N(F_2^\times)=F^\times/N(E^\times)\times F^\times/N(E^\times)$. La classe de conjugaison de $e^{\lambda X_+}$ est alors déterminée par $c=(\overline{1},-\overline{1})\in F^\times/N(E^\times)\times F^\times/N(E^\times)$ (où pour $x\in F^\times$, on désigne par $\overline{x}$ l'image de $x$ dans $F^\times/N(E^\times)$). Puisque $\Theta_\Sigma$ est un transfert de $\Theta_{\Sigma_+}\times \Theta_{\Sigma_-}$, on a, pour tout $\lambda$ assez petit,

\[\begin{aligned}
\mbox{(5)}\;\;\; \displaystyle \Theta_\Sigma(e^{\lambda X_+})D^G(e^{\lambda X_+})^{1/2}=\sum_{J_1,J_2} & \Theta_{\Sigma_+}(\zeta_\lambda(J_1))D^{G_+}(\zeta_\lambda(J_1))^{1/2} \Theta_{\Sigma_-}(\zeta_\lambda(J_2)) \\
 & D^{G_-}(\zeta_\lambda(J_2))^{1/2} \Delta_{\mu_+,\mu_-}(\zeta_\lambda(J_1),\zeta_\lambda(J_2),c)
\end{aligned}\]

\noindent où la somme porte sur les couples $(J_1,J_2)$ vérifiant $J_1\sqcup J_2=J$, $d_{J_1}=d_+$ et $d_{J_2}=d_-$. Fixons un tel couple. Supposons dans un premier temps que $d_+$ et $d_-$ sont pairs. Alors, on a $\{1,2\}\subset J_1$ ou $\{1,2\}\subset J_2$. Si $\{1,2\}\subset J_1$, alors, en remplaçant $G$ par $G_+$ dans les constructions, on peut trouver un élément $X_+^{G_+}\in\mathfrak{g}_{+,reg}(F)$ analogue à $X_+$ tel que pour tout $\lambda\in F^\times$ assez petit, la classe de conjugaison stable de $e^{\lambda X^{G_+}_+}$ soit paramétrée par $\zeta_\lambda(J_1)$. D'après le lemme 5.5.1, $\Theta_{\Sigma_+}(\zeta_\lambda(J_1))D^{G_+}(\zeta_\lambda(J_1))^{1/2}$ tend vers $0$ pour $\lambda\in F^{\times,2}$ qui tend vers $0$. D'après 4.5(1), $\Theta_{\Sigma_-}(\zeta_\lambda(J_2))D^{G_-}(\zeta_\lambda(J_2))^{1/2}$ admet aussi une limite lorsque $\lambda$ tend vers $0$, tandis que $\Delta_{\mu_+,\mu_-}(\zeta_\lambda(J_1),\zeta_\lambda(J_2),c)$ reste borné (car c'est un nombre complexe de module $1$). On en déduit que le terme indexé par $(I_1,I_2)$ dans la somme (5) tend vers $0$ pour $\lambda\in F^{\times,2}$ qui tend vers $0$. Le même résultat, avec la même preuve, est valable pour les couples $(J_1,J_2)$ avec $\{1,2\}\subset J_2$. Par conséquent, si $d_+$ et $d_-$ sont pairs, on a

$$\mbox{(6)}\;\;\; \displaystyle \lim\limits_{\lambda\in F^{\times,2},\lambda\to 0} \Theta_\Sigma(e^{\lambda X_+})D^G(e^{\lambda X_+})^{1/2}=0$$

Supposons maintenant que $d_+$ et $d_-$ sont impairs et fixons un couple $(J_1,J_2)$ comme dans la somme (5). On a alors $1\in J_1$ et $2\in J_2$ ou $2\in J_1$ et $1\in J_2$. Dans les deux cas, d'après 4.5(1), on a

$$\mbox{(7)}\;\;\; \lim\limits_{\lambda\to 0} \Theta_{\Sigma_+}(\zeta_\lambda(J_1))D^{G_+}(\zeta_\lambda(J_1))^{1/2}=w(d_+)c_{\Sigma_+}(1)$$

$$\mbox{(8)}\;\;\; \lim\limits_{\lambda\to 0} \Theta_{\Sigma_-}(\zeta_\lambda(J_2))D^{G_-}(\zeta_\lambda(J_2))^{1/2}=w(d_-)c_{\Sigma_-}(1)$$

\noindent Toujours d'après les points (A)(i) et (A)(ii) du lemme 2.3.1, pour $\lambda$ assez petit, on a l'égalité

\[\begin{aligned}
\mbox{(9)}\;\;\; \displaystyle & \Delta_{\mu_+,\mu_-}(\zeta_\lambda(J_1),\zeta_\lambda(J_2),c)= \\
 & \left\{
 \begin{array}{ll}
        \Delta_{\mu_+,\mu_-}(\zeta_\lambda(1),\zeta_\lambda(2),c) & \mbox{si } 1\in J_1 \mbox{ et } 2\in J_2 \\
        \Delta_{\mu_+,\mu_-}(\zeta_\lambda(2),\zeta_\lambda(1),c) & \mbox{si } 2\in J_1 \mbox{ et } 1\in J_2
    \end{array}
\right.
\end{aligned}\]

\noindent où on a posé $\zeta_\lambda(i)=\zeta_\lambda(\{i\})$ pour $i\in\{1,2\}$. On calcule

\[\begin{aligned}
\mbox{(10)}\;\;\; \displaystyle & \Delta_{\mu_+,\mu_-}(\zeta_\lambda(1),\zeta_\lambda(2),c)=\Delta_{\mu_+,\mu_-}(\zeta_\lambda(2),\zeta_\lambda(1),c) \\
 & =\mu_+(-1-e^{\lambda a_2})\mu_-(-1-e^{\lambda a_1}) sgn_{E/F}\big(-\nu_0\delta^{-3}(e^{\lambda a_1}-e^{\lambda a_2})(1+e^{\lambda a_1})^{-1}(1+e^{\lambda a_2})^{-1}\big)
\end{aligned}\]

\noindent On peut remplacer le $\delta^{-3}$ de l'égalité précédente par $\delta$ (car $\delta^4$ est une norme de $E^\times$). Lorsque $\lambda$ tend vers $0$, $\frac{e^{\lambda a_1}-e^{\lambda a_2}}{\lambda(a_1-a_2)}$ tend vers $1$. Par conséquent pour $\lambda\in F^{\times,2}$ assez proche de $0$, on a

\[\begin{aligned}
sgn_{E/F}(\delta(e^{\lambda a_1}-e^{\lambda a_2})) & =sgn_{E/F}(\delta\lambda(a_1-a_2)) \\
 & =sgn_{E/F}(\delta(a_1-a_2))
\end{aligned}\]

\noindent Les autres termes de l'expression (10) admettent des limites lorsque $\lambda$ tend vers $0$ qui se calculent aisément. On obtient

\[\begin{aligned}
\mbox{(11)}\;\;\; \displaystyle & \lim\limits_{\lambda\in F^{\times, 2},\lambda\to 0}\Delta_{\mu_+,\mu_-}(\zeta_\lambda(1),\zeta_\lambda(2),c)= \lim\limits_{\lambda\in F^{\times, 2},\lambda\to 0}\Delta_{\mu_+,\mu_-}(\zeta_\lambda(2),\zeta_\lambda(1),c)= \\
 & \mu_+\mu_-(-2)sgn_{E/F}(\nu_0\delta(a_2-a_1))= sgn_{E/F}(-\nu_0\delta A)
\end{aligned}\]

\noindent où à la dernière égalité, on a utilisé le fait que le caractère $\mu_+\mu_-$ est trivial sur $F^\times$. De (5), (7), (8), (9) et (11), on déduit que

$$\displaystyle \lim\limits_{\lambda\in F^{\times,2},\lambda\to 0} \Theta_\Sigma(e^{\lambda X_+}) D^G(e^{\lambda X_+})^{1/2}=\sum_{J_1,J_2} w(d_+)w(d_-)sgn_{E/F}(-\nu_0\delta A)c_{\Sigma_+}(1)c_{\Sigma_-}(1)$$

\noindent où on somme sur les mêmes couples $(J_1,J_2)$ qu'en (5). Il est facile de compter le nombre de tels couples: il y en a $2\binom{|J|-2}{(d_+-1)/2}=2\binom{d/2-1}{(d_+-1)/2}$. On vérifie aisément que l'on a l'égalité

$$2\binom{d/2-1}{(d_+-1)/2}w(d_+)w(d_-)=2d^{-1}w(d)$$

\noindent On obtient finalement que, pour $d_+$ et $d_-$ impairs, on a

$$\mbox{(12)}\;\;\; \displaystyle \lim\limits_{\lambda\in F^{\times,2},\lambda\to 0} \Theta_\Sigma(e^{\lambda X_+}) D^G(e^{\lambda X_+})^{1/2}=2d^{-1}w(d)sgn_{E/F}(-\nu_0\delta A)c_{\Sigma_+}(1)c_{\Sigma_-}(1)$$

\noindent L'égalité de l'énoncé est maintenant conséquence facile de (1), (4), (6) et (12) $\blacksquare$

\section{Transfert endoscopique}

\subsection{Définition d'une multiplicité stable}

Soient $(V,h)$ et $(V',h')$ deux espaces hermitiens de dimensions respectives $d$ et $d'$ (on ne fait ici aucune hypothèse sur leurs parités) et de groupes unitaires respectifs $G$ et $G'$. On suppose que

\vspace{3mm}

\begin{itemize}
\renewcommand{\labelitemi}{$\bullet$}
\item $G$ et $G'$ sont quasi-déployés.
\end{itemize}

\vspace{3mm}

 Soient $\Sigma$ et $\Sigma'$ des représentations virtuelles tempérées de $G(F)$ et $G'(F)$ respectivement. C'est-à-dire que $\Sigma=\sum_k a_k \sigma_k$ et $\Sigma'=\sum_k b_k \sigma_k'$ sont des combinaisons linéaires formelles finies où les $a_k$ et les $b_k$ sont des nombres complexes, les $\sigma_k$ appartiennent à $Temp(G(F))$ et les $\sigma'_k$ appartiennent à $Temp(G'(F))$. On définit par linéarité leurs caractères $\Theta_\Sigma$ et $\Theta_{\Sigma'}$. On suppose 
 
\vspace{3mm}

\begin{itemize} 
\renewcommand{\labelitemi}{$\bullet$}
\item $\Theta_\Sigma$ et $\Theta_{\Sigma'}$ sont stables (c'est-à-dire constantes sur les classes de conjugaison stable d'éléments fortement réguliers).
\end{itemize}

\vspace{3mm}

 \indent Pour tout $\underline{d}\leqslant min(d,d')$, on peut trouver une décomposition orthogonale $V=W\oplus V_\natural$ où $dim(W)=\underline{d}$ et $G_\natural$ le groupe unitaire de $V_\natural$ est quasi-déployé. Notons $G_W$ le groupe unitaire de $W$. Soit $\xi\in\Xi^*_{\underline{d},reg}$ et $x\in G_{W,reg}(F)$ un élément dont la classe de conjugaison stable est paramétrée par $\xi$ (il en existe toujours). Alors $G_x=T\times G_\natural$ où $T$ est un tore maximal de $G_W$. On peut donc définir le coefficient $c_\Sigma(x)=\sum_k a_k c_{\sigma_k}(x)$. Puisque $\Theta_\Sigma$ est stable, ce terme ne dépend pas du choix de $x$ ni de celui de la décomposition orthogonale. On le note $c_\Sigma(\xi)$. On définit de même $c_{\Sigma'}(\xi)$. Soit $\mathcal{D}(d,d')$ l'ensemble des entiers naturels $\underline{d}$ qui vérifient

\begin{itemize}
\renewcommand{\labelitemi}{$\bullet$}
\item $\underline{d}\leqslant min(d,d')$;
\item $\underline{d}\equiv d \; \mbox{mod } 2$ ou $\underline{d}\equiv  d' \; \mbox{mod } 2$.
\end{itemize}

\noindent Introduisons l'espace

$$\displaystyle \Xi^*(d,d')=\bigsqcup_{\underline{d}\in \mathcal{D}(d,d')} \Xi^*_{\underline{d},reg}$$

\noindent Soit $\mu$ un caractère continu de $E^\times$ dont la restriction à $F^\times$ coïncide avec $sgn_{E/F}^{d+d'+1}$. On pose alors

$$\mbox{(1)}\;\;\; \displaystyle S_\mu(\Sigma,\Sigma')=\lim\limits_{s\to 0^+}\int_{\Xi^*(d,d')} \mu(\delta^{-d_\xi}P_\xi(1))|C(\xi)| c_{\Sigma}(\xi) c_{\Sigma'}(\xi) D^d(\xi)^{1/2}D^{d'}(\xi)^{1/2} \Delta(\xi)^{-1/2+s} d\xi$$

\subsection{Transfert et multiplicité}

\noindent Donnons nous

\vspace{3mm}

\begin{itemize}
\renewcommand{\labelitemi}{$\bullet$}

\item $(V,h)$ et $(V',h')$ deux espaces hermitiens de dimensions respectives $d$, $d'$ et de groupes unitaires respectifs $G$ et $G'$;

\item $(V_+,h_+)$, $(V_-,h_-)$, $(V'_+,h'_+)$ et $(V'_-,h'_-)$ quatre espaces hermitiens de dimensions respectives $d_+$, $d_-$, $d'_+$, $d'_-$ et de groupes unitaires respectifs $G_+$, $G_-$, $G'_+$ et $G'_-$;

\item $\mu_+$, $\mu_-$, $\mu'_+$ et $\mu'_-$ des caractères continus de $E^\times$.
\end{itemize}

\vspace{3mm}

\noindent On suppose que

\vspace{3mm}

\begin{itemize}
\renewcommand{\labelitemi}{$\bullet$}
\item $d$ est pair et $d'$ est impair;
\item $(V,h)$ et $(V',h')$ vérifient l'hypothèse 4.1(1);
\item $(V_+,h_+)$, $(V_-,h_-)$, $(V'_+,h'_+)$ et $(V'_-,h'_-)$ vérifient l'hypothèse \textbf{(QD)};
\item $d_++d_-=d$ et $d'_++d'_-=d'$;
\item $\mu_{+,|F^\times}=sgn_{E/F}^{d_-}$, $\mu_{-,|F^\times}=sgn_{E/F}^{d_+}$, $\mu'_{+,|F^\times}=sgn_{E/F}^{d'_-}$ et $\mu'_{-,|F^\times}=sgn_{E/F}^{d'_+}$.
\end{itemize}

\vspace{3mm}

\noindent D'après 3.1, on peut donc considérer $G_+\times G_-$ comme un groupe endoscopique de $G$ et $G'_+\times G'_-$ comme un groupe endoscopique de $G'$ (les caractères $\mu_+$, $\mu_-$, $\mu'_+$ et $\mu'_-$ permettant de fixer les données endoscopiques). On normalise les facteurs de transfert comme en 3.1. Considérons

\vspace{3mm}

\begin{itemize}
\renewcommand{\labelitemi}{$\bullet$}
\item $\Sigma$ et $\Sigma'$ des représentations virtuelles tempérées de $G$ et $G'$ respectivement;
\item $\Sigma_+$, $\Sigma_-$, $\Sigma'_+$ et $\Sigma'_-$ des représentations virtuelles tempérées de $G_+$, $G_-$, $G'_+$ et $G'_-$ respectivement.
\end{itemize}

\vspace{3mm}

\noindent On suppose que

\vspace{3mm}

\begin{itemize}
\renewcommand{\labelitemi}{$\bullet$}

\item les caractères $\Theta_{\Sigma_+}$, $\Theta_{\Sigma_-}$, $\Theta_{\Sigma'_+}$, $\Theta_{\Sigma'_-}$ sont stables;

\item $\Theta_\Sigma$ est un transfert de $\Theta_{\Sigma_+}\times \Theta_{\Sigma_-}$ et $\Theta_{\Sigma'}$ est un transfert de $\Theta_{\Sigma'_+}\times \Theta_{\Sigma'_-}$.
\end{itemize}

\vspace{3mm}

\noindent En prolongeant par bilinéarité l'application $(\sigma,\sigma')\mapsto m(\sigma,\sigma')$, définie en 4.1, on définit $m(\Sigma,\Sigma')$. Enfin, on pose

$$\mu(G)= \left\{
    \begin{array}{ll}
        1 & \mbox{si } G \mbox{ est quasi-déployé} \\
        -1 & \mbox{sinon}
    \end{array}
\right.
$$

\noindent c'est-à-dire qu'avec la notation introduite en 3.1, on a $\mu(G)=\mu(V,h)$.

\begin{prop}
Sous les hypothèses précédentes, on a l'égalité

$$\displaystyle m(\Sigma,\Sigma')=\frac{1}{2}\left(S_{\mu_+\mu'_+}(\Sigma_+,\Sigma'_+)S_{\mu_-\mu'_-}(\Sigma_-,\Sigma'_-)+\mu(G)S_{\mu_+\mu'_-}(\Sigma_+,\Sigma'_-)S_{\mu_-\mu'_+}(\Sigma_-,\Sigma'_+)\right)$$
\end{prop}

\noindent\ul{Preuve}: Exprimons $m(\Sigma,\Sigma')$ grâce à la formule 4.2(1). Puisque $\mathcal{C}(h,h')$ est un revêtement de $\Xi^*(h,h')$ qui lui-même est un ouvert de $\Xi^*(d,d')$, on peut réécrire cette formule sous la forme

$$\displaystyle m(\Sigma,\Sigma')=\lim\limits_{s\to 0^+} \int_{\Xi^*(d,d')} f_1(\xi)\Delta(\xi)^{s-1/2} d\xi$$

\noindent Le membre de droite s'exprime comme limite d'une intégrale sur $\Xi^*(d_+,d'_+)\times \Xi^*(d_-,d'_-)\sqcup \Xi^*(d_+,d'_-)\times\Xi^*(d_-,d'_+)$. L'application $(\xi_1,\xi_2)\mapsto \xi_1\sqcup\xi_2$ définit un revêtement qui préserve localement les mesures d'un ouvert de complémentaire de mesure nulle de cet espace au dessus d'un ouvert de $\Xi^*(d,d')$. On a l'égalité $\Delta(\xi_1\sqcup\xi_2)=\Delta(\xi_1)\Delta(\xi_2)$, ce qui permet de réécrire le membre de droite sous la forme

$$\displaystyle \lim\limits_{s\to 0^+} \int_{\Xi^*(d,d')} f_2(\xi)\Delta(\xi)^{s-1/2} d\xi$$

\noindent Pour établir le lemme, il suffit donc de montrer l'égalité $f_1(\xi)=f_2(\xi)$ pour tout $\xi\in \Xi^*(d,d')$. Soit donc $\xi=(I,(F_{\pm i})_{i\in I},(F_i)_{i\in I}, (y_i)_{i\in I})\in \Xi^*(d,d')$. Posons $\underline{d}=d_\xi$. On construit comme en 4.2, un espace hermitien $(V_{\underline{d}},h_{\underline{d}})$ qui se plonge à la fois dans $(V,h)$ et $(V',h')$. Notons $G_{\underline{d}}$ son groupe unitaire. On dispose donc de deux plongements, bien définis à conjugaison près, de $G_{\underline{d}}$ dans $G$ et $G'$. On a alors

$$\mbox{(1)}\;\;\; \displaystyle f_1(\xi)=\sum_{c\in C(\xi)^\epsilon} D^d(\xi)^{1/2}c_\Sigma(x(\xi,c))D^{d'}(\xi)^{1/2}c_{\Sigma'}(x(\xi,c))$$

\noindent pour un certain $\epsilon\in\{\pm 1\}$. Soit $c\in C(\xi)^\epsilon$. Notons $V_\natural$ l'orthogonal de $V_{\underline{d}}$ dans $V$, $d_\natural$ sa dimension et $G_\natural$ son groupe unitaire. D'après la construction 4.2, $G_\natural$ est quasi-déployé. Soit $T_\natural$ un sous-tore maximal d'un sous-groupe de Borel de $G_\natural$ tous deux définis sur $F$. Fixons $Y_\natural\in \mathfrak{g}_{\natural,reg}(F)\cap\mathfrak{t}_\natural(F)$ qui n'admet pas $0$ comme valeur propre dans $V_\natural$. D'après 4.5(1), on a

$$\mbox{(2)}\;\;\; \displaystyle D^d(\xi)^{1/2}c_\Sigma(x(\xi,c))=w(d_\natural)^{-1}\lim\limits_{\lambda\in F^\times, \lambda\to 0} \Theta_\Sigma(x(\xi,c)e^{\lambda Y_\natural}) D^G(x(\xi,c)e^{\lambda Y_\natural})^{1/2}$$

\noindent Pour tout $\lambda\in F^\times$ assez petit, notons $\zeta_\lambda\in \Xi_{d_\natural,reg}$ l'élément qui paramètre la classe de conjugaison stable de $e^{\lambda Y_\natural}$. On distingue alors deux cas. \\

\noindent\textbullet Supposons tout d'abord que $\underline{d}$ est pair. Il existe alors une famille $(a_j)_{j\in J}$ d'éléments non nuls de $E$ telle que $\zeta_\lambda=(J,(F_{\pm j})_{j\in J}, (F_j)_{j\in J}, (y_j(\lambda))_{j\in J})$ pour tout $\lambda\in F^\times$ assez petit, où:
 
\begin{itemize}
\item $F_{\pm j}=E$ et $F_j=E\times E$;
\item $y_j(\lambda)=(e^{\lambda a_j},e^{-\lambda \overline{a}_j})$
\end{itemize}

\noindent pour tout $j\in J$. On a alors $C(\zeta_\lambda)=\{1\}$, i.e. $\zeta_\lambda$ détermine aussi la classe de conjugaison stable de $e^{\lambda Y_\natural}$. La classe de conjugaison de $x(\xi,c)e^{\lambda Y_\natural}$ est alors paramétrée par $\xi\sqcup\zeta_\lambda$ et $c\in C(\xi\sqcup\zeta_\lambda)=C(\xi)$. Pour tout $I'\subset I$, on pose $\xi(I')=(I',(F_{\pm i})_{i\in I'},(F_i)_{i\in I'},(y_i)_{i\in I'})$ et $d_{I'}=d_{\xi(I')}$. On définit de même $\zeta_\lambda(J')$ et $d_{J'}$ pour tout $J'\subset J$. Puisque $\Theta_\Sigma$ est un transfert de $\Theta_{\Sigma_+}\times \Theta_{\Sigma_-}$, on a

\[\begin{aligned}
\mbox{(3)}\;\;\; & \Theta_{\Sigma}(x(\xi,c)e^{\lambda Y_\natural})D^G(x(\xi,c)e^{\lambda Y_\natural})^{1/2}= \\
 & \sum_{\substack{I_1, I_2 \\ J_1, J_2}} \Theta_{\Sigma_+}(\xi(I_1)\sqcup \zeta_\lambda(J_1))D^{G_+}(\xi(I_1)\sqcup \zeta_\lambda(J_1))^{1/2} \Theta_{\Sigma_-}(\xi(I_2)\sqcup \zeta_\lambda(J_2)) \\
 & D^{G_-}(\xi(I_2)\sqcup \zeta_\lambda(J_2))^{1/2} \Delta_{\mu_+,\mu_-}(\xi(I_1)\sqcup \zeta_\lambda(J_1),\xi(I_2)\sqcup \zeta_\lambda(J_2),c)
\end{aligned}\]

\noindent où la somme porte sur les quadruplets $(I_1,I_2,J_1,J_2)$ satisfaisant $I_1\sqcup I_2=I$, $J_1\sqcup J_2=J$, $d_{I_1}+2|J_1|=d_+$ et $d_{I_2}+2|J_2|=d_-$. Fixons un tel quadruplet. D'après 4.5(1), on a

$$\mbox{(4)}\;\;\; \displaystyle \lim\limits_{\lambda\to 0} \Theta_{\Sigma_+}(\xi(I_1)\sqcup \zeta_\lambda(J_1))D^{G_+}(\xi(I_1)\sqcup \zeta_\lambda(J_1))^{1/2}=c_{\Sigma_+}(\xi(I_1))D^{d_+}(\xi(I_1))^{1/2}$$

$$\mbox{(5)}\;\;\; \displaystyle \lim\limits_{\lambda\to 0} \Theta_{\Sigma_-}(\xi(I_2)\sqcup \zeta_\lambda(J_2))D^{G_-}(\xi(I_2)\sqcup \zeta_\lambda(J_2))^{1/2}=c_{\Sigma_-}(\xi(I_2))D^{d_-}(\xi(I_2))^{1/2}$$

\noindent D'après les points (A)(i) et (A)(ii) du lemme 2.3.1 (appliqués plusieurs fois), pour $\lambda$ assez proche de $0$, on a

$$\mbox{(6)}\;\;\; \Delta_{\mu_+,\mu_-}(\xi(I_1)\sqcup \zeta_\lambda(J_1),\xi(I_2)\sqcup \zeta_\lambda(J_2),c)=\Delta_{\mu_+,\mu_-}(\xi(I_1),\xi(I_2),c)$$

\noindent Pour $I_1,I_2\subset I$, tels que $d_{I_1}\leqslant d_+$ et $d_{I_2}\leqslant d_-$, posons

$$B(I_1,I_2)=c_{\Sigma_+}(\xi(I_1))D^{d_+}(\xi(I_1))^{1/2}c_{\Sigma_-}(\xi(I_2)) D^{d_-}(\xi(I_2))^{1/2}$$

\noindent De (2), (3), (4), (5) et (6), on déduit que

$$\displaystyle c_\Sigma(x(\xi,c))D^d(\xi)^{1/2}=w(d_\natural)^{-1}\sum_{\substack{I_1,I_2 \\ J_1,J_2}}w(d_+-d_{I_1})w(d_--d_{I_2}) B(I_1,I_2) \Delta_{\mu_+,\mu_-}(\xi(I_1),\xi(I_2),c)$$

\noindent où la somme porte toujours sur les mêmes quadruplets qu'en (3). Examinons la somme sur $(J_1,J_2)$ à $(I_1,I_2)$ fixé. La somme est vide sauf si $I_1\sqcup I_2=I$, $d_+-d_{I_1}$ et $d_--d_{I_2}$ sont des entiers positifs pairs. Si l'on est dans ce cas, alors $(J_1,J_2)$ est soumis aux seules conditions $J_1\sqcup J_2=J$, $2|J_1|=d_+-d_{I_1}$ et $2|J_2|=d_--d_{I_2}$. Il y a $\binom{|J|}{(d_+-d_{I_1})/2}$ tels couples. On vérifie facilement l'égalité

$$\displaystyle \binom{|J|}{(d_+-d_{I_1})/2}w(d_+-d_{I_1})w(d_--d_{I_2})=w(d-\underline{d})=w(d_\natural)$$

\noindent On arrive donc à l'égalité

$$\mbox{(7)}\;\;\; \displaystyle c_\Sigma(x(\xi,c))D^d(\xi)^{1/2}=\sum_{I_1,I_2}B(I_1,I_2) \Delta_{\mu_+,\mu_-}(\xi(I_1),\xi(I_2),c)$$

\noindent où la somme porte sur les couples $(I_1,I_2)$ vérifiant $I_1\sqcup I_2=I$, $d_+-d_{I_1}$ et $d_--d_{I_2}$ positifs pairs. \\

\noindent\textbullet Supposons maintenant que $\underline{d}$ est impair. Il existe alors une famille $(a_j)_{j\in J}$ d'éléments non nuls de $E$ et $j_0\in J$ tels que $\zeta_\lambda=(J,(F_{\pm j})_{j\in J}, (F_j)_{j\in J},(y_j(\lambda))_{j\in J})$ pour tout $\lambda\in F^\times$ assez petit, où:

\begin{itemize}
\item $F_{\pm j_0}=F$, $F_{j_0}=E$, $Tr_{E/F}(a_{j_0})=0$ et $y_{j_0}(\lambda)=e^{\lambda a_{j_0}}$;
\item $F_{\pm j}=E$, $F_{j}=E\times E$ et $y_{j}(\lambda)=(e^{\lambda a_{j}},e^{-\lambda\overline{a}_j})$ pour tout $j\in J\backslash\{j_0\}$.
\end{itemize}

\noindent On a alors $C(\zeta_\lambda)=F^\times/N(E^\times)$, donc $\zeta_\lambda$ ne détermine pas entièrement la classe de conjugaison de $e^{\lambda Y_\natural}$. Celle-ci est déterminée par la donnée supplémentaire d'un élément $c_0\in F^\times/N(E^\times)$. L'hypothèse 4.1(1) sur $(V,h)$ et $(V',h')$ entraîne $c_0=\nu_0$. La classe de conjugaison de $x(\xi,c)e^{\lambda Y_\natural}$ est alors paramétrée par $\xi\sqcup\zeta_\lambda$ et $(c,c_0)\in C(\xi\sqcup \zeta_\lambda)=C(\xi)\times F^\times/N(E^\times)$. Puisque $\Theta_\Sigma$ est un transfert de $\Theta_{\Sigma_+}\times \Theta_{\Sigma_-}$, pour tout $\lambda$ assez petit, on a

\[\begin{aligned}
\mbox{(8)}\;\;\; & \Theta_{\Sigma}(x(\xi,c)e^{\lambda Y_\natural})D^G(x(\xi,c)e^{\lambda Y_\natural})^{1/2}= \\
 & \sum_{\substack{I_1, I_2 \\ J_1, J_2}} \Theta_{\Sigma_+}(\xi(I_1)\sqcup \zeta_\lambda(J_1))D^{G_+}(\xi(I_1)\sqcup \zeta_\lambda(J_1))^{1/2} \Theta_{\Sigma_-}(\xi(I_2)\sqcup \zeta_\lambda(J_2)) \\
 & D^{G_-}(\xi(I_2)\sqcup \zeta_\lambda(J_2))^{1/2} \Delta_{\mu_+,\mu_-}\big(\xi(I_1)\sqcup \zeta_\lambda(J_1),\xi(I_2)\sqcup \zeta_\lambda(J_2),(c,c_0)\big)
\end{aligned}\]

\noindent où la somme porte sur les quadruplets $(I_1,I_2,J_1,J_2)$ vérifiant $I_1\sqcup I_2=I$, $J_1\sqcup J_2=J$, $d_{I_1}+d_{J_1}=d_+$ et $d_{I_2}+d_{J_2}=d_-$. Fixons un tel quadruplet. On a toujours (4) et (5). Il y a essentiellement deux cas possibles: $j_0\in J_1$ ou $j_0\in J_2$. Encore d'après les points (A)(i) et (A)(ii) du lemme 2.3.1, pour $\lambda$ assez proche de $0$, on a

\[\begin{aligned}
\displaystyle & \Delta_{\mu_+,\mu_-}\big(\xi(I_1)\sqcup \zeta_\lambda(J_1),\xi(I_2)\sqcup \zeta_\lambda(J_2),(c, c_0)\big)= \\
 & \left\{
    \begin{array}{ll}
        \Delta_{\mu_+,\mu_-}\big(\xi(I_1)\sqcup \zeta_\lambda(j_0),\xi(I_2),(c, c_0)\big) & \mbox{si } j_0\in J_1 \\
        \Delta_{\mu_+,\mu_-}\big(\xi(I_1),\xi(I_2)\sqcup \zeta_\lambda(j_0),(c, c_0)\big) & \mbox{si } j_0\in J_2
    \end{array}
\right.
\end{aligned}\]

\noindent où on a posé $\zeta_\lambda(j_0)=\zeta_\lambda(\{j_0\})$. D'après les points (B)(i) et (C) du lemme 2.3.1, on a donc

\[\begin{aligned}
\mbox{(9)}\;\;\; \displaystyle & \lim\limits_{\lambda\to 0} \Delta_{\mu_+,\mu_-}\big(\xi(I_1)\sqcup \zeta_\lambda(J_1),\xi(I_2)\sqcup \zeta_\lambda(J_2),(c, c_0)\big)= \\
 & \left\{
    \begin{array}{ll}
        \Delta_{\mu_+,\mu_-}(\xi(I_1),\xi(I_2),c) sgn_{E/F}\big(\delta^{-d_{I_2}} \frac{P_{I_2}(1)}{P_{I_2}(-1)}\big) & \mbox{si } j_0\in J_1 \\
        \Delta_{\mu_+,\mu_-}(\xi(I_1),\xi(I_2),c) sgn_{E/F}\big(\delta^{-d_{I_1}} \frac{P_{I_1}(1)}{P_{I_1}(-1)}\big) & \mbox{si } j_0\in J_2
    \end{array}
\right.
\end{aligned}\]

\noindent où on a posé $P_{I'}=P_{\xi(I')}$ pour tout $I'\subset I$ (dans le deuxième cas, où on applique 2.3.1(C), on a $c_b=\nu=\nu_0$ et $\underline{d}$ impair). De (2), (4), (5), (8) et (9), on déduit que

\[\begin{aligned}
\mbox{(10)}\;\;\; \displaystyle c_\Sigma(x(\xi,c))D^d(\xi)^{1/2}= & w(d_\natural)^{-1}\bigg(\sum_{\substack{I_1,I_2 \\ J_1, J_2}} w(d_+-d_{I_1})w(d_--d_{I_2}) B(I_1,I_2) \Delta_{\mu_+,\mu_-}(\xi(I_1),\xi(I_2),c) \\
 & sgn_{E/F}\big(\delta^{-d_{I_2}} \frac{P_{I_2}(1)}{P_{I_2}(-1)}\big) + \sum_{\substack{I_1,I_2 \\ J_1, J_2}} w(d_+-d_{I_1})w(d_--d_{I_2}) \\
 & B(I_1,I_2) \Delta_{\mu_+,\mu_-}(\xi(I_1),\xi(I_2),c) sgn_{E/F}\big(\delta^{-d_{I_1}} \frac{P_{I_1}(1)}{P_{I_1}(-1)}\big) \bigg)
\end{aligned}\]

\noindent la première somme (resp. la deuxième) portant sur les mêmes quadruplets que (8) avec la condition supplémentaire $j_0\in J_1$ (resp. $j_0\in J_2$). Dans la première somme, à $(I_1,I_2)$ fixé, la somme sur $(J_1,J_2)$ est vide sauf si $I_1\sqcup I_2=I$, $d_+-d_{I_1}$ est positif impair et $d_--d_{I_2}$ est positif pair. Dans ce cas, elle contient $\binom{|J|-1}{(d_--d_{I_2})/2}$ termes . On vérifie que

$$\displaystyle \binom{|J|-1}{(d_--d_{I_2})/2}w(d_+-d_{I_1})w(d_--d_{I_2})=w(d-\underline{d})=w(d_\natural)$$

\noindent Dans la deuxième somme, à $(I_1,I_2)$ fixé, la somme sur $(J_1,J_2)$ est vide sauf si $I_1\sqcup I_2=I$, $d_+-d_{I_1}$ est positif pair et $d_--d_{I_2}$ est positif impair. Dans ce cas, elle contient $\binom{|J|-1}{(d_+-d_{I_1})/2}$ termes. On a de même

$$\displaystyle \binom{|J|-1}{(d_+-d_{I_1})/2}w(d_+-d_{I_1})w(d_--d_{I_2})=w(d_\natural)$$

\noindent Par conséquent, (10) devient

\[\begin{aligned}
\mbox{(11)}\;\;\; \displaystyle c_\Sigma(x(\xi,c))D^d(\xi)^{1/2} & =\sum_{I_1,I_2}B(I_1,I_2) \Delta_{\mu_+,\mu_-}(\xi(I_1),\xi(I_2),c)sgn_{E/F}\big(\delta^{-d_{I_2}} \frac{P_{I_2}(1)}{P_{I_2}(-1)}\big) \\
 & + \sum_{I_1,I_2}B(I_1,I_2) \Delta_{\mu_+,\mu_-}(\xi(I_1),\xi(I_2),c)sgn_{E/F}\big(\delta^{-d_{I_1}} \frac{P_{I_1}(1)}{P_{I_1}(-1)}\big)
\end{aligned}\]

\noindent la première (resp. la deuxième) somme portant sur les couples $(I_1,I_2)$ vérifiant $I_1\sqcup I_2=I$, $d_+-d_{I_1}$ positif impair (resp. pair) et $d_--d_{I_2}$ positif pair (resp. impair). \\

On montre de la même façon des formules analogues à (7) et (11) pour $c_{\Sigma'}(x(\xi,c))D^{d'}(\xi)^{1/2}$. Plus précisément, pour tous $I_1,I_2\subset I$ vérifiant $d'_+\geqslant d_{I_1}$ et $d'_-\geqslant d_{I_2}$, posons

$$B'(I_1,I_2)=c_{\Sigma'_+}(\xi(I_1))D^{d'_+}(\xi(I_1))^{1/2} c_{\Sigma'_-}(\xi(I_2)) D^{d'_-}(\xi(I_2))^{1/2}$$

\noindent On a alors \\

\noindent\textbullet Si $\underline{d}$ est impair,

$$\mbox{(12)}\;\;\; \displaystyle c_{\Sigma'}(x(\xi,c))D^{d'}(\xi)^{1/2}=\sum_{I_1,I_2} B'(I_1,I_2) \Delta_{\mu'_+,\mu'_-}(\xi(I_1),\xi(I_2),c)$$

\noindent où la somme porte sur les couples $(I_1,I_2)$ vérifiant $I_1\sqcup I_2=I$, $d'_+-d_{I_1}$ et $d'_--d_{I_2}$ positifs et pairs. \\

\noindent\textbullet Si $\underline{d}$ est pair,

\[\begin{aligned}
\mbox{(13)}\;\;\; \displaystyle c_{\Sigma'}(x(\xi,c))D^{d'}(\xi)^{1/2} & =\sum_{I_1,I_2}B'(I_1,I_2) \Delta_{\mu'_+,\mu'_-}(\xi(I_1),\xi(I_2),c)sgn_{E/F}\big(\delta^{-d_{I_2}} \frac{P_{I_2}(1)}{P_{I_2}(-1)}\big) \\
 & + \sum_{I_1,I_2}B'(I_1,I_2) \Delta_{\mu'_+,\mu'_-}(\xi(I_1),\xi(I_2),c)sgn_{E/F}\big(\delta^{-d_{I_1}} \frac{P_{I_1}(1)}{P_{I_1}(-1)}\big)
\end{aligned}\]

\noindent la première (resp. la deuxième) somme portant sur les couples $(I_1,I_2)$ vérifiant $I_1\sqcup I_2=I$, $d'_+-d_{I_1}$ positif impair (resp. pair) et $d'_--d_{I_2}$ positif pair (resp. impair). \\

Les formules (7), (11), (12) et (13) permettent d'exprimer $f_1(\xi)$ en fonction de $\Theta_{\Sigma_+}\times \Theta_{\Sigma_-}$ et $\Theta_{\Sigma'_+}\times \Theta_{\Sigma'_-}$. Supposons à nouveau $\underline{d}$ pair. De (1), (7) et (13), on déduit que

\[\begin{aligned}
\mbox{(14)}\;\;\; \displaystyle f_1(\xi) & =\sum_{\substack{I_1,I_2 \\ I'_1,I'_2}} B(I_1,I_2)B'(I'_1,I'_2) sgn_{E/F}\big(\delta^{-d_{I'_2}} \frac{P_{I'_2}(1)}{P_{I'_2}(-1)}\big) A(I_1,I_2,I'_1,I'_2) \\
 & +\sum_{\substack{I_1,I_2 \\ I'_1,I'_2}} B(I_1,I_2)B'(I'_1,I'_2) sgn_{E/F}\big(\delta^{-d_{I'_1}} \frac{P_{I'_1}(1)}{P_{I'_1}(-1)}\big) A(I_1,I_2,I'_1,I'_2)
\end{aligned}\]

\noindent la première (resp. la deuxième) somme portant sur les quadruplets $(I_1,I_2,I'_1,I'_2)$ vérifiant $I_1\sqcup I_2=I'_1\sqcup I'_2=I$, $d_+-d_{I_1}$ et $d_--d_{I_2}$ positifs pairs, $d'_+-d_{I'_1}$ positif impair (resp. pair) et $d'_--d_{I'_2}$ positif pair (resp. impair), et où on a posé

$$\displaystyle A(I_1,I_2,I'_1,I'_2)=\sum_{c\in C(\xi)^\epsilon} \Delta_{\mu_+,\mu_-}(\xi(I_1),\xi(I_2),c)\Delta_{\mu'_+,\mu'_-}(\xi(I'_1),\xi(I'_2),c)$$

\noindent Examinons ce dernier terme. C'est à une constante près la somme de 

$$\displaystyle \prod_{i\in I_2} sgn_{F_i/F_{\pm i}}(c_i) \prod_{i\in I'_2} sgn_{F_i/F_{\pm i}}(c_i)$$

\noindent pour $c\in C(\xi)^\epsilon$. Cette somme est donc nulle sauf si $(I_1,I_2)=(I'_1,I'_2)$ ou $(I_1,I_2)=(I'_2,I'_1)$. Dans le cas où $(I_1,I_2)=(I'_1,I'_2)$, on a

$$\Delta_{\mu_+,\mu_-}(\xi(I_1),\xi(I_2),c)\Delta_{\mu'_+,\mu'_-}(\xi(I_1),\xi(I_2),c)=\mu_+\mu'_+(P_{I_1}(-1)) \mu_-\mu'_-(P_{I_2}(-1))$$

\noindent pour tout $c\in C(\xi)^\epsilon$. D'où

$$A(I_1,I_2,I_1,I_2)=|C(\xi)^\epsilon|\mu_+\mu'_+(P_{I_1}(-1)) \mu_-\mu'_-(P_{I_2}(-1))$$

\noindent Considérons maintenant le cas $(I_1,I_2)=(I'_2,I'_1)$. La formule (6) relie $\Delta_{\mu_+,\mu_-}(\xi(I_1),\xi(I_2),c)$ au facteur de transfert relatif au groupe endoscopique $G_+\times G_-$ de $G$. Comme on le sait, permuter $G_+$ et $G_-$, ne change pas le facteur de transfert si $G$ est quasi-déployé, multiplie le facteur de transfert par $-1$ si $G$ n'est pas quasi-déployé. Par conséquent

$$\Delta_{\mu_+,\mu_-}(\xi(I_1),\xi(I_2),c)=\mu(G)\Delta_{\mu_+,\mu_-}(\xi(I_2),\xi(I_1),c)$$

\noindent On en déduit que

\[\begin{aligned}
\displaystyle \Delta_{\mu_+,\mu_-}(\xi(I_1),\xi(I_2),c) & \Delta_{\mu'_+,\mu'_-}(\xi(I_2),\xi(I_1),c) \\
 & =\mu(G)\mu_+\mu'_-(P_{I_1}(-1)) \mu_-\mu'_+(P_{I_2}(-1))
\end{aligned}\]

\noindent Puis

$$A(I_1,I_2,I_2,I_1)=\mu(G)|C(\xi)^\epsilon|\mu_+\mu'_-(P_{I_1}(-1)) \mu_-\mu'_+(P_{I_2}(-1))$$

\noindent Revenons à (14). Supposons $\xi\neq \emptyset$ et $d'_+$ de même parité que $d_+$ et $d_-$. Alors $|C(\xi)^\epsilon|=|C(\xi)|/2$ et on ne peut pas avoir $(I_1,I_2)=(I'_1,I'_2)$ ni $(I_1,I_2)=(I'_2,I'_1)$ dans la première somme. Par conséquent seule la deuxième somme contribue et, d'après ce qui précède et le fait que les cas $(I_1,I_2)=(I'_1,I'_2)$ et $(I_1,I_2)=(I'_2,I'_1)$ sont exclusifs (car $\xi\neq\emptyset$), celle-ci vaut $|C(\xi)|/2$ multiplié par

\[\begin{aligned}
\displaystyle & \sum_{I_1,I_2}B(I_1,I_2)B'(I_1,I_2)sgn_{E/F}\big(\delta^{-d_{I_1}} \frac{P_{I_1}(1)}{P_{I_1}(-1)}\big)\mu_+\mu'_+(P_{I_1}(-1)) \\
 & \mu_-\mu'_-(P_{I_2}(-1))+\mu(G)\sum_{I_1,I_2}B(I_1,I_2)B'(I_2,I_1)sgn_{E/F}\big(\delta^{-d_{I'_2}} \frac{P_{I_2}(1)}{P_{I_2}(-1)}\big) \\
 & \mu_+\mu'_-(P_{I_1}(-1)) \mu_-\mu'_+(P_{I_2}(-1))
\end{aligned}\]

\noindent la première (resp. deuxième) somme portant sur les couples $(I_1,I_2)$ vérifiant $I_1\sqcup I_2=I$, $d_{I_1}\in\mathcal{D}(d_+,d'_+)$ (resp. $d_{I_1}\in\mathcal{D}(d_+,d'_-)$) et $d_{I_2}\in\mathcal{D}(d_-,d'_-)$ (resp. $d_{I_2}\in\mathcal{D}(d_-,d'_+)$). Dans le cas que l'on considère $(\mu_+\mu'_+)_{|F^\times}=sgn_{E/F}$ et $(\mu_-\mu'_+)_{|F^\times}=sgn_{E/F}$, donc

$$sgn_{E/F}\big(\delta^{-d_{I_1}} \frac{P_{I_1}(1)}{P_{I_1}(-1)}\big)\mu_+\mu'_+(P_{I_1}(-1))=\mu_+\mu'_+(\delta^{-d_{I_1}}P_{I_1}(1))$$

\noindent et

$$sgn_{E/F}\big(\delta^{-d_{I_2}} \frac{P_{I_2}(1)}{P_{I_2}(-1)}\big)\mu_-\mu'_+(P_{I_2}(-1))=\mu_+\mu'_+(\delta^{-d_{I_2}}P_{I_2}(1))$$

\noindent Toujours dans le cas que l'on considère, $\mu_+\mu'_-$ et $\mu_-\mu'_-$ sont triviaux sur $F^\times$. On vérifie que $\delta^{-d_{I'}}P_{I'}(1)P_{I'}(-1)^{-1}$ est élément de $F^\times$ pour tout $I'\subset I$. Par conséquent

$$\mu_+\mu'_-(P_{I_1}(-1))=\mu_+\mu'_-(\delta^{-d_{I_1}}P_{I_1}(1))$$

\noindent et

$$\mu_-\mu'_-(P_{I_2}(-1))=\mu_-\mu'_-(\delta^{-d_{I_2}}P_{I_2}(1))$$

\noindent On peut donc réécrire (14) sous la forme

\[\begin{aligned}
\displaystyle & f_1(\xi)=\frac{|C(\xi)|}{2}\bigg(\sum_{\substack{I_1\sqcup I_2=I \\ d_{I_1}\in \mathcal{D}(d_+,d'_+) \\ d_{I_2}\in\mathcal{D}(d_-,d'_-)}} B(I_1,I_2)B'(I_1,I_2)\mu_+\mu'_+(\delta^{-d_{I_1}}P_{I_1}(1)) \mu_-\mu'_-(\delta^{-d_{I_2}}P_{I_2}(1)) \\
 & +\mu(G)\sum_{\substack{I_1\sqcup I_2=I \\ d_{I_1}\in \mathcal{D}(d_+,d'_-) \\ d_{I_2}\in\mathcal{D}(d_-,d'_+)}}B(I_1,I_2)B'(I_2,I_1)\mu_+\mu'_-(\delta^{-d_{I_1}}P_{I_1}(1))\mu_-\mu'_+(\delta^{-d_{I_2}}P_{I_2}(1))\bigg)
\end{aligned}\]

\noindent On vérifie facilement que le membre de droite vaut $f_2(\xi)$. Supposons maintenant $\xi=\emptyset$. Alors le membre de droite de (14) ne comprend qu'un seul terme et on trouve

$$f_1(\xi)=|C(\xi)^\epsilon|c_{\Sigma_+}(1)c_{\Sigma_-}(1)c_{\Sigma'_+}(1)c_{\Sigma'_-}(1)$$

\noindent D'un autre côté, on a

$$\displaystyle f_2(\xi)=\frac{1+\mu(G)}{2}c_{\Sigma_+}(1)c_{\Sigma_-}(1)c_{\Sigma'_+}(1)c_{\Sigma'_-}(1)$$

\noindent Dans tous les cas, on a $|C(\xi)^\epsilon|=\frac{1+\mu(G)}{2}$ ce qui établit l'égalité voulue. Le cas où $d'_+$ est de parité différente de $d_+$ et $d_-$ et le cas où $\underline{d}$ est impair se traitent de la même manière $\blacksquare$

\subsection{Transfert et facteur epsilon}

\noindent Donnons nous

\vspace{3mm}

\begin{itemize}
\renewcommand{\labelitemi}{$\bullet$}

\item $U$ et $U'$ deux $E$-espaces vectoriels de dimensions respectives $d$, $d'$ et dont les groupes tordus associés (comme en 2.2) sont notés $(M,\widetilde{M})$ et $(M',\widetilde{M}')$ respectivement;

\item $(V_+,h_+)$, $(V_-,h_-)$, $(V'_+,h'_+)$ et $(V'_-,h'_-)$ quatre espaces hermitiens de dimensions respectives $d_+$, $d_-$, $d'_+$, $d'_-$ et de groupes unitaires respectifs $G_+$, $G_-$, $G'_+$ et $G'_-$;

\item $\mu_+$, $\mu_-$, $\mu'_+$ et $\mu'_-$ des caractères continus de $E^\times$.
\end{itemize}

\vspace{3mm}

\noindent On suppose que

\vspace{3mm}

\begin{itemize}
\renewcommand{\labelitemi}{$\bullet$}
\item $d$ est pair et $d'$ est impair;
\item $(V_+,h_+)$, $(V_-,h_-)$, $(V'_+,h'_+)$ et $(V'_-,h'_-)$ vérifient l'hypothèse \textbf{(QD)};
\item $d_++d_-=d$ et $d'_++d'_-=d'$;
\item $\mu_{+|F^\times}=sgn_{E/F}^{d_-}$, $\mu_{-|F^\times}=sgn_{E/F}^{d_++1}$, $\mu'_{+|F^\times}=sgn_{E/F}^{d'_-}$ et $\mu'_{-|F^\times}=sgn_{E/F}^{d'_++1}$
\end{itemize}

\vspace{3mm}

\noindent D'après 3.2, on peut considérer $G_+\times G_-$ comme un groupe endoscopique tordu de $\widetilde{M}$ et $G'_+\times G'_-$ comme un groupe endoscopique tordu de $\widetilde{M}'$ (les caractères $\mu_+$, $\mu_-$, $\mu'_+$ et $\mu'_-$ permettant de fixer les données endoscopiques). On normalise les facteurs de transfert comme en 3.2. Considérons

\vspace{3mm}

\begin{itemize}
\renewcommand{\labelitemi}{$\bullet$}

\item $\widetilde{\Pi}$ et $\widetilde{\Pi}'$ des représentations virtuelles tempérées de $\widetilde{M}$ et $\widetilde{M}'$ respectivement;

\item $\Sigma_+$, $\Sigma_-$, $\Sigma'_+$ et $\Sigma'_-$ des représentations virtuelles tempérées de $G_+$, $G_-$, $G'_+$ et $G'_-$ respectivement.
\end{itemize}

\vspace{3mm}

\noindent On suppose

\vspace{3mm}

\begin{itemize}
\renewcommand{\labelitemi}{$\bullet$}
\item Les caractères $\Theta_{\Sigma_+}$, $\Theta_{\Sigma_-}$, $\Theta_{\Sigma'_+}$ et $\Theta_{\Sigma'_-}$ sont stables;
\item $\Theta_{\widetilde{\Pi}}$ est un transfert de $\Theta_{\Sigma_+}\times \Theta_{\Sigma_-}$ et $\Theta_{\widetilde{\Pi}'}$ est un transfert de $\Theta_{\Sigma'_+}\times \Theta_{\Sigma'_-}$.
\end{itemize}

\vspace{3mm}

\noindent En prolongeant par bilinéarité l'application $(\tilde{\pi},\tilde{\pi}')\mapsto \epsilon_{\nu_1}(\tilde{\pi},\tilde{\pi}')$, définie en 4.3, on définit $\epsilon_{\nu_1}(\widetilde{\Pi},\widetilde{\Pi}')$.

\begin{prop}
Sous ces hypothèses, on a l'égalité

$$\epsilon_{\nu_1}(\widetilde{\Pi},\widetilde{\Pi}')=sgn_{E/F}(-\nu_1)^{d_-+d'_-}S_{\mu_+\mu'_+}(\Sigma_+,\Sigma'_+)S_{\mu_-\mu'_-}(\Sigma_-,\Sigma'_-)$$
\end{prop}

\noindent\ul{Preuve}: Exprimons $\epsilon_{\nu_1}(\widetilde{\Pi},\widetilde{\Pi}')$ grâce à la formule 4.4(1). Puisque $\mathcal{X}(d,d')$ est un revêtement de $\Xi^*(d,d')$, on peut réécrire cette formule sous la forme

$$\displaystyle \epsilon_{\nu_1}(\widetilde{\Pi},\widetilde{\Pi}')=\lim\limits_{s\to 0^+}\int_{\Xi^*(d,d')} f_1(\xi) \Delta(\xi)^{s-1/2} d\xi$$

\noindent Le membre de droite est défini comme la limite d'une intégrale sur $\Xi^*(d_+,d'_+)\times \Xi^*(d_-,d'_-)$. L'application $(\xi_1,\xi_2)\mapsto \xi_1\sqcup \xi_2$ définit un revêtement qui préserve localement les mesures d'un ouvert de complémentaire de mesure nulle de cet espace sur $\Xi^*(d,d')$. Puisque $\Delta(\xi_1\sqcup\xi_2)=\Delta(\xi_1)\Delta(\xi_2)$, on peut donc réécrire le membre de droite sous la forme

$$\displaystyle\lim\limits_{s\to 0^+}\int_{\Xi^*(d,d')} f_2(\xi) \Delta(\xi)^{s-1/2} d\xi$$

\noindent Pour établir le lemme, il suffit donc de prouver l'égalité $f_1(\xi)=f_2(\xi)$ pour tout $\xi\in \Xi^*(d,d')$. Soit $\xi=(I,(F_{\pm i})_{i\in I},(F_i)_{i\in I},(y_i)_{i\in I})\in \Xi^*(d,d')$. En tenant compte du jacobien de l'application $\mathcal{X}(d,d')\to \Xi^*(d,d')$ (que l'on a calculé en 2.2), on a

$$\mbox{(1)}\;\;\; \displaystyle f_1(\xi)=|2|_F^{-(d+d')/2}\sum_{\gamma\in\Gamma(\xi)} c_{\widetilde{\Pi}}(\tilde{x}(\xi,\gamma))D^{\tilde{M}}(\tilde{x}(\xi,\gamma))^{1/2} c_{\widetilde{\Pi}'}(\tilde{x}(\xi,\gamma))D^{\tilde{M}'}(\tilde{x}(\xi,\gamma))^{1/2}$$

\noindent Posons $\underline{d}=d_\xi$. Introduisons comme en 4.4 des décompositions $U=U_{\underline{d}}\oplus W_{\underline{d}}$, $U'=U_{\underline{d}}\oplus W'_{\underline{d}}$, des formes hermitiennes $h_{\underline{d}}$ sur $W_{\underline{d}}$ et $h'_{\underline{d}}$ sur $W'_{\underline{d}}$ et le groupe tordu $(M_{\underline{d}},\widetilde{M}_{\underline{d}})$ associé à $U_{\underline{d}}$. Alors en 4.4, on a défini un plongement de $\widetilde{M}_{\underline{d}}$ dans $\widetilde{M}$ et $\widetilde{M}'$. Pour tout $\gamma\in \Gamma(\xi)$, $\tilde{x}(\xi,\gamma)$ représente une classe de conjugaison dans $\widetilde{M}_{\underline{d}}$ et les coefficients $c_{\widetilde{\Pi}}(\tilde{x}(\xi,\gamma))$, $D^{\tilde{M}}(\tilde{x}(\xi,\gamma))$, $c_{\widetilde{\Pi}'}(\tilde{x}(\xi,\gamma))$, $D^{\tilde{M}'}(\tilde{x}(\xi,\gamma))$ sont définis au moyen des plongements. \\

Soit $\gamma\in \Gamma(\xi)$. Notons $G_{\natural}$ le groupe unitaire de $(W_{\underline{d}},h_{\underline{d}})$. Ce groupe est quasi-déployé. Fixons un sous-tore maximal $T_\natural$ d'un sous-groupe de Borel de $G_\natural$ tous deux définis sur $F$. Soit $Y_\natural\in \mathfrak{t}_\natural(F)\cap \mathfrak{g}_{\natural,reg}(F)$ qui n'admet pas $0$ comme valeur propre dans $W_{\underline{d}}$. D'après 4.5(1), on a

$$\mbox{(2)}\;\;\; \displaystyle c_{\widetilde{\Pi}}(\tilde{x}(\xi,\gamma))D^{\tilde{M}}(\tilde{x}(\xi,\gamma))^{1/2}=w(d_\natural)^{-1}\lim\limits_{\substack{\lambda\in F^\times \\ \lambda\to 0}} \Theta_{\widetilde{\Pi}}(\tilde{x}(\xi,\gamma)e^{\lambda Y_\natural})D^{\tilde{M}}(\tilde{x}(\xi,\gamma)e^{\lambda Y_\natural})^{1/2}$$

\noindent où $d_\natural=d-\underline{d}$ est la dimension de $W_{\underline{d}}$. Pour tout $\lambda\in F^\times$ assez petit, la classe de conjugaison stable de $e^{\lambda Y_\natural}$ admet pour paramètre $\zeta_\lambda=(J,(F_{\pm j})_{j\in J},(F_j)_{j\in J},(y_j(\lambda))_{j\in J})$ (les trois premières données du quadruplet ne dépendent pas de $\lambda$). Alors la classe de conjugaison stable de $\tilde{x}(\xi,\gamma)e^{\lambda Y_\natural}$ est elle paramétrée par $\xi\sqcup \zeta_\lambda^2$ où $\zeta_\lambda^2=(J,(F_{\pm j})_{j\in J},(F_j)_{j\in J},(y_j(\lambda)^2)_{j\in J})$. \\

\noindent\textbullet Supposons dans un premier temps que $\underline{d}$ est pair. Alors $F_{\pm j}=E$, $F_j=E\times E$ pour tout $j\in J$ et il existe une famille $(a_j)_{j\in J}$ d'éléments distincts et de traces non nulles de $E$ telle que $y_j(\lambda)=(e^{\lambda a_j},e^{-\lambda\overline{a}_j})$, pour tout $\lambda$ assez petit. On a alors $\Gamma(\xi\sqcup \zeta_\lambda^2)=\Gamma(\xi)$ et la classe de conjugaison de $\tilde{x}(\xi,\gamma)e^{\lambda Y_\natural}$ est paramétrée par $\gamma\in \Gamma(\xi\sqcup\zeta_\lambda^2)$. Pour tout $I'\subset I$, on pose $\xi(I')=(I',(F_{\pm i})_{i\in I'}, (F_i)_{i\in I'}, (y_i)_{i\in I'})$. On définit de même $\zeta_\lambda^2(J')$ pour $J'\subset J$. Puisque $\Theta_{\widetilde{\Pi}}$ est un transfert de $\Theta_{\Sigma_+}\times \Theta_{\Sigma_-}$, pour tout $\lambda$ assez petit, on a l'égalité

\[\begin{aligned}
\mbox{(3)}\;\;\; \displaystyle & \Theta_{\widetilde{\Pi}}(\tilde{x}(\xi,\gamma)e^{\lambda Y_\natural})\cdot D_0^{\tilde{M}}(\tilde{x}(\xi,\gamma)e^{\lambda Y_\natural})^{1/2}= \\
 & \sum_{\substack{I_1, I_2\\ J_1, J_2}} \Theta_{\Sigma_+}(\xi(I_1)\sqcup \zeta_\lambda^2(J_1))\cdot D^{d_+}(\xi(I_1)\sqcup \zeta_\lambda^2(J_1))^{1/2} \cdot \Theta_{\Sigma_-}(\xi(I_2)\sqcup \zeta_\lambda^2(J_2)) \\
 & \cdot D^{d_-}(\xi(I_2)\sqcup \zeta_\lambda^2(J_2))^{1/2} \cdot \Delta_{\mu_+,\mu_-}(\xi(I_1)\sqcup \zeta_\lambda^2(J_1),\xi(I_2)\sqcup \zeta_\lambda^2(J_2),\gamma)
\end{aligned}\]

\noindent où la somme porte sur les quadruplets $(I_1,I_2,J_1,J_2)$ vérifiant $I_1\sqcup I_2=I$, $J_1\sqcup J_2=J$, $d_{I_1}+2|J_1|=d_+$ et $d_{I_2}+2|J_2|=d_-$. Fixons un tel quadruplet. D'après 6.1(1), on a

$$\mbox{(4)}\;\;\; \displaystyle \lim\limits_{\lambda\to 0} \Theta_{\Sigma_+}(\xi(I_1)\sqcup \zeta_\lambda^2(J_1))D^{d_+}(\xi(I_1)\sqcup \zeta_\lambda^2(J_1))^{1/2}=w(d_+-d_{I_1}) c_{\Sigma_+}(\xi(I_1))D^{d_+}(\xi(I_1))^{1/2}$$

$$\mbox{(5)}\;\;\; \displaystyle \lim\limits_{\lambda\to 0} \Theta_{\Sigma_-}(\xi(I_2)\sqcup \zeta_\lambda^2(J_2))D^{d_-}(\xi(I_2)\sqcup \zeta_\lambda^2(J_2))^{1/2}=w(d_--d_{I_2}) c_{\Sigma_-}(\xi(I_2))D^{d_-}(\xi(I_2))^{1/2}$$

\noindent En appliquant plusieurs fois les points (A)(iii) et (A)(iv) du lemme 2.3.1, on montre que, pour $\lambda$ assez petit, on a

$$\mbox{(6)}\;\;\; \Delta_{\mu_+,\mu_-}(\xi(I_1)\sqcup \zeta_\lambda^2(J_1),\xi(I_2)\sqcup \zeta_\lambda^2(J_2),\gamma)=\Delta_{\mu_+,\mu_-}(\xi(I_1),\xi(I_2),\gamma)$$

\noindent On vérifie facilement que

$$\mbox{(7)}\;\;\; D^{\tilde{M}}(\tilde{x})=|2|_F^dD_0^{\tilde{M}}(\tilde{x})$$

\noindent pour tout $\tilde{x}\in \widetilde{M}_{reg}(F)$. Pour tous $I_1,I_2\subset I$ tels que $d_{I_1}\leqslant d_+$, $d_{I_2}\leqslant d_-$, posons

$$B(I_1,I_2)=c_{\Sigma_+}(\xi(I_1))D^{d_+}(\xi(I_1))^{1/2}c_{\Sigma_-}(\xi(I_2))D^{d_-}(\xi(I_2))^{1/2}$$

\noindent De (2), (3), (4), (5), (6) et (7), on déduit que

\[\begin{aligned}
\mbox{(8)}\;\;\; \displaystyle & c_{\widetilde{\Pi}}(\tilde{x}(\xi,\gamma))D^{\tilde{M}}(\tilde{x}(\xi,\gamma))^{1/2}= \\
 & w(d_\natural)^{-1}|2|_F^{d/2}\sum_{\substack{I_1, I_2\\ J_1, J_2}}w(d_+-d_{I_1})w(d_--d_{I_2})B(I_1,I_2)\Delta_{\mu_+,\mu_-}(\xi(I_1),\xi(I_2),\gamma)
\end{aligned}\]

\noindent où la somme porte sur les mêmes quadruplets qu'en (3). Fixons $(I_1,I_2)$ et examinons la somme sur $(J_1,J_2)$. Elle est vide sauf si $d_+-d_{I_1}$ et $d_--d_{I_2}$ sont positifs et pairs et si tel est le cas, elle comporte $\binom{|J|}{(d_+-d_{I_1})/2}$ termes. On vérifie que

$$\displaystyle\binom{|J|}{(d_+-d_{I_1})/2}w(d_+-d_{I_1})w(d_--d_{I_2})=w(d-\underline{d})=w(d_\natural)$$

\noindent Par conséquent, (8) devient

$$\mbox{(9)}\;\;\; \displaystyle c_{\widetilde{\Pi}}(\tilde{x}(\xi,\gamma))D^{\tilde{M}}(\tilde{x}(\xi,\gamma))^{1/2}=|2|_F^{d/2}\sum_{I_1,I_2}B(I_1,I_2)\Delta_{\mu_+,\mu_-}(\xi(I_1),\xi(I_2),\gamma)$$

\noindent où la somme porte sur les couples $(I_1,I_2)$ vérifiant $I_1\sqcup I_2=I$ et $d_+-d_{I_1}$, $d_--d_{I_2}$ pairs et positifs. \\

\noindent\textbullet Supposons maintenant que $\underline{d}$ soit impair. Alors il existe une famille $(a_j)_{j\in J}$ d'éléments de $E$ et $j_0\in J$ tels que, pour $\lambda$ assez petit,

\begin{itemize}
\item $F_{\pm j_0}=F$, $F_{j_0}=E$, $Tr_{E/F}(a_{j_0})=0$ et $y_{j_0}(\lambda)=e^{\lambda a_0}$;
\item $F_{\pm j}=E$, $F_j=E\times E$ et $y_j(\lambda)=(e^{\lambda a_j},e^{-\lambda \overline{a}_j})$ pour tout $j\in J\backslash \{j_0\}$.
\end{itemize}

\noindent On a alors $\Gamma(\xi\sqcup \zeta_\lambda^2)=\Gamma(\xi)\times \Gamma(e^{2\lambda a_0})/N(E^\times)$ et la classe de conjugaison de $\tilde{x}(\xi,\gamma)e^{\lambda Y_\natural}$ est paramétrée par $(\gamma,\gamma_0(\lambda))$ pour un certain $\gamma_0(\lambda)\in \Gamma(e^{2\lambda a_0})/N(E^\times)$. D'après la façon dont $\widetilde{M}_{\underline{d}}$ a été plongé dans $\widetilde{M}$, on a $\gamma_0(\lambda)=\nu_1 e^{\lambda a_0}$. Puisque $\Theta_{\widetilde{\Pi}}$ est un transfert de $\Theta_{\Sigma_+}\times \Theta_{\Sigma_-}$, pour tout $\lambda$ assez petit, on a l'égalité

\[\begin{aligned}
\mbox{(10)}\;\;\; \displaystyle & \Theta_{\widetilde{\Pi}}(\tilde{x}(\xi,\gamma)e^{\lambda Y_\natural})D_0^{\tilde{M}}(\tilde{x}(\xi,\gamma)e^{\lambda Y_\natural})^{1/2}= \\
 & \sum_{\substack{I_1, I_2\\ J_1, J_2}} \Theta_{\Sigma_+}(\xi(I_1)\sqcup \zeta_\lambda^2(J_1))D^{d_+}(\xi(I_1)\sqcup \zeta_\lambda^2(J_1))^{1/2}\Theta_{\Sigma_-}(\xi(I_2)\sqcup \zeta_\lambda^2(J_2)) \\
 & D^{d_-}(\xi(I_2)\sqcup \zeta_\lambda^2(J_2))^{1/2} \Delta_{\mu_+,\mu_-}(\xi(I_1)\sqcup \zeta_\lambda^2(J_1),\xi(I_2)\sqcup \zeta_\lambda^2(J_2),(\gamma,\gamma_0(\lambda)))
\end{aligned}\]

\noindent où la somme porte sur les quadruplets $(I_1,I_2,J_1,J_2)$ vérifiant $I_1\sqcup I_2=I$, $J_1\sqcup J_2=J$, $d_{I_1}+d_{J_1}=d_+$ et $d_{I_2}+d_{J_2}=d_-$. Fixons un tel quadruplet. On distingue deux cas suivant que $j_0\in J_1$ ou $j_0\in J_2$. Appliquant à nouveau les points (A)(iii) et (A)(iv) du lemme 2.3.1, on montre que pour $\lambda$ assez proche de $0$, on a

\[\begin{aligned}
\displaystyle & \Delta_{\mu_+,\mu_-}(\xi(I_1)\sqcup \zeta_\lambda^2(J_1),\xi(I_2)\sqcup \zeta_\lambda^2(J_2),(\gamma,\gamma_0(\lambda)))= \\
 & \left\{
    \begin{array}{ll}
        \Delta_{\mu_+,\mu_-}(\xi(I_1)\sqcup \zeta_\lambda^2(j_0),\xi(I_2),(\gamma,\gamma_0(\lambda))) & \mbox{si } j_0\in J_1 \\
        \Delta_{\mu_+,\mu_-}(\xi(I_1),\xi(I_2)\sqcup \zeta_\lambda^2(j_0),(\gamma,\gamma_0(\lambda))) & \mbox{si } j_0\in J_2
    \end{array}
\right.
\end{aligned}\]

\noindent D'après les points (B)(ii) et (D) du lemme 2.3.1, on en déduit que

\[\begin{aligned}
\displaystyle & \lim\limits_{\lambda\to 0}\Delta_{\mu_+,\mu_-}(\xi(I_1)\sqcup \zeta_\lambda^2(J_1),\xi(I_2)\sqcup \zeta_\lambda^2(J_2),(\gamma,\gamma_0(\lambda)))= \\
 & \left\{
    \begin{array}{ll}
        \Delta_{\mu_+,\mu_-}(\xi(I_1),\xi(I_2),\gamma) sgn_{E/F}(\delta^{-d_{I_2}}\frac{P_{I_2}(1)}{P_{I_2}(-1)}) & \mbox{si } j_0\in J_1 \\
        \Delta_{\mu_+,\mu_-}(\xi(I_1),\xi(I_2),\gamma) sgn_{E/F}(\delta^{-d_{I_1}}\frac{P_{I_1}(1)}{P_{I_1}(-1)})sgn_{E/F}(-\nu_1) & \mbox{si } j_0\in J_2
    \end{array}
\right.
\end{aligned}\]

\noindent où pour tout $I'\subset I$, on a posé $P_{I'}=P_{\xi(I')}$ (dans le deuxième cas, où on applique 2.3.1(D), on a $c_b=\nu_1$ et $\underline{d}$ impair). On a toujours (4) et (5). D'après (2),(7) et (10), on en déduit que

\[\begin{aligned}
\mbox{(11)}\;\;\; \displaystyle & c_{\widetilde{\Pi}}(\tilde{x}(\xi,\gamma))D^{\tilde{M}}(\tilde{x}(\xi,\gamma))^{1/2}=|2|_F^{d/2}w(d_\natural)^{-1}\bigg(\sum_{\substack{I_1,I_2\\ J_1, J_2}} w(d_+-d_{I_1})w(d_--d_{I_2}) B(I_1,I_2) \\
 & \Delta_{\mu_+,\mu_-}(\xi(I_1),\xi(I_2),\gamma) sgn_{E/F}(\delta^{-d_{I_2}}\frac{P_{I_2}(1)}{P_{I_2}(-1)})+\sum_{\substack{I_1, I_2\\ J_1, J_2}} w(d_+-d_{I_1})w(d_--d_{I_2}) B(I_1,I_2) \\
 & \Delta_{\mu_+,\mu_-}(\xi(I_1),\xi(I_2),\gamma)sgn_{E/F}(\delta^{-d_{I_1}}\frac{P_{I_1}(1)}{P_{I_1}(-1)})sgn_{E/F}(-\nu_1)\bigg)
\end{aligned}\]

\noindent la première somme (resp. la deuxième) portant sur les mêmes quadruplets que (10) avec la condition supplémentaire $j_0\in J_1$ (resp. $j_0\in J_2$). Dans la première somme, à $(I_1,I_2)$ fixé, la somme sur $(J_1,J_2)$ est vide sauf si $d_+-d_{I_1}$ est positif impair et $d_--d_{I_2}$ est positif pair. Si tel est le cas, elle contient $\binom{|J|-1}{(d_--d_{I_2})/2}$ termes. On vérifie que

$$\mbox{(12)}\;\;\; \displaystyle \binom{|J|-1}{(d_--d_{I_2})/2}w(d_+-d_{I_1})w(d_--d_{I_2})=w(d-\underline{d})=w(d_\natural)$$

\noindent De même, dans la deuxième somme, à $(I_1,I_2)$ fixé, la somme sur $(J_1,J_2)$ est vide sauf si $d_+-d_{I_1}$ est positif pair et $d_--d_{I_2}$ est positif impair. Si tel est le cas, elle contient $\binom{|J|-1}{(d_+-d_{I_1})/2}$ termes. On vérifie que

$$\mbox{(13)}\;\;\; \displaystyle \binom{|J|-1}{(d_+-d_{I_1})/2}w(d_+-d_{I_1})w(d_--d_{I_2})=w(d-\underline{d})=w(d_\natural)$$

\noindent De (11), (12) et (13), on déduit que

\[\begin{aligned}
\mbox{(14)}\;\;\; \displaystyle & c_{\widetilde{\Pi}}(\tilde{x}(\xi,\gamma))D^{\tilde{M}}(\tilde{x}(\xi,\gamma))^{1/2}= \\
 & |2|_F^{d/2}\bigg(\sum_{I_1,I_2} B(I_1,I_2)\Delta_{\mu_+,\mu_-}(\xi(I_1),\xi(I_2),\gamma) sgn_{E/F}(\delta^{-d_{I_2}}\frac{P_{I_2}(1)}{P_{I_2}(-1)}) \\
 & +\sum_{I_1,I_2} B(I_1,I_2)\Delta_{\mu_+,\mu_-}(\xi(I_1),\xi(I_2),\gamma) sgn_{E/F}(\delta^{-d_{I_1}}\frac{P_{I_1}(1)}{P_{I_1}(-1)})sgn_{E/F}(-\nu_1)\bigg)
\end{aligned}\]

\noindent où la première (resp. la deuxième) somme porte sur les couples $(I_1,I_2)$ vérifiant $I_1\sqcup I_2=I$, $d_+-d_{I_1}$ positif impair (resp. pair) et $d_--d_{I_2}$ positif pair (resp. impair). \\

\noindent On montre de la même façon des formules analogues à (9) et (14) pour $c_{\widetilde{\Pi}'}(\tilde{x}(\xi,\gamma)) D^{\tilde{M}'}(\tilde{x}(\xi,\gamma))^{1/2}$. Plus précisément, pour tous $I_1,I_2\subset I$ vérifiant $d'_+\geqslant d_{I_1}$ et $d'_-\geqslant d_{I_2}$, posons

$$B'(I_1,I_2)=c_{\Sigma'_+}(\xi(I_1))D^{d'_+}(\xi(I_1))^{1/2} c_{\Sigma'_-}(\xi(I_2)) D^{d'_-}(\xi(I_2))^{1/2}$$

\noindent On a alors \\

\noindent\textbullet Si $\underline{d}$ est impair, alors

$$\mbox{(15)}\;\;\; \displaystyle c_{\widetilde{\Pi}'}(\tilde{x}(\xi,\gamma)) D^{\tilde{M}'}(\tilde{x}(\xi,\gamma))^{1/2}=|2|_F^{d'/2}\sum_{I_1,I_2} B'(I_1,I_2) \Delta_{\mu'_+,\mu'_-}(\xi(I_1),\xi(I_2),\gamma)$$

\noindent où la somme porte sur les couples $(I_1,I_2)$ vérifiant $I_1\sqcup I_2=I$, $d'_+-d_{I_1}$ et $d'_--d_{I_2}$ positifs pairs. \\

\noindent\textbullet Si $\underline{d}$ est pair, alors

\[\begin{aligned}
\mbox{(16)}\;\;\; \displaystyle & c_{\widetilde{\Pi}'}(\tilde{x}(\xi,\gamma)) D^{\tilde{M}'}(\tilde{x}(\xi,\gamma))^{1/2}=|2|_F^{d'/2}\bigg(\sum_{I_1,I_2} B'(I_1,I_2)\Delta_{\mu'_+,\mu'_-}(\xi(I_1),\xi(I_2),\gamma) \\
 & sgn_{E/F}(\delta^{-d_{I_2}}\frac{P_{I_2}(1)}{P_{I_2}(-1)})+\sum_{I_1,I_2} B'(I_1,I_2)\Delta_{\mu'_+,\mu'_-}(\xi(I_1),\xi(I_2),\gamma) \\
 & sgn_{E/F}(\delta^{-d_{I_1}}\frac{P_{I_1}(1)}{P_{I_1}(-1)})sgn_{E/F}(-\nu_1)\bigg)
\end{aligned}\]

\noindent où la première (resp. la deuxième) somme porte sur les couples $(I_1,I_2)$ vérifiant $I_1\sqcup I_2=I$, $d'_+-d_{I_1}$ positif impair (resp. pair) et $d'_--d_{I_2}$ positif pair (resp. impair). \\

 On peut maintenant, grâce à (1), (9), (14), (15) et (16), exprimer $f_1(\xi)$ en fonction de $\Theta_{\Sigma_+}\times \Theta_{\Sigma_-}$ et $\Theta_{\Sigma'_+}\times \Theta_{\Sigma'_-}$. Supposons à nouveau $\underline{d}$ pair. On obtient

\[\begin{aligned}
\mbox{(17)}\;\;\; \displaystyle f_1(\xi)& =\sum_{\substack{I_1,I_2 \\ I'_1,I'_2}} B(I_1,I_2)B'(I'_1,I'_2) sgn_{E/F}\big(\delta^{-d_{I'_2}} \frac{P_{I'_2}(1)}{P_{I'_2}(-1)}\big) A(I_1,I_2,I'_1,I'_2) \\
 & +\sum_{\substack{I_1,I_2 \\ I'_1,I'_2}} B(I_1,I_2)B'(I'_1,I'_2) sgn_{E/F}\big(\delta^{-d_{I'_1}} \frac{P_{I'_1}(1)}{P_{I'_1}(-1)}\big)sgn_{E/F}(-\nu_1) A(I_1,I_2,I'_1,I'_2)
\end{aligned}\]

\noindent la première (resp. la deuxième) somme portant sur les quadruplets $(I_1,I_2,I'_1,I'_2)$ vérifiant $I_1\sqcup I_2=I'_1\sqcup I'_2=I$, $d_+-d_{I_1}$ et $d_--d_{I_2}$ positifs pairs, $d'_+-d_{I'_1}$ positif impair (resp. pair) et $d'_--d_{I'_2}$ positif pair (resp. impair), et où on a posé

$$\displaystyle A(I_1,I_2,I'_1,I'_2)=\sum_{\gamma\in \Gamma(\xi)} \Delta_{\mu_+,\mu_-}(\xi(I_1),\xi(I_2),\gamma)\Delta_{\mu'_+,\mu'_-}(\xi(I'_1),\xi(I'_2),\gamma)$$

\noindent Cette somme est à une constante près la somme des

$$\displaystyle \prod_{i\in I_2} sgn_{F_i/F_{\pm i}}(\gamma_i^{-1}(1+y_i))\prod_{i\in I'_2} sgn_{F_i/F_{\pm i}}(\gamma_i^{-1}(1+y_i))$$

\noindent pour $\gamma\in\Gamma(\xi)$. Comme $\Gamma(\xi)$ est un espace principal homogène sous $C(\xi)$, cette somme est nulle sauf si $(I_1,I_2)=(I'_1,I'_2)$. Si $(I_1,I_2)=(I'_1,I'_2)$, on a

$$\Delta_{\mu_+,\mu_-}(\xi(I_1),\xi(I_2),\gamma)\Delta_{\mu'_+,\mu'_-}(\xi(I'_1),\xi(I'_2),\gamma)=\mu_+\mu'_+(P_{I_1}(-1))\mu_-\mu'_-(P_{I_2}(-1))$$

\noindent pour tout $\gamma\in \Gamma(\xi)$. On a $|\Gamma(\xi)|=|C(\xi)|$, et par conséquent

$$\mbox{(18)}\;\;\; A(I_1,I_2,I'_1,I'_2)=|C(\xi)|\mu_+\mu'_+(P_{I_1}(-1))\mu_-\mu'_-(P_{I_2}(-1))$$

\noindent Supposons que $d'_+$ est de même parité que $d_+$ et $d_-$. On en déduit que $(\mu'_+\mu_+)_{|F^\times}=sgn_{E/F}$ et $(\mu_-\mu'_-)_{|F^\times}$ est trivial. Un quadruplet de la forme $(I_1,I_2,I_1,I_2)$ ne peut intervenir que dans la deuxième somme de (17) (donc la première somme est nulle). Pour un tel quadruplet, on vérifie que $\delta^{-d_{I_2}}P_{I_2}(1)P_{I_2}(-1)^{-1}\in F^\times$. Par conséquent

$$\mu_-\mu'_-(P_{I_2}(-1))=\mu_-\mu'_-(\delta^{-d_{I_2}}P_{I_2}(1))$$

\noindent et

$$sgn_{E/F}\big(\delta^{-d_{I_1}} \frac{P_{I_1}(1)}{P_{I_1}(-1)}\big)\mu_+\mu'_+(P_{I_1}(-1))=\mu_+\mu'_+(\delta^{-d_{I_1}}P_{I_1}(1))$$

\noindent Grâce à (18), on en déduit que

$$\mbox{(19)}\;\;\; \displaystyle f_1(\xi)= sgn_{E/F}(-\nu_1)\sum_{I_1,I_2} B(I_1,I_2)B'(I_1,I_2)\mu_+\mu'_+(\delta^{-d_{I_1}}P_{I_1}(1))\mu_-\mu'_-(\delta^{-d_{I_2}}P_{I_2}(1))$$

\noindent où la somme porte sur les couples $(I_1,I_2)$ vérifiant $I_1\sqcup I_2=I$, $d_+-d_{I_1}$, $d_--d_{I_2}$, $d'_+-d_{I_1}$ positifs pairs et $d'_--d_{I_2}$ positif impair. On vérifie que les quatres dernières conditions sont équivalentes à $d_{I_1}\in\mathcal{D}(d_+,d'_+)$ et $d_{I_2}\in \mathcal{D}(d_-,d'_-)$. On voit aussi que l'on peut remplacer le terme $sgn_{E/F}(-\nu_1)$ par $sgn_{E/F}(-\nu_1)^{d_-+d'_-}$. Le membre de droite de (19) s'identifie alors à $f_2(\xi)$. D'où l'égalité $f_1(\xi)=f_2(\xi)$. Les autres cas ($d'_-$ de même parité que $d_-$ et $d_+$ puis $\underline{d}$ impair) se traitent de la même façon $\blacksquare$

\section{$L$-paquets tempérés pour les groupes unitaires et endoscopie}

\subsection{Paramètres de Langlands et représentations conjuguées-duales du groupe de Weil-Deligne}

Notons $WD_E=W_E\times SL_2(\mathbb{C})$, où $W_E$ est le groupe de Weil, le groupe de Weil-Deligne de $E$. Soit $N\geqslant 1$ un entier. On appelle paramètre de Langlands de dimension $N$ tout homomorphisme continu

$$\varphi: WD_E\to GL_N(\mathbb{C})$$

\noindent vérifiant les deux conditions suivantes:

\begin{itemize}
\renewcommand{\labelitemi}{$\bullet$}
\item $\varphi$ est semi-simple;
\item la restriction de $\varphi$ à $SL_2(\mathbb{C})$ est algébrique;
\end{itemize}

On dira de plus que $\varphi$ est tempéré si l'image de $W_E$ par $\varphi$ est relativement compacte. Deux paramètres de Langlands $\varphi$ et $\varphi'$ sont dits conjugués, ce que l'on note $\varphi\simeq\varphi'$, s'ils ont même dimension $N$ et qu'ils sont conjugués par un élément de $GL_N(\mathbb{C})$. Notons $\mathbf{\Phi}_{temp}(GL_N)$ l'ensemble des paramètres de Langlands tempérés de dimension $N$ pris à conjugaison près. La correspondance locale de Langlands pour les groupes linéaires, due à Harris-Taylor et Henniart, associe à tout $\varphi\in\mathbf{\Phi}_{temp}(GL_N)$ une représentation irréductible tempérée $\pi(\varphi)$ de $GL_N(E)$. Notons $\mathbf{\Phi}_{temp,irr}(GL_N)$ le sous-ensemble des $\varphi\in\mathbf{\Phi}_{temp}(GL_N)$ qui sont irréductibles. Fixons $t\in W_F-W_E$, cet élément agit par conjugaison sur $W_E$, on prolonge cette action à $WD_E$ en laissant agir $t$ de façon triviale sur $SL_2(\mathbb{C})$. Pour tout $\varphi\in \mathbf{\Phi}_{temp}(GL_N)$, on définit un nouveau paramètre de Langlands $\varphi^\theta\in \mathbf{\Phi}_{temp}(GL_N)$ en posant $\varphi^\theta(\tau)={}^t\varphi(t\tau t^{-1})^{-1}$ pour tout $\tau\in WD_E$. L'élément $\varphi^\theta$ ainsi construit ne dépend pas (à conjugaison près) du choix de $t$. Un paramètre $\varphi\in \mathbf{\Phi}_{temp}(GL_N)$ sera dit conjugué-dual, si on a $\varphi\simeq\varphi^\theta$. C'est équivalent à l'existence d'une forme bilinéaire non dégénérée $B:\mathbb{C}^N\times \mathbb{C}^N\to \mathbb{C}$ vérifiant \\

(1) $B(\varphi(\tau)w,\varphi(t\tau t^{-1})w')=B(w,w')$, pour tous $w,w'\in \mathbb{C}^N$ et pour tout $\tau\in WD_E$. \\

\noindent Notons $\mathbf{\Phi}_{temp}^{\theta}(GL_N)$ le sous-ensemble des paramètres $\varphi\in \mathbf{\Phi}_{temp}(GL_N)$ qui sont conjugués-duals. Soit $\epsilon\in\{\pm\}$ un signe. Un élément $\varphi\in \mathbf{\Phi}^{\theta}_{temp}(GL_N)$ sera dit conjugué-dual de signe $\epsilon$ s'il existe une forme bilinéaire non dégénérée $B:\mathbb{C}^N\times \mathbb{C}^N\to \mathbb{C}$ qui vérifie (1) et la condition supplémentaire \\

(2) $B(w,\varphi(t^2)w')=\epsilon B(w',w)$ pour tous $w,w'\in \mathbb{C}^N$. \\

\noindent Cette définition ne dépend pas du choix de $t$. Une forme bilinéaire non dégénérée $B$ vérifiant les conditions (1) et (2) sera dite $\epsilon$-conjuguée-duale. Fixons une telle forme $B$ et notons $Aut(\varphi,B)$ le groupe des éléments $g\in GL_N(\mathbb{C})$ qui préservent $B$ et commutent à l'image de $\varphi$. A automorphisme intérieur près, ce groupe ne dépend pas de $t$ ni de $B$ et c'est un groupe algébrique réductif complexe (en général non connexe). On pose alors $\mathcal{S}_{\varphi}=Aut(\varphi,B)/Aut(\varphi,B)^0$ où $Aut(\varphi,B)^0$ désigne la composante neutre. Ce groupe est abélien et est donc bien défini à un unique isomorphisme près (indépendamment de $B$ et de $s$). La notation $\mathcal{S}_{\varphi}$ est cependant quelque peu imprécise: elle ne dénote pas le même objet suivant que l'on considère $\varphi$ comme un paramètre conjugué-dual de signe $+$ ou $-$ (certain paramètres peuvent être considérés des deux façons). Néanmoins, dans ce qui suit, le contexte devrait effacer toute confusion possible (ou en tout cas l'auteur espère qu'il en est ainsi). On note $\mathbf{\Phi}_{temp}^+(GL_N)$ (resp. $\mathbf{\Phi}_{temp}^-(GL_N)$) le sous-ensemble de $\mathbf{\Phi}^{\theta}_{temp}(GL_N)$ constitué des paramètres conjugués-duaux de signe $+$ (resp. de signe $-$). On a alors $\mathbf{\Phi}_{temp}^{\pm}(GL_N)\subset \mathbf{\Phi}_{temp}^{\theta}(GL_N)$. On adopte la notation suivante $\mathbf{\Phi}_{temp,irr}^\epsilon(GL_N)=\mathbf{\Phi}_{temp,irr}(GL_N)\cap \mathbf{\Phi}_{temp}^\epsilon(GL_N)$ pour tout $\epsilon\in\{+,-,\theta\}$. On a alors l'égalité $\mathbf{\Phi}_{temp,irr}^{\theta}(GL_N)=\mathbf{\Phi}_{temp,irr}^+(GL_N)\sqcup \mathbf{\Phi}_{temp,irr}^-(GL_N)$. On posera aussi $\mathbf{\Phi}^\epsilon_{temp}=\bigsqcup_{d\geqslant 0} \mathbf{\Phi}^\epsilon_{temp}(GL_d)$ pour tout $\epsilon\in\{+,-,\theta,\emptyset\}$. On vérifie alors aisément les propriétés suivantes

\vspace{3mm}

\begin{itemize}
\renewcommand{\labelitemi}{$\bullet$}
\item Pour tous $\varphi_1\in \mathbf{\Phi}^{\epsilon_1}_{temp}$, $\varphi_2\in \mathbf{\Phi}^{\epsilon_2}_{temp}$, $\epsilon_1,\epsilon_2\in\{\pm 1\}$, on a $\varphi_1\otimes\varphi_2\in \mathbf{\Phi}^{\epsilon_1 \epsilon_2}_{temp}$;
\item Pour tout $\varphi\in \mathbf{\Phi}^{\epsilon}_{temp}$, $\epsilon\in \{\pm 1\}$, on a $det(\varphi)\in \mathbf{\Phi}^{\epsilon^n}_{temp}$.
\end{itemize}

\vspace{3mm}

\noindent Soit $\varphi\in \mathbf{\Phi}_{temp}^\theta$. Alors, il existe une décomposition, unique au choix des indices près

$$\displaystyle \varphi\simeq \bigoplus_{j\in J} \ell_j \varphi_j \oplus \bigoplus_{i\in I} \ell_i(\varphi_i\oplus \varphi_i^\theta)$$

\noindent où

\begin{itemize}
\renewcommand{\labelitemi}{$\bullet$}
\item $I$ et $J$ sont des ensembles finis (disjoints),
\item les $\ell_i$, $\ell_j$ ($i\in I,j\in J$) sont des entiers naturels non nuls,
\item pour tout $j\in J$, $\varphi_j$ appartient à $\mathbf{\Phi}_{temp,irr}^\theta$,
\item pour tout $i\in I$, $\varphi_i$ est un élément de $\mathbf{\Phi}_{temp,irr}$ qui n'est pas conjugué-dual,
\item les $\varphi_i$, $\varphi_j$ ($i\in I,j\in J$) sont distincts deux à deux.
\end{itemize}

\noindent On notera alors $J^+$ (resp. $J^-$) le sous-ensemble des $j\in J$ tels que $\varphi_j\in \mathbf{\Phi}_{temp,irr}^+$ (resp. $\varphi_j\in \mathbf{\Phi}_{temp,irr}^-$). On a $\varphi\in \mathbf{\Phi}_{temp}^+$ (resp. $\varphi\in \mathbf{\Phi}_{temp}^-$) si et seulement si $\ell_j$ est pair pour tout $j\in J^-$ (resp. pour tout $j\in J^+$). Soit $\epsilon\in \{\pm\}$ et supposons que $\varphi\in \mathbf{\Phi}_{temp,irr}^{\epsilon}$. Fixons une forme bilinéaire non dégénérée $B$ qui soit $\epsilon$-conjuguée-duale. On a alors une identification (à automorphisme intérieur près)

$$\displaystyle Aut(\varphi,B)=\prod_{j\in J^\epsilon} O(\ell_j,\mathbb{C}) \times \prod_{j\in J^{-\epsilon}} Sp(\ell_j,\mathbb{C}) \times \prod_{i\in I} GL(\ell_i,\mathbb{C})$$

\noindent d'où l'on déduit une identification (bien définie cette fois)

$$\mbox{(3)}\;\;\; \mathcal{S}_{\varphi}=\{\pm 1\}^{J^{\epsilon}}$$

\subsection{Correspondance de Langlands pour les groupes unitaires}

Si $(V,h)$ un espace hermitien de dimension finie $d$ et $G$ est son groupe unitaire, on posera $\mathbf{\Phi}_{temp}(G)=\mathbf{\Phi}_{temp}^{(-1)^{d+1}}(GL_d)$ et on notera $\mathcal{E}^G(\varphi)$ l'ensemble des caractères $\epsilon$ de $\mathcal{S}_{\varphi}$ tels que $\epsilon(z_\varphi)=\mu(V,h)$. La correspondance locale de Langlands postule le fait suivant

\vspace{4mm}

 \textbf{(CLL)} Pour tout groupe unitaire $G=U(V,h)$, il existe une décomposition en union disjointe
 
$$\displaystyle Temp(G)=\bigsqcup_{\varphi\in \mathbf{\Phi}_{temp}(G)} \Pi^G(\varphi)$$

 et des bijections
 
$$\mathcal{E}^G(\varphi)\simeq \Pi^G(\varphi)$$

$$\epsilon\mapsto \sigma(\varphi,\epsilon)$$

\vspace{4mm}

 Les ensemble finis $\Pi^G(\varphi)$ sont appelés $L$-paquets. La correspondance de Langlands fournit donc une paramétrisation des ensembles $Temp(G)$. Cette paramétrisation doit évidemment vérifier certaines conditions. On va en imposer trois, de nature endoscopique. Avant de les énoncer, introduisons une notation. Pour $G=U(V,h)$, $\varphi\in \mathbf{\Phi}_{temp}(G)$ et $s\in \mathcal{S}_{\varphi}$, on pose
 
$$\Theta_{\varphi,s}=\sum_{\epsilon\in\mathcal{E}^G(\varphi)} \epsilon(s) \Theta_{\sigma(\varphi,\epsilon)}$$

\noindent où, rappelons le, $\Theta_{\sigma(\varphi,\epsilon)}$ désigne le caractère de la représentation $\sigma(\varphi,\epsilon)$. La première condition que l'on impose est la suivante

\vspace{4mm}

 \textbf{(Stab)} Si $(V,h)$ vérifie \textbf{(QD)}, $G=U(V,h)$ et $\varphi\in \mathbf{\Phi}_{temp}(G)$ alors $\Theta_{\varphi,0}$ est stable.

\vspace{4mm}

Nos deux autres hypothèses concernent l'une une endoscopie dite "classique" et l'autre une endoscopie tordue. Soient $(V_+,h_+)$ et $(V_-,h_-)$ deux espaces hermitiens de dimensions respectives $d_+$, $d_-$ et de groupes unitaires respectifs $G_+$ et $G_-$. On suppose dans toute la suite de ce paragraphe que $(V_+,h_+)$ et $(V_-,h_-)$ vérifient \textbf{(QD)} et on pose $d=d_++d_-$. Soient $\varphi_+\in \mathbf{\Phi}_{temp}(G_+)$, $\varphi_-\in \mathbf{\Phi}_{temp}(G_-)$ et $\mu_+$, $\mu_-$ deux caractères continus de $W_E$ que l'on identifie, via la théorie du corps de classe, à des caractères de $E^\times$.

\vspace{2mm}
 
 Supposons tout d'abord que les restrictions de $\mu_+$ et $\mu_-$ à $F^\times$ coïncident respectivement avec $sgn_{E/F}^{d_-}$ et $sgn_{E/F}^{d_+}$. Soit $(V,h)$ un espace hermitien de dimension $d$ et $G=U(V,h)$. Alors $\varphi=\mu_+\varphi_+\oplus \mu_-\varphi_-$ appartient à $\mathbf{\Phi}_{temp}(G)$. Soit $s$ l'élément qui agit comme l'identité sur (l'espace sous-jacent à) $\mu_+\varphi_+$ et comme la multiplication par $-1$ sur $\mu_-\varphi_-$. Il existe sur $\mu_+\varphi_+$ et $\mu_-\varphi_-$ des formes bilinéaires non dégénérées $(-1)^{d+1}$-conjuguées-duales. Fixons en une $B_+$ sur $\mu_+\varphi_+$ et une $B_-$ sur $\mu_-\varphi_-$. Alors $B=B_+\oplus B_-$ est une forme conjuguée-duale sur $\varphi$ de signe $(-1)^{d+1}$. On a $s\in Aut(\varphi,B)$, donc $s$ détermine un élément, aussi noté $s$, de $\mathcal{S}_{\varphi}$. Comme élément de $\mathcal{S}_{\varphi}$, $s$ ne dépend pas des choix de $B_+$ et $B_-$. Le groupe $G_+\times G_-$ est un groupe endoscopique de $G$. De plus, le couple $(\mu_+,\mu_-)$ permet de fixer, comme en 3.1, la donnée endoscopique. On normalise les facteurs de transfert comme en 3.1. La troisième condition que l'on impose est alors la suivante:
 
\vspace{4mm}

 \textbf{(TE)}  Il existe un nombre complexe $\gamma^G_{\mu_+,\mu_-}(\varphi_+,\varphi_-)$ de module $1$ tel que $\gamma^G_{\mu_+,\mu_-}(\varphi_+,\varphi_-)\Theta_{\varphi,s}$  soit le transfert de $\Theta_{\varphi_+,0}\times\Theta_{\varphi_-,0}$.
  
\vspace{4mm}

 On suppose maintenant que les restrictions de $\mu_+$ et $\mu_-$ à $F^\times$ coïncident respectivement avec $sgn_{E/F}^{d_-}$ et $sgn_{E/F}^{d_++1}$. Alors $\varphi=\mu_+\varphi_+\oplus \mu_-\varphi_-$ appartient à $\mathbf{\Phi}_{temp}^\theta(GL_d)$. Introduisons un espace $U$ de dimension $d$ sur $E$ et le groupe tordu $(M,\widetilde{M})$ associé. Posons $\pi=\pi(\varphi)$. Alors $\pi$ est une représentation conjuguée-duale de $M(F)$. On en déduit comme en 3.2 un prolongement $\tilde{\pi}$ de $\pi$ à $\widetilde{M}(F)$. Le groupe $G_+\times G_-$ est un groupe endoscopique tordu de $\widetilde{M}$. De plus, le couple $(\mu_+,\mu_-)$ permet de fixer, comme en 3.2, la donnée endoscopique. On normalise les facteurs de transfert comme en 3.2. La quatrième, et dernière, condition que l'on impose est alors la suivante:
 
\vspace{4mm}
 
 \textbf{(TET)}  Il existe un nombre complexe $c_{\mu_+,\mu_-}(\varphi_+,\varphi_-)$ de module $1$ tel que $c_{\mu_+,\mu_-}(\varphi_+,\varphi_-)\Theta_{\tilde{\pi}}$ soit le transfert de $\Theta_{\varphi_+,0}\times \Theta_{\varphi_-,0}$.
  
\vspace{4mm}

\textbf{Dans tout ce qui suit, on admet l'existence de décompositions et de bijections comme en (CLL) vérifiant les conditions (Stab), (TE) et (TET)}

\subsection{Remarques sur les conjectures}

\begin{itemize}
\renewcommand{\labelitemi}{$\bullet$}

\item Les conjectures 7.2 ont été récemment démontrées pour les groupes unitaires quasi-déployés par C.P.Mok ([Mo1], [Mo2]). Les techniques développées par Arthur dans le cas des groupes orthogonaux ([A7]) permettront probablement d'établir ces conjectures dans le cas général.

\item On a introduit dans nos conjectures des constantes indéterminées (les $\gamma_{\mu_+,\mu_-}^G(\varphi_+,\varphi_-)$ et $c_{\mu_+,\mu_-}(\varphi_+,\varphi_-)$). Ceci pour deux raisons. Tout d'abord, les conjectures n'étant pas encore établies dans le cas non quasi-déployé, on ne sait pas encore avec certitudes quelles seront les valeurs de ces constantes. Ensuite, bien que C.P.Mok détermine avec précision les constantes dans le cas quasi-déployé, nos normalisations des facteurs de transfert sont différentes des siennes. En effet, on a suivi [W1] dans le choix des normalisations c'est-à-dire qu'on a fixé les facteurs de transfert par le choix d'un épinglage (correspondant aux choix d'orbites nilpotentes régulières en 3.1 et 3.2). Dans [Mo1], les facteurs de transfert vérifient la normalisation de Whittaker (c'est-à-dire qu'ils dépendent du choix d'une donnée de Whittaker). Les facteurs de transfert ne sont donc pas les même et cela a pour conséquence que nos constantes ne sont pas les même que celles de C.P.Mok (qui valent toutes $1$ dans le cas quasi-déployé, cf les propositions 7.6.1 et 8.4.1).

\item Par indépendance linéaire des caractères, la condition \textbf{(TET)} détermine entièrement la composition des $L$-paquets $\Pi^G(\varphi)$.

\item Les bijections $\epsilon\mapsto \sigma(\varphi,\epsilon)$ ne sont elles pas uniquement déterminées par nos conjectures. En effet, soit $(V,h)$ un espace hermitien de groupe unitaire $G$ et soit $\varphi\in\mathbf{\Phi}_{temp}(G)$. Soit $\epsilon_0$ un caractère de $\mathcal{S}_\varphi$ vérifiant $\epsilon_0(z_\varphi)=1$. Si on reparamètre le $L$-paquet $\Pi^G(\varphi)$ par $\epsilon\mapsto \sigma(\varphi,\epsilon\epsilon_0)$ alors les conjectures 7.2 sont encore vérifiées à condition de multiplier les constantes $\gamma_{\mu_+,\mu_-}^G(\varphi_+,\varphi_-)$ par $\epsilon_0(s)$ (où $s$ est le même élément de $\mathcal{S}_{\varphi}$ que dans \textbf{(TE)}).

\item A nouveau par indépendance linéaire des caractères, et puisque la composition des $L$-paquets est déjà bien déterminée, il est facile de voir que la condition \textbf{(TE)} détermine les bijections $\epsilon\mapsto \sigma(\varphi,\epsilon)$ à un reparamétrage près du type précédent.

\item Soit $(V,h)$ un espace hermitien vérifiant \textbf{(QD)}, $G=U(V,h)$ et $\varphi\in\mathbf{\Phi}_{temp}(G)$. Alors, une façon de fixer la bijection $\epsilon\mapsto \sigma(\varphi,\epsilon)$ pour le $L$-paquet $\Pi^G(\varphi)$ est de choisir l'image du caractère trivial c'est-à-dire de choisir un point base dans $\Pi^G(\varphi)$. C'est ce qui sera fait en 7.7.

\item Soit $(V,h)$ un espace hermitien de dimension $d$ impaire et $G$ son groupe unitaire. Supposons que $(V,h)$ ne vérifie pas \textbf{(QD)}. Soit $(\underline{V},\underline{h})$ l'unique espace hermitien de dimension $d$ qui vérifie \textbf{(QD)}. Notons $\underline{G}$ son groupe unitaire. Pour $\alpha\in F^\times\backslash N(E^\times)$, $(V,\alpha h)$ est isomorphe à $(\underline{V},\underline{h})$. On en déduit un isomorphisme bien défini à conjugaison près entre $\underline{G}$ et $G$, donc aussi une identification entre leurs représentations tempérées irréductibles. On peut considérer $\underline{G}$ comme un groupe endoscopique de $G$ (avec $\mu_+=1$). D'après les conjectures, $\gamma^G_{1,\mu_-}(\varphi,0)\Theta^G_{\varphi,0}$ est le transfert de $\Theta^{\underline{G}}_{\varphi,0}$. Ici le transfert correspond à l'identification entre classes de conjugaison stable de $G$ et de $\underline{G}$ et les facteurs de transfert valent $1$. On en déduit que via l'isomorphisme entre $G$ et $\underline{G}$, on a l'égalité $\gamma^G_{1,\mu_-}(\varphi,0)\Theta^G_{\varphi,0}=\Theta^{\underline{G}}_{\varphi,0}$. D'après l'indépendance des caractères, on en déduit que $\gamma^G_{1,\mu_-}(\varphi,0)=1$ et $\Pi^G(\varphi)=\Pi^{\underline{G}}(\varphi)$. Les conjectures fournissent un paramètrage de cet ensemble par $\mathcal{E}^{\underline{G}}(\varphi)$. Définissons un caractère $\epsilon_0$ de $\mathcal{S}_\varphi$ donné, via l'identification 7.1(3), par

$$\displaystyle \epsilon_0((e_j)_{j\in J^1})=\prod_{j\in J^1} e_j^{d_j}$$

\noindent où pour tout $j\in J^1$, $d_j$ est la dimension de $\varphi_j$. On a alors $\epsilon_0(z_\varphi)=-1$. Pour tout $\epsilon\in \mathcal{E}^G(\varphi)$, posons $\sigma(\varphi,\epsilon)=\sigma(\varphi,\epsilon\epsilon_0)$. On vérifie que pour cette paramétrisation de $\Pi^G(\varphi)$, les conjectures 7.2 sont encore valables et que l'on a les égalités

$$\mbox{(1)}\;\;\; \gamma^G_{\mu_+,\mu_-}(\varphi_+,\varphi_-)=\gamma^{\underline{G}}_{\mu_+,\mu_-}(\varphi_+,\varphi_-)$$

\noindent Dans la suite, on supposera que la paramétrisation de $\Pi^G(\varphi)$ est ainsi obtenue.

\item Il y a quelques cas où on peut déterminer les constantes:

\begin{itemize}

\item Tout d'abord la condition \textbf{(TET)} dans le cas $d=0$ donne formellement

$$\mbox{(2)}\;\;\; c_{\mu_+,\mu_-}(0,0)=1$$

\noindent (Signalons ici que l'on adoptera la convention suivante: on notera $1$ la représentation triviale de n'importe quel groupe $G$ sauf lorsque $G$ est lui-même trivial, auquel cas on notera $0$ la représentation triviale).

\item Appliquons \textbf{(TE)} dans le cas où $(V,h)$ vérifie \textbf{(QD)}, $s=0$ et $\mu_+$ est le caractère trivial. On trouve alors que $\gamma^G_{1,\mu_-}(\varphi,0)\Theta_{\varphi,0}$ est le transfert de $\Theta_{\varphi,0}$. Puisque dans ce cas, les facteurs de transfert valent $1$, on en déduit que

$$\mbox{(3)}\;\;\; \gamma^G_{1,\mu_-}(\varphi,0)=1$$

\item Examinons ce qu'il se passe pour la condition \textbf{(TE)} lorsque l'on échange $G_+$ et $G_-$. Alors $s$ est remplacer par $sz_\varphi$ et on a $\Theta_{\varphi,sz_\varphi}=\mu(V,h)\Theta_{\varphi,s}$. D'un autre côté, les facteurs de transfert sont inchangés si le cocyle choisi pour les définir est trivial, ils sont multipliés par $-1$ sinon. Comme on a choisi le cocycle de sorte qu'il soit trivial si et seulement si $\mu(V,h)=1$, les facteurs de transfert sont aussi multipliés par $\mu(V,h)$. Par conséquent, on a

$$\mbox{(4)}\;\;\; \gamma^G_{\mu_+,\mu_-}(\varphi_+,\varphi_-)=\gamma^G_{\mu_-,\mu_+}(\varphi_-,\varphi_+)$$

\end{itemize}

\item Soit $(V,h)$ un espace hermitien de dimension $1$ et $G=U(V,h)$. On peut complétement déterminer les $L$-paquets dans ce cas. La correspondance de Langlands pour $GL_1$ est juste la théorie du corps de classe, et les représentations tempérées irréductibles de $G(F)$ sont les caractères continus de $Ker N_{E/F}$. Soit $\mu\in \Phi_{temp}(G)$. Via la théorie du corps de classe, on identifie $\mu$ à un caractère de $E^\times$ trivial sur $F^\times$. On a alors $\pi(\mu)=\mu$. Notons $\mu'$ le caractère de $Ker N_{E/F}$ défini par

$$\mu'(x/\overline{x})=\mu(x)$$

\noindent pour tout $x\in E^\times$. On vérifie que $\tilde{\mu}(\tilde{x}(\xi,\gamma))=\mu'(a)$ pour tout $\xi=(F,E,a)\in \Xi_{1,reg}$ et pour tout $\gamma\in \Gamma(\xi)$. On a $\mathcal{S}_\mu=\{\pm 1\}$. D'après \textbf{(CLL)}, $\Pi^G(\mu)$ est donc réduit à un élément. Notons $\mu_G$ cet élément. D'après \textbf{(TET)}, $c_{1,\mu_-}(\mu,0)\mu_G$ est le transfert de $\tilde{\mu}$. Les facteurs da transfert valent alors $1$ et on en déduit que $\mu_G=\mu'$. En particulier, le $L$-paquet correspondant au caractère trivial $1\in \mathbf{\Phi}_{temp}(G)$ ne contient que le caractère trivial de $G(F)$.

\end{itemize}

\subsection{Facteurs epsilons de représentations du groupe de Weil-Deligne}

Soient $d,d'\geqslant 0$ deux entiers. Pour toutes représentations admissibles irréductibles $\pi$ de $GL_d(E)$ et $\pi'$ de $GL_{d'}(E)$, on pose

$$\epsilon(\pi)=\epsilon(1/2,\pi,\psi_E^\delta)$$
$$\epsilon(\pi\times\pi')=\epsilon(1/2,\pi\times \pi',\psi_E^\delta)$$

\noindent Pour tout $\varphi\in \mathbf{\Phi}_{temp}(GL_d)$, on pose

$$\epsilon(\varphi)=\epsilon(\varphi,\psi_E^\delta)$$

\noindent où $\epsilon(\varphi,\psi_E^\delta)$ est défini comme dans [GGP] §5. D'après les résultats d'Haris-Taylor et Henniart, on a les identités

$$\mbox{(1)}\;\;\; \epsilon(\varphi)=\epsilon(\pi(\varphi))$$

$$\mbox{(2)}\;\;\; \epsilon(\varphi\otimes \varphi')=\epsilon(\pi(\varphi)\times \pi(\varphi'))$$

\noindent pour tous $\varphi,\varphi'\in \mathbf{\Phi}_{temp}$. Pour $\varphi\in \mathbf{\Phi}_{temp}^+$, la proposition 5.2(2) de [GGP] affirme l'égalité

$$\mbox{(3)}\;\;\; \epsilon(\varphi)=1$$

\noindent Puisque pour $\varphi\in \mathbf{\Phi}_{temp}^\theta$, on a $\varphi\oplus\varphi\in \mathbf{\Phi}_{temp}^+$, on en déduit que

$$\mbox{(4)}\;\;\; \epsilon(\varphi)^2=\epsilon(\varphi\oplus\varphi)=1$$

\noindent pour tout $\varphi\in \mathbf{\Phi}_{temp}^\theta$.

\subsection{Une application de la proposition 6.3.1}

Dans ce paragraphe, on se donne les objets suivants:

\vspace{3mm}

\begin{itemize}
\renewcommand{\labelitemi}{$\bullet$}

\item $\underline{d}\geqslant 0$ et $\underline{d}'\geqslant 0$ deux entiers naturels;

\item $(V_+,h_+)$, $(V_-,h_-)$, $(V'_+,h'_+)$ et $(V'_-,h'_-)$  quatre espaces hermitiens;

\item $d_+$, $d_-$, $d'_+$ et $d'_-$ les dimensions respectives de ces quatre espaces hermitiens;

\item $G_+$, $G_-$, $G'_+$ et $G'_-$ les groupes unitaires respectifs de ces quatre espaces hermitiens;

\item $\varphi_+\in \mathbf{\Phi}_{temp}(G_+)$, $\varphi_-\in \mathbf{\Phi}_{temp}(G_-)$, $\varphi'_+\in \mathbf{\Phi}_{temp}(G'_+)$ et $\varphi'_-\in \mathbf{\Phi}_{temp}(G'_-)$ des paramètres de Langlands;

\item $\mu_+$, $\mu_-$, $\mu'_+$, $\mu'_-$ des caractères continus de $E^\times$;

\end{itemize}

\vspace{3mm}

\noindent On suppose que ces données vérifient les conditions suivantes:

\vspace{3mm}

\begin{itemize}
\renewcommand{\labelitemi}{$\bullet$}

\item $\underline{d}$ est pair et $\underline{d}'$ est impair;

\item $\underline{d}=d_++d_-$ et $\underline{d}'=d'_++d'_-$;

\item $(V_+,h_+)$, $(V_-,h_-)$, $(V'_+,h'_+)$ et $(V'_-,h'_-)$ vérifient tous la condition \textbf{(QD)};

\item $\mu_{+|F^\times}=sgn_{E/F}^{d_-}$, $\mu_{-|F^\times}=sgn_{E/F}^{d_++1}$, $\mu'_{+|F^\times}=sgn_{E/F}^{d'_-}$ et $\mu'_{-|F^\times}=sgn_{E/F}^{d'_++1}$.

\end{itemize}

\vspace{3mm}

\noindent Posons 

$$\underline{\varphi}=\mu_+\varphi_+\oplus\mu_-\varphi_-\in \mathbf{\Phi}_{temp}^\theta$$

\noindent et

$$\underline{\varphi}'=\mu'_+\varphi'_+\oplus\mu'_-\varphi'_-\in\mathbf{\Phi}_{temp}^\theta$$

\noindent D'après la définition donnée en 4.3, on a l'égalité

$$\mbox{(1)}\;\;\; \epsilon_{\nu_1}(\pi(\underline{\varphi}),\pi(\underline{\varphi}'))=\omega_{\pi(\underline{\varphi})}(\nu_1)\omega_{\pi(\underline{\varphi}')}(-\nu_1)\epsilon(\pi(\underline{\varphi})\times\pi(\underline{\varphi}'))$$

\noindent où l'on rappelle que $\omega_{\pi(\underline{\varphi})}$ et $\omega_{\pi(\underline{\varphi}')}$ désignent les caractères centraux de $\pi(\underline{\varphi})$ et $\pi(\underline{\varphi}')$ respectivement. Par simple calcul de parité, on voit que $\mu_-\varphi_-\otimes\mu'_+\varphi'_+\in\mathbf{\Phi}_{temp}^+$ et $\mu_+\varphi_+\otimes\mu'_-\varphi'_-\in \mathbf{\Phi}_{temp}^+$. D'après 7.4(2) et 7.4(3), on en déduit que

\[\begin{aligned}
\mbox{(2)}\;\;\; \epsilon(\pi(\underline{\varphi})\times\pi(\underline{\varphi}')) & =\epsilon(\underline{\varphi}\otimes\underline{\varphi}') \\
 & =\epsilon(\mu_+\varphi_+\otimes\mu'_+\varphi'_+)\epsilon(\mu_-\varphi_-\otimes\mu'_-\varphi'_-)
\end{aligned}\]

\noindent D'après les résultats d'Harris-Taylor et Henniart, on a

$$\omega_{\pi(\underline{\varphi})}=det \underline{\varphi}$$

$$\omega_{\pi(\underline{\varphi}')}=det \underline{\varphi}'$$

\noindent Or, $det(\underline{\varphi})=det(\mu_+\varphi_+)det(\mu_-\varphi_-)$ et $det(\underline{\varphi}')=det(\mu'_+\varphi'_+)det(\mu'_-\varphi'_-)$. Puisque $\mu_+\varphi_+\in \mathbf{\Phi}_{temp}^{-1}(GL_{d_+})$, on a $det(\mu_+\varphi_+)\in \mathbf{\Phi}_{temp}^{(-1)^{d_+}}(GL_1)$ donc $det(\mu_+\varphi_+)_{|F^\times}=sgn_{E/F}^{d_+}$. De la même façon, on montre que

\begin{center}
$det(\mu_-\varphi_-)_{|F^\times}=det(\mu'_+\varphi'_+)_{|F^\times}=1$ et $det(\mu'_-\varphi'_-)_{|F^\times}=sgn_{E/F}^{d'_-}$.
\end{center}

\noindent On en déduit que

$$\mbox{(3)}\;\;\; (\omega_{\pi(\underline{\varphi})})_{|F^\times}=sgn_{E/F}^{d_+}$$

$$\mbox{(4)}\;\;\; (\omega_{\pi(\underline{\varphi}')})_{|F^\times}=sgn_{E/F}^{d'_-}$$

\noindent Posons

$$\displaystyle \Sigma_+=\sum_{\sigma\in \Pi^{G_+}(\varphi_+)} \sigma$$

\noindent On définit de même $\Sigma_-$, $\Sigma'_+$ et $\Sigma'_-$. La condition 7.2\textbf{(TET)}, la proposition 6.4.1, ainsi que les points (1) à (4) précédents entraînent

\[\begin{aligned}
\mbox{(5)}\;\;\; & S_{\mu_+\mu'_+}(\Sigma_+,\Sigma'_+)S_{\mu_-\mu'_-}(\Sigma_-,\Sigma'_-)= \\
 & c_{\mu_+,\mu_-}(\varphi_+,\varphi_-)c_{\mu'_+,\mu'_-}(\varphi'_+,\varphi'_-)sgn_{E/F}(-1)^{d_-} \epsilon(\mu_+\varphi_+\otimes \mu'_+\varphi'_+) \epsilon(\mu_-\varphi_-\otimes\mu'_-\varphi'_-)
\end{aligned}\]

\subsection{Détermination des constantes pour le changement de base}

\noindent On se donne les objets suivants

\vspace{3mm}

\begin{itemize}
\renewcommand{\labelitemi}{$\bullet$}

\item $(V,h)$ et $(V',h')$ deux espaces hermitiens;

\item $d$ et $d'$ leurs dimensions respectives;

\item $G$ et $G'$ leurs groupes unitaires respectifs;

\item $\varphi\in \mathbf{\Phi}_{temp}(G)$ et $\varphi'\in\mathbf{\Phi}_{temp}(G')$ des paramètres de Langlands;

\item On pose 

\begin{center}
$\displaystyle \Sigma=\sum_{\sigma\in \Pi^G(\varphi)} \sigma$ et $\displaystyle \Sigma'=\sum_{\sigma'\in \Pi^{G'}(\varphi')} \sigma'$
\end{center}

\item $\mu$ et $\mu'$ deux caractères continus de $E^\times$.

\end{itemize}

\vspace{3mm}

\noindent On suppose vérifiées les conditions suivantes:

\vspace{3mm}

\begin{itemize}
\renewcommand{\labelitemi}{$\bullet$}

\item $(V,h)$ et $(V',h')$ vérifient la condition \textbf{(QD)};

\item $\mu_{|F^\times}=sgn_{E/F}^{d'}$ et $\mu'_{|F^\times}=sgn_{E/F}^{d+1}$.

\end{itemize}

\vspace{3mm}

\begin{prop}
Sous ces hypothèses on a les égalités suivantes

\begin{enumerate}[(i)]

\item $$\displaystyle \frac{S_{\mu\mu'}(\Sigma,\Sigma')}{\epsilon(\mu\varphi\otimes\mu'\varphi')}=\left\{
    \begin{array}{ll}
        \gamma_{\psi}(N_{E/F})^{-1}sgn_{E/F}(-2) & \mbox{si } d \mbox{ et } d' \mbox{ sont impairs} \\
        1 & \mbox{sinon.}
    \end{array}
\right.
$$

\item $$\displaystyle c_{\mu,\mu'}(\varphi,\varphi')=\left\{
    \begin{array}{ll}
        \gamma_{\psi}(N_{E/F})^{-1}sgn_{E/F}(2) & \mbox{si } d \mbox{ et } d' \mbox{ sont impairs} \\
        1 & \mbox{sinon.}
    \end{array}
\right.
$$

\item 

$$c_\Sigma(1)=c_{\Sigma'}(1)=1$$

\end{enumerate}
\end{prop}

\noindent\ul{Preuve}: La preuve consiste principalement à appliquer l'égalité 7.5(5) dans divers cas particuliers. On remarquera que l'intersection des notations du paragraphe 7.5 et de la proposition est vide.

\begin{itemize}
\renewcommand{\labelitemi}{$\bullet$}

\item Supposons $d$ et $d'$ de parités différentes. On peut alors appliquer 7.5(5) au cas où $\underline{d}=0$, $(V,h)=(V'_+,h'_+)$, $(V',h')=(V'_-,h'_-)$, $\varphi'_+=\varphi$, $\varphi'_-=\varphi'$, $\mu'_+=\mu$ et $\mu'_-=\mu'$. On obtient alors l'égalité

$$\mbox{(1)}\;\;\; c_{\mu,\mu'}(\varphi,\varphi')=c_{\Sigma}(1)c_{\Sigma'}(1)$$

\noindent D'après le résultat de Rodier ([Ro]), $c_{\Sigma}(1)$ compte le nombre de représentations de $\Pi^{G}(\varphi)$ admettant un modèle de Whittaker divisé par le nombre de types de modèles de Whittaker pour $G$. Dans tout les cas c'est un demi-entier positif et on a le même résultat pour $c_{\Sigma'}(1)$. Puisque $c_{\mu,\mu'}(\varphi,\varphi')$ est un nombre complexe de module $1$, on a par conséquent

$$\mbox{(2)}\;\;\; c_{\mu,\mu'}(\varphi,\varphi')=1$$

\noindent C'est la deuxième égalité du (ii) de la proposition dans le cas où $d$ et $d'$ sont de parités différentes.

\item Supposons $d=0$ et $d'$ impair. Alors, on a trivialement $c_\Sigma(1)=1$. Les égalités (1) et (2) entraînent

$$\mbox{(3)}\;\;\; c_{\Sigma'}(1)=0$$

\item Supposons $d$ pair et $d'$ impair. Les points (1), (2) et (3) alliés entraînent

$$\mbox{(4)}\;\;\; c_\Sigma(1)=1$$

\item On déduit facilement des égalités (3) et (4) le (iii) de la proposition (quitte à opérer une permutation entre $\Sigma$ et $\Sigma'$).

\item Supposons $d$ pair et $d'$ impair. Alors, on peut appliquer 7.5(5) au cas où $(V_+,h_+)=(V'_-,h'_-)=0$, $(V_-,h_-)=(V,h)$, $\mu_-=\mu$ et $\varphi_-=\varphi$. On obtient alors

$$c_{\Sigma}(1)c_{\Sigma'_+}(1)=c_{\mu_+,\mu}(0,\varphi)c_{\mu'_+,\mu'_-}(\varphi'_+,0)$$

\noindent où, rappelons le, $0$ désigne la représentation triviale du groupe trivial. D'après le (iii) de la proposition, on a $c_{\Sigma}(1)=c_{\Sigma'_+}(1)=1$. D'après l'égalité (2), on a aussi $c_{\mu'_+,\mu'_-}(\varphi'_+,0)=1$. On en déduit que

$$\mbox{(5)}\;\;\; c_{\mu_+,\mu}(0,\varphi)=1$$

\item Supposons toujours $d$ pair et $d'$ impair. On peut alors appliquer 7.5(5) au cas où $(V_+,h_+)=(V'_+,h'_+)=0$, $(V_-,h_-)=(V,h)$, $(V'_-,h'_-)=(V',h')$, $\mu_-=\mu$, $\mu'_-=\mu'$, $\varphi_-=\varphi$ et $\varphi'_-=\varphi'$. On obtient

$$S_{\mu\mu'}(\Sigma,\Sigma')=c_{\mu_+,\mu}(0,\varphi)c_{\mu'_+,\mu'}(0,\varphi')\epsilon(\mu\varphi\otimes\mu'\varphi')$$

\noindent D'après (5), on a $c_{\mu_+,\mu}(0,\varphi)=1$ et d'après (2), on a $c_{\mu'_+,\mu'}(0,\varphi')=1$. On en déduit que

$$\mbox{(6)}\;\;\; \displaystyle \frac{S_{\mu\mu'}(\Sigma,\Sigma')}{\epsilon(\mu\varphi\otimes\mu'\varphi')}=1$$

\noindent C'est le (i) de la proposition dans le cas $d$ pair et $d'$ impair.

\item Supposons que $d$ et $d'$ sont pairs. On peut appliquer 7.5(5) au cas où $(V_+,h_+)=(V,h)$, $(V_-,h_-)=(V',h')$, $(V'_-,h'_-)=0$, $\mu_+=\mu$, $\mu_-=\mu'$, $\varphi_+=\varphi$ et $\varphi_-=\varphi'$. On obtient alors

$$c_{\Sigma'}(1)S_{\mu\mu'_+}(\Sigma,\Sigma'_+)=c_{\mu,\mu'}(\varphi,\varphi') c_{\mu'_+,\mu'_-}(\varphi'_+,0) \epsilon(\mu\varphi\otimes\mu'_+ \varphi'_+)$$

\noindent D'après le (iii) de la proposition, on a $c_{\Sigma'}(1)$. D'après (2), on a $c_{\mu'_+,\mu'_-}(\varphi'_+,0)=1$ et d'après l'égalité (6), on a aussi $S_{\mu\mu'_+}(\Sigma,\Sigma'_+)=\epsilon(\mu\varphi\otimes\mu'_+ \varphi'_+)$. On en déduit que

$$\mbox{(7)}\;\;\; c_{\mu,\mu'}(\varphi,\varphi')=1$$

\noindent C'est le (ii) de la proposition dans le cas où $d$ et $d'$ sont pairs.

\item Supposons $d$ impair et $d'$ pair. On peut appliquer 7.5(5) au cas où $(V_+,h_+)=(V',h')$, $(V'_+,h'_+)=(V,h)$, $(V_-,h_-)=(V'_-,h'_-)=0$, $\mu_+=\mu'$, $\mu'_+=\mu$, $\varphi_+=\varphi'$ et $\varphi'_+=\varphi$. On obtient alors

$$S_{\mu\mu'}(\Sigma,\Sigma')=c_{\mu',\mu_-}(\varphi',0)c_{\mu,\mu'_-}(\varphi,0)\epsilon(\mu\varphi\otimes\mu'\varphi')$$

\noindent D'après (7), on a $c_{\mu',\mu_-}(\varphi',0)=1$ et d'après (2),on a aussi $c_{\mu,\mu'_-}(\varphi,0)=1$. On en déduit que

$$\displaystyle \mbox{(8)}\;\;\; \frac{S_{\mu\mu'}(\Sigma,\Sigma')}{\epsilon(\mu\varphi\otimes\mu'\varphi')}=1$$

\noindent C'est le (i) de la proposition dans le cas où $d$ impair et $d'$ pair.

\item Supposons que $d$ est $d'$ sont pairs. On peut appliquer 7.5(5) au cas où $(V_+,h_+)=(V,h)$, $(V'_+,h'_+)=(V',h')$, $(V_-,h_-)=0$, $\mu_+=\mu$, $\mu'_+=\mu'$, $\varphi_+=\varphi$ et $\varphi'_+=\varphi'$. On obtient alors

$$S_{\mu\mu'}(\Sigma,\Sigma')c_{\Sigma'_-}(1)=c_{\mu,\mu_-}(\varphi,0)c_{\mu',\mu'_-}(\varphi',\varphi'_-)\epsilon(\mu\varphi\otimes\mu'\varphi')$$

\noindent D'après le (iii) de la proposition, on a $c_{\Sigma'_-}(1)=1$. D'après (7), on a $c_{\mu,\mu_-}(\varphi,0)=1$ et d'après (2) on a aussi $c_{\mu',\mu'_-}(\varphi',\varphi'_-)=1$. On en déduit que

$$\mbox{(9)} \displaystyle \;\;\; \frac{S_{\mu\mu'}(\Sigma,\Sigma')}{\epsilon(\mu\varphi\otimes\mu'\varphi')}=1$$

\noindent C'est le (i) de la proposition dans le cas où $d$ et $d'$ sont pairs.

\item Supposons que $d$ et $d'$ sont impairs. On peut appliquer 7.5(5) au cas où $(V_+,h_+)=(V,h)$, $d_-=d'_-=1$, $d'_+=0$, $\mu_+=\mu'_-=\mu$, $\mu_-=1$, $\varphi_+=\varphi$ et $\varphi_-=\varphi'_-=1$ (où ici $1$ désigne la représentation triviale). En 7.3, on a complétement déterminé la correspondance de Langlands en dimension $1$. On a vu en particulier que le $L$-paquet correspondant à la représentation triviale ne contient que la représentation triviale. On a donc $\Sigma_-=\Sigma'_-=1$. L'égalité 7.5(5) s'écrit donc

$$c_{\Sigma}(1)S_{\mu}(1,1)=c_{\mu,1}(\varphi,1)c_{\mu'_+,\mu}(0,1)sgn_{E/F}(-1)\epsilon(\mu)$$

\noindent D'après le (iii) de la proposition, on a $c_{\Sigma}(1)=1$. D'après (2), on a $c_{\mu'_+,\mu}(0,1)=1$. On a donc

$$\displaystyle c_{\mu,1}(\varphi,1)=\frac{sgn_{E/F}(-1)S_{\mu}(1,1)}{\epsilon(\mu)}$$

\noindent D'après le lemme A.1 de l'appendice, le membre de droite de l'égalité ci-dessus est à une constante réelle strictement positive près égal à $sgn_{E/F}(2)\gamma_{\psi}(N_{E/F})^{-1}$. Ce terme est de valeur absolue $1$ et dans nos conjectures, on a supposé qu'il en était de même de $c_{\mu,1}(\varphi,1)$. Par conséquent, on a l'égalité

$$\mbox{(10)}\;\;\; c_{\mu,1}(\varphi,1)=sgn_{E/F}(2)\gamma_{\psi}(N_{E/F})^{-1}$$

\item Supposons toujours que $d$ et $d'$ sont impairs. On peut appliquer 7.5(5) au cas où $(V_+,h_+)=(V,h)$, $(V'_+,h'_+)=(V',h')$, $d_-=1$, $d'_-=0$, $\mu_+=\mu$, $\mu'_+=\mu'$, $\mu_-=1$, $\varphi_+=\varphi$, $\varphi'_+=\varphi'$ et $\varphi_-=1$. Comme on l'a rappelé plus haut, le $L$-paquet correspondant à $\varphi_-$ ne contient alors que la représentation triviale. On obtient donc

$$S_{\mu\mu'}(\Sigma,\Sigma')=c_{\mu,1}(\varphi,1)c_{\mu',\mu'_-}(\varphi',0)sgn_{E/F}(-1)\epsilon(\mu\varphi\otimes\mu'\varphi')$$

\noindent D'après (2), on a $c_{\mu',\mu'_-}(\varphi',0)=1$ et d'après (10), on a $c_{\mu,1}(\varphi,1)=sgn_{E/F}(2)\gamma_{\psi}(N_{E/F})^{-1}$. On en déduit

$$\displaystyle \mbox{(11)}\;\;\; \frac{S_{\mu\mu'}(\Sigma,\Sigma')}{\epsilon(\mu\varphi\otimes\mu'\varphi')}=sgn_{E/F}(-2)\gamma_{\psi}(N_{E/F})^{-1}$$

\noindent C'est le (i) de la proposition dans le cas où $d$ et $d'$ sont impairs.

\item Supposons toujours que $d$ et $d'$ sont impairs. On peut appliquer 7.5(5) au cas où $(V_+,h_+)=(V,h)$, $(V_-,h_-)=(V'_+,h'_+)=(V',h')$, $d'_-=0$, $\mu_+=\mu$, $\mu_-=\mu'_+=\mu'$, $\varphi_+=\varphi$ et $\varphi_-=\varphi'_+=\varphi'$. On obtient alors

$$S_{\mu\mu'}(\Sigma,\Sigma')c_{\Sigma'}(1)=c_{\mu,\mu'}(\varphi,\varphi')c_{\mu',\mu'_-}(\varphi',0)sgn_{E/F}(-1)\epsilon(\mu\varphi\otimes\mu'\varphi')$$

\noindent D'après (2), on a $c_{\mu',\mu'_-}(\varphi',0)=1$ et d'après le (iii) de la proposition, on a $c_{\Sigma'}(1)=1$. De (11), on déduit donc que

$$c_{\mu,\mu'}(\varphi,\varphi')=sgn_{E/F}(2)\gamma_{\psi}(N_{E/F})^{-1}$$

\noindent C'est le (ii) de la proposition dans le cas où $d$ et $d'$ sont impairs. C'était le dernier cas à traiter $\blacksquare$

\end{itemize}

\vspace{5mm}

De la proposition précédente, on déduit en particulier que $c_{\mu,\mu'}(\varphi,\varphi')$ ne dépend que de $\varphi$ et $\varphi'$. On note dorénavant $c(\varphi,\varphi')$ cette constante.

\subsection{Paramétrisation des $L$-paquets dans le cas quasi-déployé et modèles de Whittaker}

\noindent Dans cette section, on fixe

\begin{itemize}
\renewcommand{\labelitemi}{$\bullet$}

\item $(V,h)$ un espace hermitien de dimension $d$ et de groupe unitaire $G$;

\item $\varphi\in \mathbf{\Phi}_{temp}(G)$ un paramètre de Langlands.
\end{itemize}

\begin{lem}
Supposons que $(V,h)$ vérifie la condition \textbf{(QD)}, alors

\begin{enumerate}[(i)]

\item Si $d$ est impair, il existe une et une seule représentation du $L$-paquet $\Pi^G(\varphi)$ admettant un modèle de Whittaker.

\item Si $d$ est pair et $\eta\in Ker\; Tr_{E/F}\backslash\{0\}$, il existe une et une seule représentation du $L$-paquet $\Pi^G(\varphi)$ admettant un modèle de Whittaker de type $\eta$.

\end{enumerate}
\end{lem}

\noindent\ul{Preuve}: Posons $\displaystyle \Sigma=\sum_{\sigma\in \Pi^G(\varphi)} \sigma$.

\begin{enumerate}[(i)]

\item D'après la proposition 7.6.1(iii), on a

$$\displaystyle c_\Sigma(1)=\sum_{\sigma\in \Pi^G(\varphi)} c_{\sigma}(1)=1$$

\noindent D'après un résultat de Rodier ([Ro]), pour tout $\sigma\in \Pi^G(\varphi)$, on a

$$c_{\sigma}(1)=\left\{
    \begin{array}{ll}
        1 & \mbox{si } \sigma \mbox{ admet un modèle de Whittaker} \\
        0 & \mbox{sinon.}
    \end{array}
\right.
$$

\noindent Le (i) du lemme s'en déduit.

\item D'après la proposition 5.6.1, appliquée à $\Sigma_+=\Sigma$, et la proposition 7.6.1(iii), on a

$$m(\Sigma,\eta)=c_{\Sigma}(1)=1$$

\noindent Or, $m(\Sigma,\eta)$ compte précisément le nombre de représentations du $L$-paquet $\Pi^G(\varphi)$ admettant un modèle de Whittaker de type $\eta$. Le (ii) du lemme s'en suit $\blacksquare$
\end{enumerate}

\vspace{4mm}

On peut maintenant fixer la paramétrisation du $L$-paquet $\Pi^G(\varphi)$, c'est-à-dire la bijection $\epsilon\in \mathcal{E}^G(\varphi)\mapsto \sigma(\varphi,\epsilon)$ de 7.2\textbf{(CLL)}, dans les cas suivants:

\vspace{3mm}

\begin{itemize}
\renewcommand{\labelitemi}{$\bullet$}

\item Si $(V,h)$ vérifie \textbf{(QD)}, alors on a remarqué en 7.3 qu'il suffisait de choisir un point base dans $\Pi^G(\varphi)$ (la représentation correspondant à $\epsilon=1$). Si $d$ est impair, on choisit l'unique représentation du $L$-paquet admettant un modèle de Whittaker. Si $d$ est pair, on choisit l'unique représentation du $L$-paquet admettant un modèle de Whittaker de type $\nu_0\delta$.

\item Si $(V,h)$ ne vérifie pas \textbf{(QD)} mais que $d$ est impair. Soit $(\underline{V},\underline{h})$ l'unique espace hermitien de dimension $d$ qui vérifie \textbf{(QD)} et $\underline{G}$ son groupe unitaire. Alors on a fixé dans le premier point la paramétrisation du $L$-paquet $\Pi^{\underline{G}}(\varphi)$. On a vu en 7.3 comment en déduire une paramétrisation de $\Pi^G(\varphi)$ et on a supposé que les paramétrisation de $\Pi^G(\varphi)$ et $\Pi^{\underline{G}}(\varphi)$ vérifiaient cette compatibilité. La paramétrisation du $L$-paquet$\Pi^G(\varphi)$ est donc maintenant aussi fixée.
\end{itemize}

\vspace{4mm}

On remarquera qu'il reste à fixer les paramétrisations dans le cas où $(V,h)$ ne vérifie pas \textbf{(QD)} et $d$ est pair. Ce sera fait en 8.3.

\section{Preuve de la conjecture de Gan-Gross-Prasad}

\subsection{Définition de deux caractères}

\noindent Donnons nous

\vspace{3mm}

\begin{itemize}
\renewcommand{\labelitemi}{$\bullet$}

\item $(V,h)$ et $(V',h')$ deux espaces hermitiens de dimensions respectives $d$, $d'$ et de groupes unitaires respectifs $G$, $G'$;

\item $\varphi\in \mathbf{\Phi}_{temp}(G)$ et $\varphi'\in \mathbf{\Phi}_{temp}(G')$ des paramètres de Langlands.

\end{itemize}

\vspace{2mm}

\noindent On suppose que

\vspace{2mm}

\begin{itemize}
\renewcommand{\labelitemi}{$\bullet$}

\item $d$ est pair et $d'$ est impair;

\item $(V,h)$ et $(V',h')$ vérifient la condition 4.1(1).
\end{itemize}

\vspace{2mm}

\noindent On va définir deux caractères

\begin{center}
$\epsilon_{\varphi,\varphi'}^G:\mathcal{S}_{\varphi}\to\{\pm 1\}$ et $\epsilon_{\varphi,\varphi'}^{G'}: \mathcal{S}_{\varphi'}\to\{\pm 1\}$
\end{center}

Soit $s\in \mathcal{S}_{\varphi}$. Choisissons une forme conjuguée-duale non dégénérée de signe $(-1)^{d+1}$ sur $\varphi$ et notons la $B$. Relevons $s$ en un élément, encore noté $s$, de $Aut(\varphi,B)$ vérifiant $s^2=1$ (ceci est toujours possible). On pose

$$\epsilon_{\varphi,\varphi'}^G(s)=\epsilon(\varphi^{s=-1}\otimes\varphi')$$

\noindent où $\varphi^{s=-1}$ désigne la sous-représentation de $\varphi$ où $s$ agit comme $-Id$. Alors, cette définition ne dépend ni du choix de $B$ ni du relèvement choisi: cela découle essentiellement de 7.4(3) (cf theorem 6.1 de [GGP]). De plus, $\epsilon^G_{\varphi,\varphi'}$ est bien un caractère de $\mathcal{S}_{\varphi}$ (cf theorem 6.1 de [GGP]). On définit de la même manière le caractère $\epsilon_{\varphi,\varphi'}^{G'}$, en échangeant les rôles de $\varphi$ et $\varphi'$.

\vspace{3mm}

\noindent Remarquons que l'on a

$$\mbox{(1)}\;\;\; \epsilon_{\varphi,\varphi'}^G(z_\varphi)=\epsilon_{\varphi,\varphi'}^{G'}(z_{\varphi'})=\epsilon(\varphi\otimes\varphi')$$

\subsection{Utilisation des résultats précédents}

\noindent Dans ce paragraphe, on se donne:

\vspace{3mm}

\begin{itemize}
\renewcommand{\labelitemi}{$\bullet$}

\item $(V,h)$ et $(V',h')$ deux espaces hermitiens de dimensions respectives $d$, $d'$ et de groupes unitaires respectifs $G$, $G'$;

\item $(V_+,h_+)$, $(V_-,h_-)$, $(V'_+,h'_+)$ et $(V'_-,h'_-)$ quatre espaces hermitiens de dimensions respectives $d_+$, $d_-$, $d'_+$, $d'_-$ et de groupes unitaires respectifs $G_+$, $G_-$, $G'_+$, $G'_-$;

\item $\varphi_+\in\mathbf{\Phi}_{temp}(G_+)$, $\varphi_-\in\mathbf{\Phi}_{temp}(G_-)$, $\varphi'_+\in\mathbf{\Phi}_{temp}(G'_+)$ et $\varphi'_-\in\mathbf{\Phi}_{temp}(G'_-)$ des paramètres de Langlands;

\item $\mu_+$, $\mu_-$, $\mu'_+$ et $\mu'_-$ des caractères continus de $E^\times$.
\end{itemize}

\vspace{3mm}

\noindent On suppose que:

\vspace{3mm}

\begin{itemize}
\renewcommand{\labelitemi}{$\bullet$}

\item $d$ est pair et $d'$ est impair;

\item $(V,h)$ et $(V',h')$ vérifient la condition 4.1(1);

\item $d=d_++d_-$ et $d'=d'_++d'_-$;

\item $(V_+,h_+)$, $(V_-,h_-)$, $(V'_+,h'_+)$ et $(V'_-,h'_-)$ vérifient tous la condition \textbf{(QD)};

\item $\mu_{+|F^\times}=sgn_{E/F}^{d_-}$, $\mu_{-|F^\times}=sgn_{E/F}^{d_+}$, $\mu'_{+|F^\times}=sgn_{E/F}^{d'_-}$ et $\mu'_{-|F^\times}=sgn_{E/F}^{d'_+}$.

\end{itemize}

\vspace{3mm}

\noindent On pose alors

$$\varphi=\mu_+\varphi_+\oplus\mu_-\varphi_-\in\mathbf{\Phi}_{temp}(G)$$

\noindent et

$$\varphi'=\mu'_+\varphi'_+\oplus\mu'_-\varphi'_-\in\mathbf{\Phi}_{temp}(G')$$

\noindent Définissons deux éléments $s\in\mathcal{S}_\varphi$ et $s'\in\mathcal{S}_{\varphi'}$ comme suit

\begin{enumerate}[(1)]

\item Fixons des formes conjuguées-duales non dégénérées de signe $(-1)^{d+1}$ $B_+$ et $B_-$ sur $\mu_+\varphi_+$ et $\mu_-\varphi_-$ respectivement. Posons $B=B_+\oplus B_-$. Alors $s\in \mathcal{S}_{\varphi}=\pi_0\left(Aut(\varphi,B)\right)$ est l'image de l'élément qui agit trivialement sur $\mu_+\varphi_+$ et comme $-Id$ sur $\mu_-\varphi_-$.

\item De même, fixons des formes conjuguées-duales non dégénérées de signe $(-1)^{d'+1}$ $B'_+$ et $B'_-$ sur $\mu'_+\varphi'_+$ et $\mu'_-\varphi'_-$ respectivement. Posons $B'=B'_+\oplus B'_-$. Alors $s'\in \mathcal{S}_{\varphi'}=\pi_0\left(Aut(\varphi',B')\right)$ est l'image de l'élément qui agit trivialement sur $\mu'_+\varphi'_+$ et comme $-Id$ sur $\mu'_-\varphi'_-$.

\end{enumerate}

\vspace{2mm}

\noindent Alors $s$ et $s'$ ne dépendent pas des choix de $B_+$, $B_-$, $B'_+$ et $B'_-$. Posons

$$\displaystyle m(\varphi,s;\varphi',s')=\sum\limits_{\substack{\epsilon\in \mathcal{E}^G(\varphi) \\ \epsilon'\in \mathcal{E}^{G'}(\varphi')}} \epsilon(s)\epsilon'(s')m(\sigma(\varphi,\epsilon),\sigma(\varphi',\epsilon'))$$

\noindent Introduisons la représentation virtuelle suivante

$$\displaystyle \Sigma_+=\sum_{\sigma_+\in \Pi^{G_+}(\varphi_+)} \sigma_+$$

\noindent On définit de même les représentations virtuelles $\Sigma_-$, $\Sigma'_+$ et $\Sigma'_-$. D'après la conjecture 7.2\textbf{(TE)} et la proposition 6.2.1, on a l'égalité

\[\begin{aligned}
\displaystyle\mbox{(3)}\;\;\; \gamma^G_{\mu_+,\mu_-}(\varphi_+, & \varphi_-) \gamma^{G'}_{\mu'_+,\mu'_-}(\varphi'_+,\varphi'_-) m(\varphi,s;\varphi',s') \\
 & =\frac{1}{2}\left(S_{\mu_+\mu'_+}(\Sigma_+,\Sigma'_+)S_{\mu_-\mu'_-}(\Sigma_-,\Sigma'_-)+\mu(G)S_{\mu_+\mu'_-}(\Sigma_+,\Sigma'_-)S_{\mu_-\mu'_+}(\Sigma_-,\Sigma'_+)\right)
\end{aligned}\]

\noindent où, rappelons le, on a posé $\mu(G)=1$ si $G$ est quasi-déployé, $-1$ sinon. Fixons un caractère continu $\mu$ de $E^\times$ dont la restriction à $F^\times$ coïncide avec $sgn_{E/F}$. D'après la proposition 7.6.1(i) et (ii), on a les égalités

\[\begin{aligned}
S_{\mu_+\mu'_+}(\Sigma_+,\Sigma'_+)S_{\mu_-\mu'_-}(\Sigma_-, & \Sigma'_-) \\
 & =c(\varphi_+,\varphi_-)sgn_{E/F}(-1)^{d_-}\epsilon(\mu_+\varphi_+\otimes \mu'_+\varphi'_+)\epsilon(\mu_-\mu\varphi_-\otimes \mu'_-\mu^{-1}\varphi'_-) \\
 & =c(\varphi_+,\varphi_-)sgn_{E/F}(-1)^{d_-}\epsilon(\mu_+\varphi_+\otimes \mu'_+\varphi'_+)\epsilon(\mu_-\varphi_-\otimes \mu'_-\varphi'_-)
\end{aligned}\]

\noindent et

\[\begin{aligned}
S_{\mu_+\mu'_-}(\Sigma_+,\Sigma'_-)S_{\mu_-,\mu'_+}(\Sigma_-, & \Sigma'_+) \\
 & =c(\varphi_+,\varphi_-)sgn_{E/F}(-1)^{d_-}\epsilon(\mu_+\varphi_+\otimes \mu'_-\varphi'_-)\epsilon(\mu_-\mu\varphi_-\otimes \mu'_+\mu^{-1}\varphi'_+) \\
 & =c(\varphi_+,\varphi_-)sgn_{E/F}(-1)^{d_-}\epsilon(\mu_+\varphi_+\otimes \mu'_-\varphi'_-)\epsilon(\mu_-\varphi_-\otimes \mu'_+\varphi'_+)
\end{aligned}\]

\noindent On en déduit que le membre de droite de (3) est égal à

\[\begin{aligned}
\mbox{(4)}\;\;\; \displaystyle \frac{c(\varphi_+,\varphi_-)sgn_{E/F}(-1)^{d_-}}{2}  & \bigg(\epsilon(\mu_+\varphi_+\otimes \mu'_+\varphi'_+)\epsilon(\mu_-\varphi_-\otimes \mu'_-\varphi'_-) \\
 & +\mu(G)\epsilon(\mu_+\varphi_+\otimes \mu'_-\varphi'_-)\epsilon(\mu_-\varphi_-\otimes\mu'_+\varphi'_+)\bigg)
\end{aligned}\]

\noindent D'après 7.4(4), la forme bilinéaire $(\varphi_1,\varphi_2)\in\big(\mathbf{\Phi}_{temp}^\theta\big)^2\mapsto \epsilon(\varphi_1\otimes \varphi_2)$ prend ses valeurs dans $\{\pm 1\}$. On peut donc réécrire (4), sous la forme

\[\begin{aligned}
\mbox{(5)}\;\;\; \displaystyle c(\varphi_+,\varphi_-)sgn_{E/F}(-1)^{d_-}\epsilon(\mu_-\varphi_-\otimes \varphi')\epsilon(\varphi\otimes \mu'_-\varphi'_-) \left(\frac{\epsilon(\varphi\otimes\varphi')+\mu(G)}{2}\right)
\end{aligned}\]

\noindent D'après les définitions 8.1, on a $\epsilon(\mu_-\varphi_-\otimes \varphi')=\epsilon_{\varphi,\varphi'}^G(s)$ et $\epsilon(\varphi\otimes \mu'_-\varphi'_-)=\epsilon_{\varphi,\varphi'}^{G'}(s')$. Par conséquent, l'égalité (3) se réécrit

\[\begin{aligned}
\mbox{(6)}\;\;\; \displaystyle & \gamma^G_{\mu_+,\mu_-}(\varphi_+,\varphi_-)\gamma^{G'}_{\mu'_+,\mu'_-}(\varphi'_+,\varphi'_-) m(\varphi,s;\varphi',s')= \\
 & c(\varphi_+,\varphi_-)sgn_{E/F}(-1)^{d_-}\epsilon_{\varphi,\varphi'}^G(s)\epsilon_{\varphi,\varphi'}^{G'}(s')\big(\frac{\epsilon(\varphi\otimes\varphi')+\mu(G)}{2}\big)
\end{aligned}\]

\subsection{Paramétrisation des $L$-paquets dans le cas non quasi-déployé}

\noindent Donnons nous

\vspace{3mm}

\begin{itemize}
\renewcommand{\labelitemi}{$\bullet$}

\item $(V,h)$ un espace hermitien de dimension $d$ et de groupe unitaire $G$;

\item $\varphi\in\mathbf{\Phi}_{temp}(G)$ un paramètre de Langlands.
\end{itemize}

\vspace{3mm}

\noindent On suppose que

\begin{itemize}
\renewcommand{\labelitemi}{$\bullet$}

\item $d$ est pair;

\item $(V,h)$ ne vérifie pas \textbf{(QD)}.
\end{itemize}

\vspace{3mm}

C'est-à-dire que l'on se place dans le seul cas qui n'a pas été traité en 7.7. On va maintenant fixer la paramétrisation du $L$-paquet $\Pi^G(\varphi)$. Il existe un vecteur $v_0\in V$ tel que $h(v_0)=\nu_0$. Fixons un tel vecteur. Notons $V'$ l'orthogonal de $v_0$ dans $V$ et $G'$ son groupe unitaire.

\begin{lem}
Il existe $\varphi'\in \mathbf{\Phi}_{temp}(G')$ tel que l'application $(\sigma,\sigma')\mapsto m(\sigma,\sigma')$ ne soit pas identiquement nulle sur $\Pi^G(\varphi)\times \Pi^{G'}(\varphi')$. De plus, pour un tel paramètre $\varphi'$, il existe alors un unique couple $(\sigma,\sigma')\in \Pi^G(\varphi)\times \Pi^{G'}(\varphi')$ tel que $m(\sigma,\sigma')=1$ et on a $\epsilon(\varphi\otimes\varphi')=-1$.
\end{lem}
 
\noindent\ul{Preuve}: Soit $\sigma\in \Pi^G(\varphi)$ et munissons $E_\sigma$ d'un produit scalaire invariant. Pour $e_1,e_2\in E_\sigma$, notons $f_{e_1,e_2}$ la fonction sur $G'(F)$ définie par $g\mapsto (e_1,\sigma(g)e_2)$. D'après [B1] lemme 12.0.5, cette fonction appartient à l'espace de Schwartz-Harish-Chandra de $G'(F)$. Pour $e_1$ et $e_2$ convenables, $f_{e_1,e_2}$ est non nulle. D'après la formule de Plancherel-Harish-Chandra, il existe $\sigma'\in Temp(G')$ telle que ${\sigma'}^\vee(f_{e_1,e_2})$ soit non nulle. Munissons aussi $E_{\sigma'}$ d'un produit scalaire invariant. Il existe alors $e'_1,e'_2\in E_{\sigma'}$ tels que

$$\displaystyle \int_{G'(F)} (e_1,\sigma(g)e_2) (\sigma'(g)e'_1,e'_2)dg \neq 0$$

\noindent D'après [B1] section 14, cela entraîne $m(\sigma,\sigma')=1$. Il ne reste plus qu'à prendre pour $\varphi'$ l'unique paramètre de Langlands tel que $\sigma'\in \Pi^{G'}(\varphi')$. Cela montre la première partie du lemme. \\

Soit $\varphi'\in \mathbf{\Phi}_{temp}(G')$ tel que l'application $(\sigma,\sigma')\mapsto m(\sigma,\sigma')$ ne soit pas identiquement nulle sur $\Pi^G(\varphi)\times \Pi^{G'}(\varphi')$. Appliquons l'égalité 8.2(6) au cas où $s=s'=0$ (c'est-à-dire le cas où $d_-=d'_-=0$). On a ici $\mu(G)=-1$. Dans le membre de gauche de 8.2(6) n'apparaissent que des termes de module $1$ sauf eventuellement $m(\varphi,0;\varphi',0)$. Ce terme compte le nombre de couples $(\sigma,\sigma')\in \Pi^G(\varphi)\times \Pi^{G'}(\varphi')$ vérifiant $m(\sigma,\sigma')=1$, il est donc non nul par hypothèse. Dans le membre de droite de 8.2(6) n'apparaissent que des termes de module $1$ sauf eventuellement $\displaystyle \frac{\epsilon(\varphi\otimes\varphi')-1}{2}$. Ce terme ne peut être nul, puisque le membre de gauche de 8.2(6) ne l'est pas. D'après 7.4(4), on a $\epsilon(\varphi\otimes\varphi')^2=1$. On en déduit que $\epsilon(\varphi\otimes\varphi')=-1$ et que le membre de droite de 8.2(6) ne contient que des termes de module $1$. Par conséquent $m(\varphi,0;\varphi',0)=1$ c'est-à-dire qu'il existe un unique couple $(\sigma,\sigma')\in \Pi^G(\varphi)\times \Pi^{G'}(\varphi')$ tel que $m(\sigma,\sigma')=1$. $\blacksquare$

\vspace{2mm}

\noindent Fixons un paramètre $\varphi'\in \mathbf{\Phi}_{temp}(G')$ qui vérifie la conclusion du lemme. Il existe donc un unique couple $(\sigma,\sigma')\in \Pi^G(\varphi)\times \Pi^{G'}(\varphi')$ vérifiant $m(\sigma,\sigma')=1$. D'après 8.1(1), on a

$$\epsilon_{\varphi,\varphi'}^G(z_\varphi)=\epsilon(\varphi\otimes\varphi')=-1$$

\noindent Donc $\epsilon_{\varphi,\varphi'}^G\in \mathcal{E}^G(\varphi)$. On peut maintenant fixer la paramétrisation du $L$-paquet $\Pi^G(\varphi)$ en associant au caractère $\epsilon_{\varphi,\varphi'}^G$ la représentation $\sigma$.

\subsection{Détermination des constantes $\gamma^G_{\mu_+,\mu_-}(\varphi_+,\varphi_-)$}

\noindent Donnons nous

\vspace{3mm}

\begin{itemize}
\renewcommand{\labelitemi}{$\bullet$}

\item $(\underline{V},\underline{h})$ un espace hermitien de dimension $\underline{d}$ et de groupe unitaire $\underline{G}$;

\item $(\underline{V}_+,\underline{h}_+)$ et $(\underline{V}_-,\underline{h}_-)$ deux espaces hermitiens de dimension $\underline{d}_+$, $\underline{d}_-$ et de groupes unitaires $\underline{G}_+$, $\underline{G}_-$;

\item $\underline{\varphi}_+\in\mathbf{\Phi}_{temp}(\underline{G}_+)$ et $\underline{\varphi}_-\in \mathbf{\Phi}_{temp}(\underline{G}_-)$ des paramètres de Langlands;

\item $\underline{\mu}_+$ et $\underline{\mu}_-$ des caractères continus de $E^\times$.
\end{itemize}

\vspace{3mm}

\noindent On suppose que

\vspace{3mm}

\begin{itemize}
\renewcommand{\labelitemi}{$\bullet$}

\item $\underline{d}=\underline{d}_++\underline{d}_-$;

\item $(\underline{V}_+,\underline{h}_+)$ et $(\underline{V}_-,\underline{h}_-)$ vérifient \textbf{(QD)};

\item $\underline{\mu}_{+|F^\times}=sgn_{E/F}^{\underline{d}_-}$ et $\underline{\mu}_{-|F^\times}=sgn_{E/F}^{\underline{d}_+}$.
\end{itemize}

\noindent C'est tout ce qu'il nous faut pour définir une constante $\gamma^{\underline{G}}_{\underline{\mu}_+,\underline{\mu}_-}(\underline{\varphi}_+,\underline{\varphi}_-)$. Puisque maintenant les paramétrisations des $L$-paquets ont été fixées (en 8.3 et 7.7), cette constante est bien définie.

\begin{prop}
\begin{enumerate}[(i)]

\item Si $(\underline{V},\underline{h})$ vérifie \textbf{(QD)}, alors on a

$$\gamma^{\underline{G}}_{\underline{\mu}_+,\underline{\mu}_-}(\underline{\varphi}_+,\underline{\varphi}_-)=\left\{
    \begin{array}{ll}
        \gamma_{\psi}(N_{E/F})^{-1}sgn_{E/F}(-2) & \mbox{si } \underline{d}_+ \mbox{ et } \underline{d}_- \mbox{ sont impairs} \\
        1 & \mbox{sinon.}
    \end{array}
\right.
$$

\item Si $(\underline{V},\underline{h})$ ne vérifie pas \textbf{(QD)}, alors on a

$$\gamma^{\underline{G}}_{\underline{\mu}_+,\underline{\mu}_-}(\underline{\varphi}_+,\underline{\varphi}_-)=\left\{
    \begin{array}{ll}
        -\gamma_{\psi}(N_{E/F})^{-1}sgn_{E/F}(-2) & \mbox{si } \underline{d}_+ \mbox{ et } \underline{d}_- \mbox{ sont impairs} \\
        -1 & \mbox{si } \underline{d}_+ \mbox{ et } \underline{d}_- \mbox{ sont pairs} \\
        1 & \mbox{sinon.}
    \end{array}
\right.
$$

\end{enumerate}
\end{prop}

\noindent\ul{Preuve}: Posons $\underline{\varphi}=\underline{\mu}_+\underline{\varphi}_+\oplus\underline{\mu}_-\underline{\varphi}_-\in\mathbf{\Phi}_{temp}(\underline{G})$. On va distinguer les cas suivants la parité de $\underline{d}$. On changera nos notations de la façon suivante: les objets ne seront plus soulignés et porteront un exposant $'$ lorsque $\underline{d}$ est impair. Ainsi pour $\underline{d}$ pair, on aura: $(V,h)$, $(V_+,h_+)$, $(V_-,h_-)$, $G$, $G_+$, $G_-$, $d$, $d_+$, $d_-$... et pour $\underline{d}$ impair, on aura: $(V',h')$, $(V'_+,h'_+)$, $(V'_-,h'_-)$, $G'$, $G'_+$, $G'_-$, $d'$, $d'_+$, $d'_-$...

\begin{enumerate}[(i)]

\item

\begin{itemize}
\renewcommand{\labelitemi}{$\bullet$}

\item Supposons $\underline{d}=d'$ impair. On définit comme en 8.2(2), un élément $s'\in\mathcal{S}_{\varphi'}$. Appliquons l'égalité 8.2(6) au cas où  $d=0$, $\mu_+=1$ et $\varphi_+=\varphi_-=0$. D'après 7.3(2) et 7.3(3), on a $c(0,0)=\gamma^G_{1,\mu_-}(0,0)=1$. D'après la définition 8.1, on a aussi $\epsilon_{0,\varphi'}^{G'}(s')=1$. Enfin, $G$ étant dans ce cas le groupe trivial, on a $\mu(G)=1$. L'identité 8.2(6) s'écrit donc

$$\mbox{(1)}\;\;\; \gamma^{G'}_{\mu'_+,\mu'_-}(\varphi'_+,\varphi'_-)m(0,0;\varphi',s')=1$$

\noindent Revenant à la définition, on a

$$\displaystyle m(0,0;\varphi',s')=\sum_{\epsilon'\in\mathcal{E}^{G'}(\varphi')} \epsilon'(s')m\left(0,\sigma(\varphi',\epsilon')\right)$$

\noindent où la multiplicité $m\left(0,\sigma(\varphi',\epsilon')\right)$ vaut $1$ si $\sigma(\varphi',\epsilon')$ admet un modèle de Whittaker, $0$ sinon. D'après les paramétrisations 7.7, on a donc $m(0,0;\varphi',s')=1$. On déduit alors de (1) que

$$\gamma^{G'}_{\mu'_+,\mu'_-}(\varphi'_+,\varphi'_-)=1$$

\item Supposons $\underline{d}=d$ pair. On définit comme en 8.2(1), un élément $s\in\mathcal{S}_{\varphi}$. Posons

$$\displaystyle m(\varphi,s,\nu_0\delta)=\sum_{\epsilon\in\mathcal{E}^G(\varphi)} \epsilon(s)m\left(\sigma(\varphi,\epsilon),\nu_0\delta\right)$$

$$\displaystyle \Sigma_+=\sum_{\sigma_+\in \Pi^{G_+}(\varphi_+)} \sigma_+$$

$$\displaystyle \Sigma_-=\sum_{\sigma_-\in \Pi^{G_-}(\varphi_-)} \sigma_-$$

\noindent D'après la condition 7.2\textbf{(TE)} et la proposition 5.6.1, on a alors

\[\begin{aligned}
\mbox{(2)}\;\;\; \displaystyle \gamma_{\mu_+,\mu_-}^G(\varphi_+ & ,\varphi_-)m(\varphi,s,\nu_0\delta) \\
 & =\left\{
    \begin{array}{ll}
        c_{\Sigma_+}(1)c_{\Sigma_-}(1) & \mbox{si } d_+ \mbox{ et } d_- \mbox{ sont pairs} \\
        sgn_{E/F}(-2)\gamma_{\psi}(N_{E/F})^{-1}c_{\Sigma_+}(1)c_{\Sigma_-}(1) & \mbox{si } d_+ \mbox{ et } d_- \mbox{ sont impairs}
    \end{array}
\right.
\end{aligned}\]

\noindent Or, d'après le choix des paramétrisations 7.7, on a $m(\varphi,s,\nu_0\delta)=1$ et d'après la proposition 7.6.1(iii), on a $c_{\Sigma_+}(1)=c_{\Sigma_-}(1)=1$. L'égalité (2) devient alors celle de l'énoncé.
\end{itemize}

\item

\begin{itemize}
\renewcommand{\labelitemi}{$\bullet$}

\item Supposons $\underline{d}=d'$ impair. Il existe un unique espace hermitien $(\underline{V}',\underline{h}')$ de dimension $d'$ vérifiant \textbf{(QD)}. Notons $\underline{G}'$ son groupe unitaire. D'après 7.3(1), on a

$$\gamma_{\mu'_+,\mu'_-}^{G'}(\varphi'_+,\varphi'_-)=\gamma_{\mu'_+,\mu'_-}^{\underline{G}'}(\varphi'_+,\varphi'_-)$$

\noindent Dans ce cas, l'égalité de l'énoncé découle donc du (i) de la proposition.

\item Supposons $\underline{d}=d$ pair. Il existe $v_0\in V$ tel que $h(v_0)=\nu_0$. Notons $V'$ l'orthogonal de $v_0$ dans $V$ et $G'$ son groupe unitaire. Pour fixer la paramétrisation du $L$-paquet $\Pi^G(\varphi)$ en 8.3, on a choisi un paramètre $\varphi'\in \mathbf{\Phi}_{temp}(G')$ vérifiant la conclusion du lemme 8.3.1. Définissons un élément $s\in\mathcal{S}_{\varphi}$ comme en 8.2(1). Appliquons l'égalité 8.2(6) au cas où $(V'_-,h'_-)=0$, $\mu'_+=1$ et $\varphi'_+=\varphi'$. On a alors $s'=0$. D'après 7.3(1) et 7.3(3), on a $\gamma_{1,\mu'_-}^{G'}(\varphi',0)=1$. De plus, puisque $(V,h)$ ne vérifie pas \textbf{(QD)}, on a $\mu(G)=-1$. Enfin, d'après le lemme 8.3.1, on a $\epsilon(\varphi\otimes\varphi')=-1$. L'identité 8.2(6) s'écrit donc

$$\mbox{(3)}\;\;\; \gamma^G_{\mu_+,\mu_-}(\varphi_+,\varphi_-)m(\varphi,s;\varphi',0)=-c(\varphi_+,\varphi_-)sgn_{E/F}(-1)^{d_-}\epsilon_{\varphi,\varphi'}^G(s)$$

\noindent Revenant à la définition, on a

$$\displaystyle m(\varphi,s;\varphi',0)=\sum_{\epsilon\in \mathcal{E}^G(\varphi), \sigma'\in \Pi^{G'}(\varphi')} \epsilon(s) m\left(\sigma(\varphi,\epsilon),\sigma'\right)$$

\noindent D'après le lemme 8.3.1, il existe un unique couple $(\sigma,\sigma')\in \Pi^G(\varphi)\times \Pi^{G'}(\varphi')$ vérifiant $m(\sigma,\sigma')=1$ et d'après le choix des paramétrisations 8.3, on a alors $\sigma=\sigma(\varphi, \epsilon_{\varphi,\varphi'}^G)$. Par conséquent, on a

$$m(\varphi,s;\varphi',0)=\epsilon_{\varphi,\varphi'}^G(s)$$

\noindent Ainsi, l'égalité (3) se réécrit

$$\gamma^G_{\mu_+,\mu_-}(\varphi_+,\varphi_-)=-c(\varphi_+,\varphi_-)sgn_{E/F}(-1)^{d_-}$$

\noindent Maintenant l'égalité de l'énoncé découle de la proposition 7.6.1(ii). $\blacksquare$
\end{itemize}

\end{enumerate}

\subsection{Le théorème}

\noindent Dans ce paragraphe, on se donne:

\vspace{3mm}

\begin{itemize}
\renewcommand{\labelitemi}{$\bullet$}

\item $(V,h)$ et $(V',h')$ deux espaces hermitiens de dimensions respectives $d$, $d'$ et de groupes unitaires respectifs $G$, $G'$;

\item $\varphi\in\mathbf{\Phi}_{temp}(G)$ et $\varphi'\in \mathbf{\Phi}_{temp}(G')$ des paramètres de Langlands.
\end{itemize}

\vspace{3mm}

\noindent De plus, on suppose que

\vspace{3mm}

\begin{itemize}
\renewcommand{\labelitemi}{$\bullet$}

\item $d$ est pair et $d'$ est impair;

\item $(V,h)$ et $(V',h')$ vérifient la condition 4.1(1).

\end{itemize}

\vspace{3mm}

\noindent Dans cette situation, on a défini en 8.1 deux caractères

\begin{center}
$\epsilon_{\varphi,\varphi'}^G:\mathcal{S}_{\varphi}\to\{\pm 1\}$ et $\epsilon_{\varphi,\varphi'}^{G'}:\mathcal{S}_{\varphi'}\to\{\pm 1\}$
\end{center}

\noindent Rappelons que l'on avait alors observer que

$$\epsilon_{\varphi,\varphi'}^G(z_{\varphi})=\epsilon_{\varphi,\varphi'}^{G'}(z_{\varphi'})=\epsilon(\varphi\otimes\varphi')$$

\begin{theo}
Admettons les conjectures 7.2 et fixons les paramétrages des $L$-paquets comme en 7.7 et 8.3. On a alors la dichotomie suivante:

\begin{enumerate}[(i)]

\item Si $\epsilon(\varphi\otimes\varphi')\neq \mu(G)$, alors $m(\sigma,\sigma')=0$ pour tout $(\sigma,\sigma')\in \Pi^G(\varphi)\times \Pi^{G'}(\varphi')$.

\item Si $\epsilon(\varphi\otimes\varphi')=\mu(G)$, alors pour tout $(\epsilon,\epsilon')\in\mathcal{E}^G(\varphi)\times \mathcal{E}^{G'}(\varphi')$, on a

$$m(\sigma(\varphi,\epsilon),\sigma(\varphi',\epsilon'))=\left\{
    \begin{array}{ll}
        1 & \mbox{si } (\epsilon,\epsilon')=(\epsilon_{\varphi,\varphi'}^G,\epsilon_{\varphi,\varphi'}^{G'}) \\
        0 & \mbox{sinon.}
    \end{array}
\right.
$$
\end{enumerate}
\end{theo}

\vspace{4mm}

\noindent \textbf{\ul{Remarque}}: C'est la conjecture 17.3 de [GGP] restreinte aux représentations tempérées dans le cas des groupes unitaires et pour les modèles de Bessel. Dans cette référence, les auteurs fixent la correspondance de Langlands p.67. On remarquera que cette normalisation correspond exactement à la paramétrisation des $L$-paquets fixée en 7.7.

\subsection{Preuve du théorème 8.5.1}

Soient $s\in\mathcal{S}_\varphi$ et $s'\in\mathcal{S}_{\varphi'}$. On peut alors fixer des données $(V_+,h_+)$, $(V_-,h_-)$, $(V'_+,h'_+)$, $(V'_-,h'_-)$, $\varphi_+$, $\varphi_-$, $\varphi'_+$, $\varphi'_-$, $\mu_+$, $\mu_-$, $\mu'_+$ et $\mu'_-$ vérifiant les conditions de la section 8.2 et de sorte que $s$ et $s'$ soient donnés par les constructions 8.2(1) et 8.2(2) respectivement. D'après les propositions 7.6.1(ii) et 8.4.1, on a

$$\gamma^{G'}_{\mu'_+,\mu'_-}(\varphi'_+,\varphi'_-)=1$$

$$\gamma^G_{\mu_+,\mu_-}(\varphi_+,\varphi_-)=\mu(G)c(\varphi_+,\varphi_-)sgn_{E/F}(-1)^{d_-}$$

\noindent L'égalité 8.2(6) se réécrit donc sous la forme

$$\mbox{(1)}\;\;\; \displaystyle m(\varphi,s;\varphi',s')=\epsilon^G_{\varphi,\varphi'}(s)\epsilon^{G'}_{\varphi,\varphi'}(s') \frac{1+\epsilon(\varphi\otimes\varphi')\mu(G)}{2}$$

\noindent Cette égalité étant valable pour tout $(s,s')\in\mathcal{S}_{\varphi}\times \mathcal{S}_{\varphi'}$. Distinguons maintenant deux cas:

\vspace{3mm}

\begin{enumerate}[(i)]

\item Si $\epsilon(\varphi\otimes \varphi')\neq \mu(G)$. D'après 7.4(4), on a alors $\epsilon(\varphi\otimes\varphi')\mu(G)=-1$. L'égalité (1) appliquée à $s=s'=0$ devient alors

$$m(\varphi,0;\varphi',0)=0$$

\noindent Or, par définition, $m(\varphi,0;\varphi',0)$ compte le nombre de couples $(\sigma,\sigma')\in \Pi^G(\varphi)\times \Pi^{G'}(\varphi')$ vérifiant $m(\sigma,\sigma')=1$. Le (i) du théorème 8.5.1 en découle aussitôt.

\item Si $\epsilon(\varphi\otimes \varphi')= \mu(G)$. D'après 7.4(4), on a alors $\epsilon(\varphi\otimes\varphi')\mu(G)=1$. L'égalité (1) appliquée à $s=s'=0$ devient alors

$$m(\varphi,0;\varphi',0)=1$$

\noindent Or, par définition, $m(\varphi,0;\varphi',0)$ compte le nombre de couples $(\sigma,\sigma')\in \Pi^G(\varphi)\times \Pi^{G'}(\varphi')$ vérifiant $m(\sigma,\sigma')=1$. Il n'existe donc qu'un seul tel couple et on le note $(\sigma,\sigma')$. Soient $\epsilon\in\mathcal{E}^G(\varphi)$ et $\epsilon'\in\mathcal{E}^{G'}(\varphi')$ les caractères tels que $\sigma=\sigma(\varphi,\epsilon)$ et $\sigma'=\sigma(\varphi',\epsilon')$. D'après la définition de $m(\varphi,s;\varphi',s')$, on a alors

$$m(\varphi,s;\varphi',s')=\epsilon(s)\epsilon'(s')$$

\noindent pour tout $(s,s')\in\mathcal{S}_{\varphi}\times \mathcal{S}_{\varphi'}$. De l'égalité (1), on déduit donc

$$\epsilon(s)\epsilon'(s')=\epsilon^G_{\varphi,\varphi'}(s)\epsilon^{G'}_{\varphi,\varphi'}(s')$$

\noindent pour tout $(s,s')\in\mathcal{S}_{\varphi}\times \mathcal{S}_{\varphi'}$, c'est-à-dire $\epsilon=\epsilon^G_{\varphi,\varphi'}$ et $\epsilon'=\epsilon^{G'}_{\varphi,\varphi'}$. C'est exactement l'assertion du (ii) du théorème 8.5.1. $\blacksquare$

\end{enumerate}

\begin{appendices}

\section{Calcul d'un facteur epsilon par une formule intégrale}

Dans cette section, on considère un caractère continu $\mu$ de $E^\times$ dont la restriction à $F^\times$ coïncide avec $sgn_{E/F}$. Rappelons qu'en 7.4, on a posé

$$\epsilon(\mu)=\epsilon(\frac{1}{2},\mu,\psi_E^\delta)$$

Appliquons la définition 6.1 au cas où $d=d'=1$ et $\Sigma=\Sigma'=1$ sont les caractères triviaux. On a alors une identification naturelle $\Xi^*(1,1)=\Xi^*_{1,reg}=Ker\; N_{E/F}$ et on obtient

$$\mbox{(1)}\;\;\; \displaystyle S_\mu(1,1)=2\lim\limits_{s\to 0^+}\int_{Ker N_{E/F}} \mu\left(\delta^{-1}(1-x)\right) |1-x|_E^{s-1/2} dx$$

\noindent où $dx$ est la mesure de Haar sur $Ker\; N_{E/F}$ de masse totale $1$.

\vspace{2mm}

Adoptons la notation suivante: pour deux nombres complexes non nul $z$ et $z'$, on notera

$$z\propto z'$$

\noindent pour signifier que $z/z'$ est un réel strictement positif.

\begin{lema}
On a

$$\epsilon(\mu)\propto sgn_{E/F}(-2)\gamma_\psi(N_{E/F}) S_{\mu}(1,1)$$

\end{lema}

\noindent\ul{Preuve}: Munissons $E$ de la mesure de Haar qui donne à $\mathcal{O}_F+\delta\mathcal{O}_F$ la mesure $1$. On note $dz$ cette mesure. Pour $\Phi\in C_c^\infty(E)$ et $s\in\mathbb{C}$, posons

$$\displaystyle\zeta(\Phi,\mu,s)=\int_E \Phi(z) \mu(z) |z|_E^{s-1} dz$$

\noindent L'intégrale est convergente pour $Re(s)>0$ et la fonction $s\mapsto \zeta(\Phi,\mu,s)$ admet un prolongement méromorphe à tout le plan complexe. On pose aussi pour tout $s\in\mathbb{C}$,

$$L(\mu,s)=\left\{
    \begin{array}{ll}
        \left(1-\mu(\varpi_E)q_E^{-s}\right)^{-1} & \mbox{si } \mu \mbox{ est non ramifié} \\
        1 & \mbox{sinon.}
    \end{array}
\right.
$$

\noindent On définit de même $\zeta(\Phi,\mu^{-1},s)$ et $L(\mu^{-1},s)$. On a alors l'équation fonctionnelle

$$\mbox{(2)}\;\;\; \displaystyle \frac{\zeta\left(\widehat{\Phi},\mu^{-1},1-s\right)}{L\left(\mu^{-1},1-s\right)}=\epsilon(s,\mu,\psi_E^\delta) \frac{\zeta\left(\Phi,\mu,s\right)}{L\left(\mu,s\right)}$$

\noindent pour tout $\Phi\in C_c^\infty(E)$, où $\widehat{\Phi}$ désigne la transformée de Fourier de $\Phi$ par rapport au caractère $\psi_E^\delta$ et à la mesure autoduale associée que l'on note $dz'$. L'hypothèse sur $\mu$ entraîne l'identité $L(\mu,.)=L(\mu^{-1},.)$. Soient $N$ et $N'$ deux entiers naturels non nuls. On définit $\Phi_{N,N'}$ comme la fonction caractéristique de l'ensemble des $z\in E$ qui s'écrivent $z=a+b\delta$, $a,b\in F$, avec $val_F(a-1)\geqslant N$ et $val_F(b)\geqslant N'$. Dans tout ce qui suit, on suppose que $N$ et $N'$ sont assez grands. On a l'égalité

$$\zeta(\Phi_{N,N'},\mu,s)=q_F^{-N-N'}$$

\noindent pour tout $s\in\mathbb{C}$. Notons $n(\psi)$ le conducteur de $\psi$, i.e. le plus petit entier tel que $\psi$ soit triviale sur $\varpi_F^{n(\psi)}\mathcal{O}_F$. Soit $C$ tel que $dz'=Cdz$. Posons $N_0=n(\psi)-val_F(2\delta^2)$. On calcule

$$\widehat{\Phi}_{N,N'}(a+b\delta)=\left\{
    \begin{array}{ll}
        Cq_F^{-N-N'}\psi(2\delta^2 b) & \mbox{si } val_F(b)\geqslant -N+N_0 \mbox{ et } val_F(a)\geqslant -N'+N_0  \\
        0 & \mbox{sinon.}
    \end{array}
\right.
$$

Appliquant (2) à $\Phi=\Phi_{N,N'}$, on en déduit que pour tout $s\in\mathbb{C}$ tel que $Re(s)<1$, on a

$$\mbox{(3)}\;\;\; \displaystyle \epsilon\left(s,\mu,\psi_E^\delta\right)\frac{L\left(\mu,1-s\right)}{L\left(\mu,s\right)}=C\int_{val_F(a)\geqslant -N'+N_0; val_F(b)\geqslant -N+N_0} \psi(2\delta^2b)\mu(a+b\delta)^{-1}|a+b\delta|_E^{-s}dadb$$

\noindent Pour $Re(s)>1/2$, le membre de droite ci-dessus converge, lorsque $N'\to\infty$, vers

$$\displaystyle\int_{val_F(b)\geqslant -N+N_0}\int_F \psi(2\delta^2b)\mu(a+b\delta)^{-1}|a+b\delta|_E^{-s}dadb$$

\noindent Puisque le membre de gauche de (3) ne dépend pas de $N'$, on a l'égalité

$$\mbox{(4)}\;\;\; \displaystyle \epsilon\left(s,\mu,\psi_E^\delta\right)\frac{L\left(\mu^{-1},1-s\right)}{L\left(\mu,s\right)}=C\int_{val_F(b)\geqslant -N+N_0}\int_F \psi(2\delta^2b)\mu(a+b\delta)^{-1}|a+b\delta|_E^{-s}dadb$$

\noindent pour $Re(s)>1/2$. Examinons le membre de droite. Après le changement de variable $a\mapsto ab$, celui-ci est égal au produit des intégrales

$$\mbox{(5)}\;\;\; \displaystyle \int_{val_F(b)\geqslant -N+N_0} \psi(2\delta^2b)sgn_{E/F}(b)|b|_F^{1-2s} db$$

\noindent et

$$\mbox{(6)}\;\;\; \displaystyle \int_F \mu(a+\delta)^{-1}|a+\delta|_E^{-s} da$$

\noindent Or, on a 

$$\displaystyle\int_{val_F(b)\geqslant -N+N_0} \psi(2\delta^2b)db=0$$

\noindent pour $N$ assez grand. Par conséquent l'expression (5) évaluée en $s=1/2$ est égale à

$$\displaystyle \int_{val_F(b)\geqslant -N+N_0, b\in N_{E/F}(E^\times)} \psi(2\delta^2b)db$$

\noindent A une constante multiplicative positive près cette intégrale est égale à

$$\displaystyle \int_{z\in E; val_F(N_{E/F}(z))\geqslant -N+N_0} \psi(2\delta^2 N_{E/F}(z)) dz$$

\noindent Par définition de la constante de Weil, cette intégrale est elle-même égale à une constante multiplicative positive près à $\gamma_{\psi}(4\delta^2N_{E/F})=sgn_{E/F}(-1)\gamma_\psi(N_{E/F})$. Examinons maintenant l'intégrale (6). Par le changement de variable $x\in Ker N_{E/F}\mapsto a=\frac{\delta+\delta x}{1-x}$ dont le jacobien est proportionnel à $|1-x|_E^{-1}$, on obtient

$$\displaystyle \int_F \mu(a+\delta)^{-1}|a+\delta|_E^{-s} da \propto sgn_{E/F}(2)|2\delta|_E^{-s}\int_{Ker N_{E/F}} \mu(\delta^{-1}(1-x))|1-x|_E^{s-1} dx$$

\noindent On en déduit, d'après (1),

$$\displaystyle \lim\limits_{s\to (1/2)^+}\int_F \mu(a+\delta)^{-1}|a+\delta|_E^{-s} da \propto sgn_{E/F}(2)S_\mu(1,1)$$

\noindent Le résultat du lemme est alors conséquence de ce qui précède et de l'égalité (4) évaluée en $s=1/2$ $\blacksquare$

\end{appendices}

\bigskip

{\bf Bibliographie}

\bigskip

[A] J. Arthur: {\it The endoscopic classification of representations : orthogonal and
symplectic groups}, Colloquium Publication Series, AMS. To appear
  
[AGRS] A. Aizenbud, D. Gourevitch, S. Rallis, G. Schiffmann: {\it Multiplicity one theorems}, Ann. of Math. (2) 172 (2010),  no. 2, 1407-1434.

[B1] R. Beuzart-Plessis: {\it La conjecture locale de Gross-Prasad pour les représentations tempérées de groupes unitaires}, prépublication 2012

[B2] ---------------------: {\it Expression d'un facteur epsilon de paire par une formule intégrale}, prépublication 2012

[C] L.Clozel: {\it Characters of nonconnected, reductive p-adic groups}, Canad. J. Math. 39 (1987), no. 1, 149-167

[GGP]  W. T. Gan, B. Gross, D. Prasad:{\it Symplectic local root numbers, central critical $L$-values and restriction problems in the representation theory of classical groups}, astérisque 346 (2012)

[HCDS]  Harish-Chandra: {\it Admissible invariant distributions on reductive $p$-adic groups}, notes par S. DeBacker et P. Sally, University Lecture series 16, AMS (1999)

[HCvD]  Harish-Chandra: {\it Harmonic analysis on reductive $p$-adic groups}, notes par
G. van Dijk, Springer Lecture Notes 162 (1970)

[JPSS] H. Jacquet, I.I. Piatetskii-Shapiro, J. Shalika: {\it Rankin-Selberg convolutions},  Amer. J. Math.  105  (1983),  no. 2, 367-464

[K] T. Konno: {\it Twisted endoscopy and the generic packet conjecture},  Israel J. Math. 129 (2002), 253-289

[KS] R. Kottwitz, D. Shelstad: {\it Foundations of twisted endoscopy}, Astérisque No 255 (1999)

[LS] R. Langlands, D. Shelstad: {\it On the definition of transfer factors},  Math. Ann.  278  (1987),  no. 1-4, 219-271

[Mo1]: C.P. Mok: {\it Endoscopic classification of representations of quasi-split unitary groups I}, prépublication 2012

[Mo2]: C.P. Mok: {\it Endoscopic classification of representations of quasi-split unitary groups II}, prépublication 2012

[Ro] F. Rodier: {\it Modèle de Whittaker et caractères de représentations}, in Non commutative harmonic analysis, J. Carmona, J. Dixmier, M. Vergne éd. Springer LN 466 (1981), p.151-171

[Sh] D. Shelstad: {\it A formula for regular unipotent germs}, dans "Orbites unipotentes et représentations, II", Astérisque  No. 171-172 (1989), 275-277

[W1] J.-L. Waldspurger: {\it Les facteurs de transfert pour les groupes classiques: un formulaire},  Manuscripta Math.  133  (2010),  no. 1-2, 41-82

[W2] -------------------: {\it Intégrales orbitales nilpotentes et endoscopie pour les groupes classiques non ramifiés}, Astérisque No. 269  (2001)

[W3] -------------------: {\it Le lemme fondamental implique le transfert},  Compositio Math.  105  (1997),  no. 2, 153-236

[W4] -------------------: {\it La conjecture locale de Gross-Prasad pour les représentations tempérées des groupes spéciaux orthogonaux}, Asterique 347 (2012)

\bigskip

Institut de mathématiques de Jussieu 2 place Jussieu 75005 Paris \\
 e-mail: rbeuzart@gmail.com

\end{document}